\documentclass[11pt,a4paper,reqno]{amsart}
\usepackage{amsfonts,amsmath,amscd,amssymb,amsbsy,amsthm,amstext,amsopn,amsxtra,mathtools}
\usepackage{color,fullpage,mathrsfs,verbatim,subfigure}
\usepackage{stmaryrd}
\usepackage{extarrows}
\usepackage[all]{xy}
\usepackage{bm}
\usepackage{hyperref}
\usepackage{tikz}
\usepackage{tikz-cd}
\usetikzlibrary{shapes.geometric,calc,matrix,arrows,decorations.pathmorphing,arrows.meta}
\usepackage{etex}
\usepackage{cleveref}
\usepackage{nameref}
\usepackage[enableskew, vcentermath]{youngtab}
\usepackage{caption}
\usepackage{enumitem}

\hypersetup{colorlinks,linkcolor=[rgb]{0,0,1},citecolor=[rgb]{1,0,0},urlcolor=black}
 
\numberwithin{equation}{section}
\newtheorem{theorem}{Theorem}[section]
\newtheorem{lemma}[theorem]{Lemma}
\newtheorem{proposition}[theorem]{Proposition}
\newtheorem{corollary}[theorem]{Corollary}
\newtheorem{conjecture}[theorem]{Conjecture}

\newtheorem*{thm:surj_onto_homs}{Theorem \ref*{prop:surj_onto_homs}}
\newtheorem*{cor:tilting_quiver_rank2}{Corollary \ref*{thm:tilting_quiver_rank2}}
\newtheorem*{thm:complete_relns_rank2}{Theorem \ref*{thm:complete_relns_rank2}}
\newtheorem*{thm:grn2_isomorphism}{Theorem \ref*{thm:grn2_isomorphism}}

\theoremstyle{definition}
\newtheorem{definition}[theorem]{Definition}
\newtheorem{example}[theorem]{Example}

\newtheorem{remark}[theorem]{Remark}

\newtheorem{notation}[theorem]{Notation}

\newcommand{\kk}{\ensuremath{\Bbbk}}

\newcommand{\NN}{\ensuremath{\mathbb{N}}} 
\newcommand{\PP}{\ensuremath{\mathbb{P}}}
\newcommand{\QQ}{\ensuremath{\mathbb{Q}}}

\renewcommand{\SS}{\ensuremath{\mathbb{S}}}

\newcommand{\vv}{\ensuremath{\bm{\mathrm v}}}
\newcommand{\ZZ}{\ensuremath{\mathbb{Z}}}
\newcommand{\cM}{\mathcal{M}}
\newcommand{\cO}{\mathcal{O}}

\newcommand{\cW}{\mathcal{W}}
\newcommand{\cB}{\mathcal{B}}

\newcommand{\one}{\ensuremath{(\mathrm{i})}}
\newcommand{\two}{\ensuremath{(\mathrm{ii})}}
\newcommand{\three}{\ensuremath{(\mathrm{iii})}}
\newcommand{\four}{\ensuremath{(\mathrm{iv})}}

\DeclareMathOperator{\Coh}{Coh}
\DeclareMathOperator{\col}{col}

\DeclareMathOperator{\End}{End}

\DeclareMathOperator{\GL}{GL} 
\DeclareMathOperator{\Gr}{Gr} 
\DeclareMathOperator{\Hom}{Hom}
\DeclareMathOperator{\Homy}{Hom_{\cO_Y}}

\DeclareMathOperator{\hd}{h}

\DeclareMathOperator{\im}{im}

\DeclareMathOperator{\rank}{rank}

\DeclareMathOperator{\Sym}{Sym}

\DeclareMathOperator{\sign}{sgn}
\DeclareMathOperator{\tl}{t}
\DeclareMathOperator{\Young}{Young}

\newcommand{\xrightarrowdbl}[2][]{%
  \xrightarrow[#1]{#2}\mathrel{\mkern-14mu}\rightarrow
}

\def \lrn{c_{\lambda,\gamma}^{\mu}}

\def \moduli{\cM(A,\vv,\theta)} 
\def \modulie{\cM(E)} 
\def \lambdat{\widetilde{\lambda}} 
\def \modlam{\lvert \lambda \rvert} 
\def \modlamt{\lvert \lambdat \rvert} 
\def \modmu{\lvert \mu \rvert} 
\def \modgamma{\lvert \gamma \rvert} 
\def \bw{\bigwedge}
\def \TQ{Q} 
\def \ww{\wedge}

\def \uve{\underline{\xi_\sigma}}

\def \ione{i_1}
\def \itwo{i_2}
\def \ithree{i_3}

\begin{document}

\title{Reconstructing the Grassmannian of lines from Kapranov's tilting bundle}
 
\author{James Green} 
\address{Department of Mathematical Sciences, 
University of Bath, 
Claverton Down, 
Bath BA2 7AY, 
United Kingdom.}
\email{james.j.green@bath.edu}

\begin{abstract}
Let $E$ be the tilting bundle on the Grassmannian $\Gr(n,r)$ of $r$-dimensional quotients of $\kk^n$ constructed by Kapranov \cite{Kap84}. Buchweitz, Leuschke and Van den Bergh \cite{BLV16} introduced a quiver $\TQ$ and a surjective $\kk$-algebra homomorphism $\Phi\colon\kk\TQ\rightarrow A=\End(E)$, together with a recipe on how the kernel may be computed. In this paper, for the case $r=2$ we give a new, direct proof that $\Phi$ is surjective and then complete the picture by calculating the ideal of relations explicitly. As an application, we then use this presentation to show that $\Gr(n,2)$ is isomorphic to a fine moduli space of certain stable $A$-modules, just as $\PP^n$ can be recovered from the endomorphism algebra of Beilinson's tilting bundle $\bigoplus_{0\leq i\leq n}\cO_{\PP^n}(i)$.
\end{abstract}

\maketitle

\section{Introduction}
Let $\Bbbk$ be an algebraically closed field of characteristic zero, let $n$ be a positive integer and set $V=\kk^n$. For $1<r\leq n$, denote by $\Gr(n,r)$ the Grassmannian of $r$-dimensional quotients of $V$ and let $\cW$ be the rank $r$ tautological quotient bundle on $\Gr(n,r)$. Denote by $\Young(n-r,r)\subset\ZZ^r$ the set of partitions $\lambda=(n-r\geq \lambda_1\geq \dots \geq \lambda_r\geq 0)$ and write $\SS^\lambda\cW$ for the vector bundle whose fibre over each point is the irreducible $\GL(r)$-module of highest weight $\lambda$. Kapranov \cite{Kap84} proved that the vector bundle
\begin{align}
\label{eqn:tilt_bundle_grass}
E:=\bigoplus_{\lambda\in \Young(n-r,r)} \SS^\lambda \cW
\end{align}
is a \emph{tilting bundle} on $\Gr(n,r)$, i.e.\ the bounded derived category of coherent sheaves on $\Gr(n,r)$ is equivalent to bounded derived category of finite dimensional $A$-modules where $A:=\End(E)$. We refer to $A$ as \emph{Kapranov's tilting algebra}.

Recently, Buchweitz, Leuschke and Van den Bergh \cite{BLV16} used a result of Segal \cite{Seg08} to give a presentation of the tilting algebra $A$ using a quiver $\TQ$ with relations, which we call the \emph{tilting quiver}. In other words, they wrote down a surjective $\kk$-algebra homomorphism
\[
\Phi\colon\kk \TQ\xrightarrowdbl{\phantom{aaa}} A,
\]
where $\kk\TQ$ is the path algebra of $\TQ$, and provided a recipe on how the relations may be calculated from the kernels of certain linear maps. In the first half of this paper, for the specific case that $r=2$ we provide a new, direct proof of this result which bypasses Segal's result, and moreover explicitly write down the kernel of $\Phi$ as an ideal of relations in $\kk\TQ$.

First, we show for all $1<r\leq n$ that for two Young diagrams $\lambda$ and $\mu$ with $\lambda$ contained in $\mu$ (write $\lambda\leq\mu$) we have $\Hom(\SS^\lambda \cW,\SS^\mu \cW)\cong\SS^{\mu/\lambda}V$, where $\SS^{\mu/\lambda}V$ is the skew-Schur power of $V$ corresponding to the skew diagram $\mu/\lambda$. Writing $e_1,\dots,e_r$ for the standard basis of $\ZZ^r$, an immediate corollary is that for all $\lambda,\mu\in\Young(n-r,r)\subset \ZZ^r$ with $\mu=\lambda+e_i$ for some $1 \leq i \leq r$, we have $\Hom(\SS^\lambda\cW,\SS^\mu\cW)\cong V$.

Henceforth fix $r=2$ and a choice of basis $u_1,\dots,u_n$ of $V$. Then we get a basis of each of the homomorphism spaces $\Hom(\SS^\lambda\cW,\SS^{\lambda+e_i}\cW)$ which we denote by $f^\lambda_{u_1},\dots,f^\lambda_{u_n}$ if $i=1$ or $g^\lambda_{u_1},\dots,g^\lambda_{u_n}$ if $i=2$; we call these \emph{$f$-type} and \emph{$g$-type} maps. We then write down these maps explicitly in \Cref{prop:schurmaps} and use them to prove the following key surjectivity result. Given $\lambda<\mu\in\Young(n-2,2)$, define $m_1=\mu_1-\lambda_1$ and $m_2=\mu_2-\lambda_2$, and for $0\leq k\leq m_1+m_2$ let $\lambda=\tau_0<\tau_1<\dots<\tau_{m_1+m_2}=\mu$ be the unique sequence of partitions including $(\mu_1,\lambda_2)$.

\begin{thm:surj_onto_homs}
Let $r=2$ and let $\lambda<\mu\in\Young(n-2,2)$. Let $\tau_k$ be the sequence of partitions defined above. Then the composition map
\[
\Theta_{\lambda,\mu}\colon\bigotimes_{k=1}^{m_1+m_2}\Hom(\SS^{\tau_{k-1}}\cW, \SS^{\tau_{k}}\cW) \longrightarrow \Hom(\SS^\lambda \cW, \SS^\mu\cW)
\]
is surjective.
\end{thm:surj_onto_homs}

This implies that any element of $A$ may be decomposed as a linear combination of $f$-type and $g$-type maps. In light of this, we define the tilting quiver $\TQ$ to have vertex set corresponding to the irreducible summands of $E$, namely $\SS^\lambda\cW$ for each $\lambda\in\Young(n-2,2)$, and arrow set corresponding to the collection of $f$-type or $g$-type maps between each summand. This produces a staircase-like diagram for $\TQ$, see \Cref{fig:grn2tiltQ} (page~\pageref{fig:grn2tiltQ}). Then we can define a $\kk$-algebra homomorphism
\[
\Phi\colon\kk\TQ\longrightarrow A
\]
by mapping each arrow in the quiver to the appropriate $f$-type or $g$-type map. In particular, we obtain the following.

\begin{cor:tilting_quiver_rank2}
Let $r=2$. Then the $\kk$-algebra homomorphism $\Phi\colon\kk\TQ\rightarrow A$ is surjective.
\end{cor:tilting_quiver_rank2}

In the next section we describe the ideal of relations explicitly. By observing that $\ker(\Phi)$ is generated by paths with the same head and tail and length at least two, it is enough to find a basis $K_{\lambda,\mu}$ for the kernels of the induced maps $\Phi_{\lambda,\mu} \colon e_\mu\kk\TQ e_\lambda \twoheadrightarrow \Hom(\SS^\lambda\cW,\SS^\mu \cW)$ 
where $e_\mu,e_\lambda$ are idempotents corresponding to the length zero paths at the vertices $\mu,\lambda\in\TQ_0$ respectively. Define $P:=\{(\lambda,\mu)\in{\TQ_0}^2 \mid \ \lambda<\mu, \modmu\geq\modlam+2\}$. Then we have
\[
\ker(\Phi)=\left(\bigcup_{(\lambda,\mu)\in P} K_{\lambda,\mu}\right).
\]
First we consider the subset $P_2:=\{(\lambda,\mu)\in{\TQ_0}^2 \mid \ \lambda<\mu, \modmu=\modlam+2\}\subset P$ and write down each set $K_{\lambda,\mu}$ for $(\lambda,\mu)\in P_2$ explicitly. Then we define the ideal 
\[
I\colon=\left(\bigcup_{(\lambda,\mu)\in P_2} K_{\lambda,\mu}\right)\subseteq\ker(\Phi),
\]
and by considering the remaining pairs $(\lambda,\mu)\in P\setminus P_2$ we prove that $I=\ker(\Phi)$. Denote the arrows of $\TQ$ by their images under $\Phi$ (i.e.\ the $f$-type and $g$-type maps). Since we record the vertex at the tail of an arrow in the superscript but the $f$-type or $g$-type notation implies whether the next arrow is horizontal or vertical, we may omit all other superscripts without ambiguity; see \Cref{not:superscript_drop}\one. Then the presentation of Kapranov's tilting algebra $A$ is as follows.

\begin{thm:complete_relns_rank2}
For $r=2$, let $E$ be the tilting bundle  \eqref{eqn:tilt_bundle_grass} and let $A=\End(E)$. Let $\TQ$ be the quiver defined in \Cref{defn:tilting_quiver_gr_n2}. Then the $\kk$-algebra $A$ is isomorphic to $\kk\TQ/I$ for the ideal
\[
I=\left(\bigcup_{(\lambda,\mu)\in P_2} K_{\lambda,\mu}\right),
\]
where
\begin{itemize}
	\item[\one] if $\mu=(\lambda_1+2,\lambda_2)$,   $K_{\lambda,\mu}=\left\{ f_{u_j} f^\lambda_{u_i}-f_{u_i} f^\lambda_{u_j} \mid   1\leq i,j\leq n \right\}$.
	\item[\two] if $\mu=(\lambda_1,\lambda_2+2)$,  $K_{\lambda,\mu}=\left\{ g_{u_j} g^\lambda_{u_i}-g_{u_i} g^\lambda_{u_j} \mid 1\leq i,j\leq n \right\}$.
	\item[\three] if $\lambda_1=\lambda_2$ and $\mu=(\lambda_1+1,\lambda_1+1)$,  $K_{\lambda,\mu}=\left\{ g_{u_j} f^\lambda_{u_i}+g_{u_i} f^\lambda_{u_j} \mid 1\leq i,j\leq n \right\}$.
	\item[\four] if $\lambda_1>\lambda_2$ and  $\mu=(\lambda_1+1,\lambda_2+1)$,  
	\[
	K_{\lambda,\mu}=\left\{ \left( \lambda_1-\lambda_2 \right)g_{u_j} f^\lambda_{u_i}-\left( \lambda_1-\lambda_2+1\right)f_{u_i} g^\lambda_{u_j} + f_{u_j} g^\lambda_{u_i} \mid 1\leq i,j\leq n \right\}.
	\] 
\end{itemize}
\end{thm:complete_relns_rank2}

We illustrate this result by computing the relations for $\Gr(5,2)$.

\smallskip
In the second half of the paper we reconstruct $\Gr(n,2)$ from Kapranov's tilting algebra $A$ using our presentation from \Cref{thm:complete_relns_rank2}. Let $E_0,\dots,E_n$ denote the summands of the tilting bundle $E$ with $E_0=\cO_{\Gr(n,2)}$ and let $\vv=(v_i)$ be the dimension vector given by $v_i=\rank(E_i)$. Since $A$ is a finite dimensional associative $\kk$-algebra and $\vv$ is indivisible, King constructs for any generic $\theta$ the fine moduli space of $\theta$-stable $A$-modules with dimension $\vv$, denoted $\moduli$, as a GIT quotient; see \cite{Kin94}. Now, since $E_0,\dots,E_n$ are globally generated, Craw, Ito and Karmazyn~\cite{CIK17} choose a particular stability parameter $\theta$ which is generic and  construct the universal morphism $f_E:\Gr(n,2)\rightarrow \moduli$. This generalises the classical morphism from a scheme with a basepoint-free line bundle to its linear series. Using \cite[Theorem~2.6,~Remark~2.8]{CIK17} we deduce that the morphism
\begin{equation}
\label{eqn:closed_imm_intro}
f_E\colon \Gr(n,2)\longrightarrow \modulie:= \moduli
\end{equation}
to the \emph{multigraded linear series} $\modulie$ is a closed immersion; see \Cref{thm:closed_immersion}. Hence we may embed $\Gr(n,2)$ into an auxiliary, ambient moduli space $\modulie$; in fact, we have the following.

\begin{thm:grn2_isomorphism}
The closed immersion $f_E:\Gr(n,2)\rightarrow \modulie$ is an isomorphism.
\end{thm:grn2_isomorphism}

The proof follows a similar strategy to that used by Craw and the author in \cite{CG18} to prove an analogous result for toric quiver flag varieties. First we describe the image of the closed immersion $f_E$ by observing that, as a quiver flag variety, points of $\modulie$ may be described by two systems of matrices: one which must satisfy the stability condition of \cite[Lemma~2.1]{Cra11}, and one which must satisfy the relations from \Cref{thm:complete_relns_rank2}. The key remaining point is to prove that  $f_E$ is surjective. To do this, we consider an arbitrary closed point $w\in\modulie$ and by working through the tilting quiver vertex by vertex, use the quiver relations and stability  conditions to show that $w$ is isomorphic to a closed point of the form $f_E(y)$ for some $y\in \Gr(n,2)$. We do this by induction with $\Gr(4,2)$ forming the base case and then complete the proof using \cite[Theorem~1.13]{BP08}.

\smallskip
By generalising the above construction to $\Gr(n,r)$, \Cref{thm:grn2_isomorphism} combined with the work of Bergman-Proudfoot~\cite{BP08} suggests the following:

\begin{conjecture}
\label{conj:general_gr_iso}
For any $1\leq r<n$ the morphism $f_E:\Gr(n,r)\rightarrow \modulie$ is an isomorphism.
\end{conjecture}

Due to \cite[Theorem~6.9]{BLV16}, we already have the tilting quiver and a surjective homomorphism $\Phi\colon \kk\TQ\twoheadrightarrow A$ for $\Gr(n,r)$ in general. To prove the above conjecture using the methods in this paper it remains to complete two main steps: firstly, write down the ideal of relations $\ker(\Phi)$ explicitly, and secondly, take a similar approach to the proof in \Cref{chap:moduli_grn2} to get the result. Both steps pose a far greater combinatorial challenge than the $r=2$ case. We must write down a system of maps that generalise the $f$-type and $g$-type maps above; such a system is called a \emph{compatible Pieri system} for the tilting quiver. An example of a map between Schur powers including a partition with three rows can be found in \cite[Section~7.1]{Gre18}. To complete this task, Buchweitz, Leuschke and Van den Bergh mention that Olver was the first to write down such a system in the preprint \cite{Olv82}; see also \cite{ABW82}, \cite{MO92}, \cite{SW11}, \cite{Sam09}. Only after a Pieri system can be written down in a practical way can $\ker(\Phi)$ be calculated. Consequently, little can be said about a potential proof of \Cref{thm:grn2_isomorphism} in the general case. We do however suspect that such a proof, while combinatorially challenging, would be quite similar to the $\Gr(n,2)$ case.

Despite complications with the general Grassmannian case, evidence that the closed immersion $f_E\colon Y \rightarrow \modulie$ is an isomorphism when $Y=\Gr(n,2)$ or $Y$ is any toric quiver flag variety (see \cite{CG18}), along with the work of Bergman and Proudfoot \cite{BP08} which identifies any quiver flag variety $Y$ with a connected component of $\modulie$, leads us to the following:

\begin{conjecture}
\label{conj:all_qfvs}
Let $Y$ be any quiver flag variety and $E$ the tilting bundle from \cite[Theorem~4.5]{Cra11}. Then the morphism $f_E:Y\rightarrow \modulie$ is an isomorphism.
\end{conjecture}

\subsection*{Acknowledgements} Many thanks to Alastair Craw, whose constant support and inspiration meant this project was possible. Thanks to Ed Segal for examining the author's PhD Thesis \cite{Gre18}, which contains this paper's results, and also for giving some useful insights with the final proof. Thanks to David Calderbank for many helpful conversations, and to Jerzy Weyman and Victor Reiner for some helpful emails. The author was supported by a doctoral studentship from the EPSRC.

\section{The tilting quiver of $\Gr(n,2)$}
\label{chap:tiltingalgebra}

\subsection{Schur powers}
\label{subsec:schurpowers}
We begin by establishing our conventions for Young diagrams,  Littlewood-Richardson numbers and Schur powers. Let $\lambda \in \ZZ^r$ be a weakly decreasing finite sequence of non-negative integers $\lambda_1\geq \dots \geq \lambda_r\geq0$. We call such $\lambda$ a \textit{partition} with $r$ \textit{parts}, even if $\lambda$ ends in a trail of zeroes. For every partition we associate a \textit{Young diagram}, a finite collection of boxes arranged into left-justified rows in descending order. For example, the Young diagram representing $\lambda=(4,3,1,1)$ is
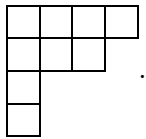
\[
\yng(4,3,1,1).
\]
Denote the number of boxes in a Young diagram by $\modlam:=\sum_{i=1}^r \lambda_i$, and let $\Young(n,r)$ be the set of Young diagrams with at most $n$ columns and $r$ rows, i.e.
\[
\Young(n,r):=\left\{\lambda=(\lambda_1,\dots,\lambda_r)\in \ZZ^r \mid n\geq \lambda_1 \geq \dots \geq \lambda_r \geq 0 \right\}.
\]
Note that we will often refer to partitions and Young diagrams as one and the same.

Let $\lambda,\mu$ be partitions with $r$ parts and suppose $\mu$ \textit{contains} $\lambda$, i.e.\ $\lambda_i\leq\mu_i$ for all $i=1,\dots,r$, and write $\lambda\leq \mu$ (or $\lambda<\mu$ for strict containment). The \textit{skew diagram} $\mu / \lambda$ is given by the Young diagram $\mu$ with $\lambda$ removed from the top left corner. For example, if $\lambda=(3,2,2,1)$ and $\mu=(6,4,4,2)$ then $\mu / \lambda$ is as follows:

$$\begin{array}{cccc}
\phantom{aa} \yng(6,4,4,2) \phantom{aa} & \phantom{aa} \yng(3,2,2,1) \phantom{a} & \young(:::~~~,::~~,::~~,:~) \\
\mu & \lambda & \mu / \lambda \\
\end{array}$$
Observe that any Young diagram $\mu$ is also a skew diagram $\mu / \lambda$ where $\lambda=(0)$.

A \textit{filling} of a skew diagram is the insertion of a positive integer into each box. A \textit{(semi-standard) skew tableau} is a skew diagram with a filling such that each row is weakly increasing and each column is strictly increasing. We say that a skew tableau $\mu / \lambda$ has \textit{content} $\gamma=(\gamma_1,\dots,\gamma_k)\in \NN^k$ if $\mu / \lambda$ contains $\gamma_1$ $1$'s, $\gamma_2$ $2$'s, and so on up to $\gamma_k$ $k$'s. For example, taking $\lambda$ and $\mu$ as above, one possible skew tableau with shape $\mu / \lambda$ and content $\gamma=(4,2,2)$ is 
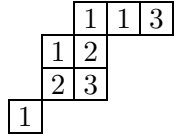
\[
\young(:::113,::12,::23,:1) \ .
\]

The sequence of integers given by concatenating the rows of a skew tableau from top to bottom and in reverse order is called the \textit{reverse word} of the tableau, e.g. the reverse word of the tableau above is $3,1,1,2,1,3,2,1$. We say that a word is a \textit{lattice word} if the content of every initial sequence is a partition, i.e.\ for every initial sequence of the word there must contain at least as many $1$'s as $2$'s, at least as many $2$'s as $3$'s, and so on. The reverse word above is not a lattice word, for example, as the first number is a $3$ so there are more $3$'s than $1$'s in the first letter of the reverse word. However, if we swap the $3$ in the top right box with the $1$ in the bottom left box, we now have a skew tableau whose reverse word  $1,1,1,2,1,3,2,3$ is a lattice word.

\begin{definition}
\label{def:lrt}
A \textit{Littlewood-Richardson tableau} is a skew tableau whose reverse word is a lattice word.
\end{definition}

\begin{definition}
\label{def:lrn}
Let $\lambda,\mu,\gamma$ be partitions such that $\lambda\leq \mu$ and $\modlam+\modgamma=\modmu$. The \textit{Littlewood-Richardson number} $\lrn\in \NN$ is equal to the number of Littlewood-Richardson tableaux of shape $\mu/\lambda$ and content $\gamma$.
\end{definition}

\begin{remark}
\label{rem:basiclrn}
\one \ If $\lambda,\gamma,\mu$ have at most two parts then $\lrn$ is equal to either $0$ or $1$. This is because only $1$'s may be placed in the top row of $\mu/\lambda$, and all the $2$'s must be right aligned in the bottom row. The remaining $1$'s must then be placed left of the $2$'s.\\
\two \ When the skew diagram $\mu/\lambda$ consists of a single row or column of size $k$, the reverse lattice word condition implies that $\lrn=1$ when $\gamma=(k)$ or $\gamma=(1,\dots,1)\in\ZZ^k$ respectively, otherwise $\lrn=0$.\\
\three \ Let $\lambda,\mu,\gamma$ be partitions such that $\modlam+\modgamma=\modmu$. Then $\lrn=c_{\gamma,\lambda}^\mu$, and moreover, $\lrn=0$ if either $\lambda$ or $\gamma$ is not contained in $\mu$; see \cite[\S 5.2~Corollary.~2]{Ful97}.
\end{remark}

We now construct the \emph{skew-Schur functor} using \cite[Ex~6.19]{FH91}. Let $V=\kk^n$ and fix a basis $u_1,\dots,u_n$ of $V$. Let $\lambda<\mu$ be partitions with at most $r$ parts and  set $d=\modmu-\modlam$. For any filling of the skew diagram $\mu/\lambda$ with entries in $\{1,\dots,n\}$, there is a corresponding basis vector of $V^{\otimes d}$ given by reading the content of the boxes from left to right, top to bottom. For example, if $\mu/\lambda=(3,2)/(1,0)$ and $n=4$ then
\[
\young(:31,43) \ \longleftrightarrow \ u_3\otimes u_1 \otimes u_4 \otimes u_3 \in V^{\otimes 4}.
\]
Consider the action of $S_d$ (the permutation group on $\{1,\dots,d\}$) on $V^{\otimes d}$ by permuting the indices, i.e.\ for $\sigma\in S_d$ and $v_i\in\cB$ we have $(v_1\otimes\cdots\otimes v_d)\cdot \sigma= v_{\sigma(1)}\otimes\cdots\otimes v_{\sigma(d)}$. Define the subgroups
\begin{align*}
P_{\text{row}}&=\left\{\sigma \in S_d \mid \sigma \  \text{preserves the content of each row of} \  \mu/\lambda\right\}, \\
P_{\text{col}}&=\left\{\sigma \in S_d \mid \sigma \  \text{preserves the content of each column of} \  \mu/\lambda\right\}.
\end{align*}
Now consider the group algebra $\kk S_d$ with generators $e_\sigma$ and define the elements
\begin{align*}
a_{\mu/\lambda}:=\sum_{\sigma\in P_{\text{row}}}e_\sigma, \quad
b_{\mu/\lambda}:=\sum_{\sigma\in P_{\text{col}}}\sign(\sigma) e_\sigma.
\end{align*}
These define endomorphisms on $V^{\otimes d}$ by setting $e_\sigma(v)=v\cdot\sigma$, and we have
\begin{align*}
\im(a_{\mu/\lambda})&\cong\Sym^{\mu_1-\lambda_1}V\otimes\dots \otimes \Sym^{\mu_r-\lambda_r}V, \\
\im(b_{\mu/\lambda})&\cong\bw^{\mu'_1-\lambda'_1}V\otimes\dots \otimes \bw^{\mu'_k-\lambda'_k}V,
\end{align*}
where $\lambda',\mu'$ are the conjugate partitions to $\lambda,\mu$ (list the heights of the columns left to right instead of the lengths of the rows top to bottom).

\begin{definition}
\label{defn:young_sym_skewSchur}
The \emph{Young symmetrizer} with respect to the skew diagram $\mu/\lambda$ is defined by
\[
c_{\mu/\lambda}:=b_{\mu/\lambda} a_{\mu/\lambda}.
\]
The image of $c_{\mu/\lambda}$ on a vector in $V^{\otimes d}$ is given by summing over the symmetrization of rows followed by the anti-symmetrization of columns. This defines an endomorphism on $V^{\otimes d}$; we call its image the \emph{skew-Schur power} and write
\begin{align*}
\SS^{\mu/\lambda}V:=\im(c_{\mu/\lambda}).
\end{align*}
By \cite[Ex~6.19]{FH91}, the basis $u_1,\dots,u_n$ of $V$ induces a basis of $\SS^{\mu/\lambda}V$ given by the images of semi-standard skew tableau with shape $\mu/\lambda$ filled with integers from $\{1,\dots,n\}$ under $c_{\mu/\lambda}$.
\end{definition}

\begin{example}
\label{exm:young_syms}
Let $\lambda=(1,0)$, $\mu=(2,2)$ and $u_1,u_2,u_3$ be a basis of $V$. We have $d=\modmu-\modlam=3$. Consider the semi-standard skew tableau $\scriptsize\young(:1,23)$ \ , which corresponds to the basis vector $u_1\otimes u_2\otimes u_3\in V^{\otimes 3}$. Symmetrizing the rows, we have $a_{\mu/\lambda}=\scriptsize\young(:1,23) \ + \ \scriptsize\young(:1,32) \ $, and now anti-symmetrizing the columns gives
\[
c_{\mu/\lambda}\left( \ \young(:1,23) \ \right)= \ \young(:1,23) \ - \ \young(:3,21) \ + \ \young(:1,32) \ - \ \young(:2,31) \ .
\]
Thus, $c_{\mu/\gamma}(u_1\otimes u_2\otimes u_3)=u_1\otimes u_2\otimes u_3-u_3\otimes u_2\otimes u_1+u_1\otimes u_3\otimes u_2-u_2\otimes u_3\otimes u_1$. In general therefore, $\SS^{(2,2)/(1,0)}V$ is the subspace of $V^{\otimes 3}$ spanned by vectors of the form
\[
u_{i_1}\otimes u_{i_2}\otimes u_{i_3}-u_{i_3}\otimes u_{i_2}\otimes u_{i_1}+u_{i_1}\otimes u_{i_3}\otimes u_{i_2}-u_{i_2}\otimes u_{i_3}\otimes u_{i_1}
\]
where $\scriptsize\young(:\ione,\itwo\ithree)$ is a semi-standard skew tableau with each $i_j\in\{1,2,3\}$.
\end{example}

When $\lambda=(0)$, the skew-Schur power reduces to the \emph{Schur power} $\SS^\mu V$, and in this case we write $c_\mu$ rather than $c_{\mu/\lambda}$ for the corresponding Young symmetrizer. The Schur power $\SS^\mu V$ is precisely the irreducible polynomial representation of $\GL(V)$ with highest weight $\mu$, and moreover these are all of the irreducible polynomial representations of $\GL(V)$; see \cite[Theorem~2,~p.114]{Ful97}. Since $\SS^\mu V$ consists of all irreducible polynomial representations of $\GL(V)$ as $\mu$ varies, we have the following decomposition results.

\begin{proposition}[{\cite[\S 6.1~Eqn.~6.7,~Ex.~6.19]{FH91}}]
\label{prop:skewschurdecomp}
Let $V$ be an $r$-dimensional vector space and let $\lambda,\mu$ be partitions with at most $r$ parts.\\
\one \ \emph{(Littlewood-Richardson rule)} We have
\[
\SS^\lambda V \otimes \SS^\mu V \cong \bigoplus_{\gamma}(\SS^\gamma V)^{\oplus c_{\lambda,\mu}^{\gamma}},
\]
where $\gamma$ ranges over all partitions satisfying $\modgamma = \modlam + \modmu$.\\
\two \ Suppose $\lambda\leq \mu$. Then $\SS^{\mu/\lambda} V$ is a polynomial representation of $\GL(V)$ with irreducible decomposition
\[
\SS^{\mu/\lambda} V \cong \bigoplus_{\gamma}(\SS^\gamma V)^{\oplus \lrn},
\]
where $\gamma$ ranges over all partitions satisfying $\modgamma = \modmu - \modlam$.
\end{proposition}

The \emph{Pieri rules} are two special cases of the Littlewood-Richardson rule.

\begin{proposition}[Pieri rules {\cite[Eqn~6.8-9]{FH91}}]
\label{prop:pieri}
Let $\lambda$ be a partition and $m\in \NN$. Then
\[
\SS^\lambda V\otimes\Sym^m V \cong \bigoplus_{\gamma} \SS^\gamma V \quad \text{and} \quad \SS^\lambda V\otimes\bw^m V \cong \bigoplus_{\mu} \SS^\mu V
\]
where $\gamma$ and $\mu$ range over all partitions formed by adding $m$ boxes to $\lambda$ with no two new boxes in the same column or row respectively.
\end{proposition}

Note that while in this section we have focussed on vector spaces, Schur powers can be constructed similarly in many other categories, in particular vector bundles; if $\mathcal{V}$ is a vector bundle whose fibre at a point $p$ is the vector space $V$, then $\SS^\lambda \mathcal{V}$ is the vector bundle whose fibre at $p$ is $\SS^\lambda V$.

\subsection{Generators of Kapranov's tilting algebra}
\label{sec:tilt_alg_gens}

Let $1<r\leq n$ be positive integers and fix $V=\kk^n$. Denote by $\Gr(n,r)$ the Grassmannian of $r$-dimensional quotients of $V$. Since $\Gr(n,r)\cong \Gr(n,n-r)$, we will assume that $n\geq 4$ and $1<r\leq n/2$. Let $\cW$ be the rank $r$ tautological quotient bundle on $\Gr(n,r)$. 
In \cite{Kap84}, Kapranov gave a presentation of $D^b(\Coh(\Gr(n,r)))$ by proving that the vector bundle
\begin{align}
\label{eqn:tilt_bundle_grass1}
E=\bigoplus_{\lambda\in \Young(n-r,r)} \SS^\lambda \cW
\end{align}
is a tilting bundle for $\Gr(n,r)$, and as such we refer to the algebra $A:=\End(E)$ as \emph{Kapranov's tilting algebra}. This algebra may be decomposed as the direct sum of the spaces $\Hom(\SS^\lambda \cW,\SS^\mu \cW)$ for all pairs $\lambda,\mu \in \Young(n-r,r)$. We will first prove a slightly more refined version of Kapranov's presentation of these homomorphism spaces.

\begin{proposition}[{\cite[p.189~3.0]{Kap84}}]
\label{prop:betterkap}
Let $\cW$ be the rank $r$ tautological quotient bundle of $\Gr(n,r)$ and let $\lambda\leq\mu\in\Young(n-r,r)$. Then
\begin{align}
\label{eqn:betterkap}
\Hom(\SS^\lambda \cW,\SS^\mu \cW)\cong\SS^{\mu/\lambda}V.
\end{align}
\begin{proof}
Kapranov's original result states that $\Hom(\SS^\lambda \cW,\SS^\mu \cW)\cong\bigoplus_{\gamma}\SS^{\gamma} V$
where $\gamma$ ranges over the positive summands in the decomposition of $\SS^{(-\lambda_r,\dots,-\lambda_1)}\cW \otimes \SS^{(\mu_1,\dots,\mu_r)}\cW$ into irreducibles. We will find precisely the multiplicities of each $\gamma$ in this decomposition. For partitions with negative entries we use the identity \cite[Eqn~0.1]{Kap84},
\begin{align}
\label{eqn:Wdual}
\SS^{(-\lambda_r,\dots,-\lambda_1)}\cW \cong \SS^\lambda (\cW^\vee) \cong (\SS^\lambda \cW)^\vee ,
\end{align}
where $\cW^\vee$ denotes the dual bundle of $\cW$. We may deal with Schur powers of bundles with negative entries by multiplying and then dividing by a line bundle; see \cite[2.1,~p.187]{Kap84}. For $m\in\ZZ$, this implies the identity
\begin{align}
\label{eqn:detW}
\SS^\lambda \cW \cong \det(\cW)^{-m}\otimes \SS^{(\lambda_1+m,\dots,\lambda_r+m)} \cW.
\end{align}
As a consequence of \cite[Proposition~3.6.7]{Har77}, since $\SS^\lambda \cW$ is a vector bundle we have the isomorphism
\[
\Hom(\SS^\lambda \cW,\SS^\mu \cW)\cong H^0(\Gr(n,r),(\SS^\lambda \cW)^\vee \otimes \SS^\mu\cW).
\]
By combining \eqref{eqn:Wdual} and \eqref{eqn:detW} with $m=\lambda_1$, we have
\begin{align*}
(\SS^\lambda \cW)^\vee \otimes \SS^\mu\cW &\cong \SS^{(-\lambda_r,\dots,-\lambda_1)}\cW \otimes \SS^{(\mu_1,\dots,\mu_r)}\cW \\
&\cong \det(\cW)^{-\lambda_1}\otimes (\SS^{\lambdat} \cW \otimes \SS^\mu \cW)
\end{align*}
where $\lambdat:=(\lambda_1-\lambda_r,\lambda_1-\lambda_{r-1},\dots,\lambda_1-\lambda_2,0)$. The decomposition of $\SS^{\lambdat} \cW \otimes \SS^\mu \cW$ into irreducibles ranges over partitions of size $\modlamt+\modmu$, but then multiplying back by $\det(\cW)^{-\lambda_1}$ results in partitions $\gamma$ of size $\modmu-\modlam$, many of which contain negative entries. However, when taking global sections these vanish by \cite[Lem~3.2a]{Kap88}, and so we are left with those $\gamma$ containing only non-negative entries and satisfying $\modgamma=\modmu-\modlam$.

It remains to find the multiplicities of each summand. The multiplicity of $\SS^\gamma\cW$ in $(\SS^\lambda \cW)^\vee \otimes \SS^\mu\cW$ is given by
\begin{align*}
\dim(\Hom(\SS^\gamma \cW, (\SS^\lambda \cW)^\vee \otimes \SS^\mu\cW))&
=\dim(\Hom(\SS^\gamma \cW\otimes \SS^\lambda \cW, \SS^\mu\cW))\\
&=\dim(\Hom(\oplus_{\mu^\prime}(\SS^{\mu^\prime}\cW)^{\oplus c_{\lambda,\gamma}^{\mu^\prime}}, \SS^\mu\cW))\\
&=\dim(\oplus_{\mu^\prime}\Hom((\SS^{\mu^\prime}\cW)^{\oplus c_{\lambda,\gamma}^{\mu^\prime}}, \SS^\mu\cW))\\
&=\dim(\oplus_{\mu^\prime}\Hom(\SS^{\mu^\prime}\cW, \SS^\mu\cW)^{\oplus c_{\lambda,\gamma}^{\mu^\prime}})\\
&=\dim(\Hom(\SS^\mu\cW, \SS^\mu\cW)^{\oplus\lrn})\\
&=\lrn
\end{align*}
where $\mu^\prime$ ranges over $\lvert \mu^\prime \rvert=\modgamma+\modlam=\modmu$ and the fifth and sixth equalities follow from \cite[3.5,~p.490]{Kap88}. Thus, we have shown that $\Hom(\SS^\lambda \cW,\SS^\mu \cW)=\bigoplus_{\gamma}(\SS^\gamma V)^{\oplus \lrn}$, and the identity from \Cref{prop:skewschurdecomp} completes the proof.
\end{proof}
\end{proposition}

\begin{corollary}
\label{cor:length1arrows_and_syms}
Let $\lambda\in\Young(n-r,r)$ and let $e_1,\dots,e_r$ denote the standard basis of $\ZZ^r$. Then for all $1\leq i\leq r$ and all $m>0$ such that $\lambda+me_i\in \Young(n-r,r)$, we have
\[
\Hom(\SS^\lambda \cW, \SS^{\lambda+me_i}\cW)\cong \Sym^m V.
\]
\begin{proof}
By \Cref{prop:betterkap} we have $\Hom(\SS^\lambda \cW, \SS^{\lambda+me_i}\cW)\cong \bigoplus_{\gamma}(\SS^\gamma V)^{\oplus c_{\lambda,\gamma}^{\lambda+me_i}}$, where $\gamma$ ranges over all partitions with $m$ boxes. The skew diagram $(\lambda+me_i)/\lambda$ however is just a single row of length $m$, hence by \Cref{rem:basiclrn}\two \ the only non-zero Littlewood-Richardson number in this decomposition occurs when $\gamma=(m)$, in which case $c_{\lambda,(m)}^{\lambda+me_i}=1$ and the result follows since $\SS^{(m)}V=\Sym^m V$.
\end{proof}
\end{corollary}

As a result of \Cref{prop:betterkap}, $\Hom(\SS^\lambda \cW, \SS^\mu \cW)$ depends only on the shape of the skew diagram $\mu/\lambda$ and so we may add  redundant columns to the left of both $\lambda$ and $\mu$ without changing $\Hom(\SS^\lambda\cW,\SS^\mu\cW)$. More precisely, we have the following.

\begin{corollary}
\label{cor:invariant_homs}
Let $\lambda,\mu \in\Young(n-r,r)$ with $\lambda\leq\mu$ and let $c\in \ZZ$ such that $c\geq -\lambda_r$. Then if both $(\lambda_1+c,\dots,\lambda_r+c),(\mu_1+c,\dots,\mu_r+c)\in \Young(n-r,r)$, we have $$\Hom(\SS^\lambda \cW, \SS^\mu \cW)\cong\Hom(\SS^{(\lambda_1+c,\dots,\lambda_r+c)}\cW,\SS^{(\mu_1+c,\dots,\mu_r+c)}\cW).$$ 
\end{corollary}

\subsection{Homomorphisms of adjacent summands}
\label{sec:schurmaps}
Hereafter we restrict our attention to $\Gr(n,2)$, i.e.\ $r=2$ and all partitions considered have at most two parts. Suppose $\lambda\in\Young(n-2,2)$. It will often be convenient to instead use the following alternative construction of $\SS^\lambda \cW$ (see \cite[\S 8]{Ful97}),
\begin{align}
\label{eqn:schur_2partsform}
\SS^\lambda \cW\cong \frac{(\bw^2 \cW)^{\otimes \lambda_2}\otimes \Sym^{\lambda_1-\lambda_2}\cW}{E_\lambda},
\end{align}
where $E_\lambda$ is the sub-bundle of \emph{exchange relations}, described below. In general therefore, write sections of $\SS^\lambda \cW$ as
\begin{align}
\label{eqn:w_lambda}
w^\lambda:=x_{1,1} \ww x_{1,2} \otimes \cdots \otimes x_{\lambda_2,1} \ww x_{\lambda_2,2} \otimes y_1\cdots y_{\lambda_1-\lambda_2}.
\end{align}
Note that each variable comprising $w^\lambda$ corresponds to a box of the Young diagram $\lambda$ in the obvious way. We may summarise the exchange relations as follows:
\begin{itemize}
\setlength\itemsep{0em}
	\item[(E1):] Take a single box from any column (where the symmetric part $\Sym^{(\lambda_1-\lambda_2)}\cW$ counts as $\lambda_1-\lambda_2$ distinct columns) and sum over both ways of swapping that box with the boxes in a height two column to the left of it. For example, using $w_\lambda$ choose $y_1$ and the first column (containing $x_{1,1} \ww x_{1,2}$) to get
\begin{align*}
w^\lambda_{E1}:=& \ y_1 \ww x_{1,2} \otimes \cdots \otimes x_{\lambda_2,1} \ww x_{\lambda_2,2} \otimes x_{1,1} y_2\cdots y_{\lambda_1-\lambda_2} \\
+& \ x_{1,1} \ww y_1 \otimes \cdots \otimes x_{\lambda_2,1} \ww x_{\lambda_2,2} \otimes x_{1,2} y_2\cdots y_{\lambda_1-\lambda_2}.
\end{align*}
Then modulo the exchange relations, we have $w^\lambda=w^\lambda_{E1}$.
\item[(E2):] Take both boxes in a height two column and swap them with another column of height two, maintaining the order of both boxes. For example, using $w_\lambda$ choose the first column and the one containing $x_{\lambda_2,1} \ww x_{\lambda_2,2}$ to get
\begin{align*}
w^\lambda_{E2}:=x_{\lambda_2,1} \ww x_{\lambda_2,2} \otimes \cdots \otimes x_{1,1} \ww x_{1,2} \otimes y_1\cdots y_{\lambda_1-\lambda_2}.
\end{align*}
Then modulo the exchange relations, we have $w^\lambda=w^\lambda_{E2}$.
\end{itemize}

Fix a basis $u_1,\dots,u_n$ of $V$. We now explicitly write down the homomorphisms between adjacent summands in the tilting bundle, i.e.\ those defining $\Hom(\SS^\lambda \cW,\SS^{\lambda+e_i}\cW)$ for all pairs $\lambda,\lambda+e_i\in\Young(n-2,2)$ with $i\in\{1,2\}$. By \Cref{cor:length1arrows_and_syms} these spaces are all isomorphic to $V$.

First consider $\lambda=(0,0)$ and $i=1$, so $\lambda+e_i=(1,0)$. Then
\begin{equation}
\label{eqn:sections}
\Hom(\cO_{\Gr(n,2)},\cW)\cong H^0({\Gr(n,2)},\cW)\cong V,
\end{equation}
so given a basis vector $v$ there is a homomorphism $s_v\colon \cO_{\Gr(n,2)} \rightarrow \cW$ and a uniquely determined global section $z_v:=s_v(1)$. Next suppose $\lambda=(1,0)$ and $\lambda+e_i=(1,1)$. We also have $\Hom(\cW,\bw^2 \cW)\cong V$ by \Cref{cor:length1arrows_and_syms}, so there is a homomorphism $s^\prime_v\colon \cW\rightarrow \bw^2 \cW$ which we can define using the same section: $s^\prime_v(x)=x \ww z_v$.

Unfortunately, writing down homomorphisms $\SS^\lambda \cW\rightarrow\SS^{\lambda+e_i} \cW$ is in general not as straightforward as adding in a new variable $z_v$ where required; we need to consider well-definedness with respect to the exchange relations (E1) and (E2).

\begin{proposition}
\label{prop:schurmaps}
Let $\lambda,\lambda+e_i\in\Young(n-2,2)$ where $i\in\{1,2\}$ and let $w^\lambda\in\SS^\lambda \cW$ be as in \eqref{eqn:w_lambda}. Let $v$ and $z_v$ be as above. Then the following are well-defined linear homomorphisms.
\begin{itemize}
	\setlength\itemsep{0em}
	\item[\one] For $i=1$, the map $f^\lambda_v\colon \SS^\lambda \cW\rightarrow \SS^{\lambda+e_1}\cW$ given by
	\[
	f^\lambda_v(w^\lambda)=x_{1,1} \ww x_{1,2} \otimes \cdots \otimes x_{\lambda_2,1} \ww x_{\lambda_2,2} \otimes y_1\cdots y_{\lambda_1-\lambda_2} z_v.
	\]
	\item[\two] For $i=2$, the map $g^\lambda_v\colon \SS^\lambda \cW\rightarrow \SS^{\lambda+e_2}\cW$ given by
	\[
	g^\lambda_v(w^\lambda)=\sum_{k=1}^{\lambda_1-\lambda_2} x_{1,1} \ww x_{1,2} \otimes \cdots \otimes x_{\lambda_2,1} \ww x_{\lambda_2,2} \otimes y_k \ww z_v \otimes \prod_{j\neq k} y_j.
	\]
\end{itemize}
Moreover, $f^\lambda_{u_1},\dots,f^\lambda_{u_n}$ or $g^\lambda_{u_1},\dots,g^\lambda_{u_n}$ form a basis for $\Homy(\SS^\lambda \cW, \SS^{\lambda+e_i}\cW)$ where $i=1$ or $2$ respectively.
\begin{proof}
The maps are clearly linear so we just prove they are well-defined with respect to the exchange relations on $\SS^\lambda\cW$. Using $w^\lambda$ from \eqref{eqn:w_lambda} and $w^\lambda_{E1},w^\lambda_{E2}$ defined above, we will show that $f^\lambda_v(w^\lambda)=f^\lambda_v(w^\lambda_{E1})=f^\lambda_v(w^\lambda_{E2})$ and $g^\lambda_v(w^\lambda)=g^\lambda_v(w^\lambda_{E1})=g^\lambda_v(w^\lambda_{E2})$. Any other choice of exchange is of the form (E1) or (E2) and the proof is almost identical.

\noindent \one \ If $i=1$, the symmetric power in $\SS^{\lambda+e_1} \cW$ compared to $\SS^\lambda \cW$ is increased by one while the alternating part is left unchanged; see \eqref{eqn:schur_2partsform}. It turns out that simply inserting $z_v$ into the symmetric part is well-defined. Firstly, we have
\begin{align*}
f^\lambda_v(w^\lambda_{E1})&=y_1 \ww x_{1,2} \otimes \cdots \otimes x_{\lambda_2,1} \ww x_{\lambda_2,2} \otimes x_{1,1} y_2\cdots y_{\lambda_1-\lambda_2} z_v \\
&\quad+x_{1,1} \ww y_1 \otimes \cdots \otimes x_{\lambda_2,1} \ww x_{\lambda_2,2} \otimes x_{1,2} y_2\cdots y_{\lambda_1-\lambda_2} z_v\\&=x_{1,1} \ww x_{1,2} \otimes \cdots \otimes x_{\lambda_2,1} \ww x_{\lambda_2,2} \otimes y_1\cdots y_{\lambda_1-\lambda_2} z_v \\&=f^\lambda_v(w^\lambda)
\end{align*}
where for the second equality we perform the inverse to an (E1) exchange. Secondly, since the exchange defining $w^\lambda_{E2}$ has no effect on the symmetric part where $z_v$ is added, $f^\lambda_v(w^\lambda_{E2})=f^\lambda_v(w^\lambda)$ is immediate.

\noindent \two \ If $i=2$, the symmetric power in $\SS^{\lambda+e_2} \cW$ decreases by one while the alternating power increases by one. The new height two column in $\lambda+e_2$ requires two variables; one of these will be $z_v$ while the other will be a variable removed from the symmetric part. To make this well-defined with respect to the exchange relations, we must sum over every choice of variable we remove from the symmetric part to pair with $z_v$, which leads to the definition of $g^\lambda_v$. For $w^\lambda_{E1}$, we have
\begin{align*}
g_v(w^\lambda_{E1})&=\sum_{k=1}^{\lambda_1-\lambda_2} y_1 \ww x_{1,2} \otimes \cdots \otimes x_{\lambda_2,1} \ww x_{\lambda_2,2} \otimes y_k \ww z_v \otimes x_{1,1} \prod_{j\neq k,1} y_j \\&\quad+\sum_{k=1}^{\lambda_1-\lambda_2} x_{1,1} \ww y_1 \otimes \cdots \otimes x_{\lambda_2,1} \ww x_{\lambda_2,2} \otimes y_k \ww z_v \otimes x_{1,2} \prod_{j\neq k,1} y_j \\&=\sum_{k=1}^{\lambda_1-\lambda_2} x_{1,1} \ww x_{1,2} \otimes \cdots \otimes x_{\lambda_2,1} \ww x_{\lambda_2,2} \otimes y_k \ww z_v \otimes \prod_{j\neq k} y_j \\&=g^\lambda_v(w^\lambda),
\end{align*}
where for the second equality we perform the inverses to the (E1)-type exchanges that move $y_1$ into the first column on each pair of terms from the two sums in turn. As for $w^\lambda_{E2}$, this is similar to \one \ because the exchange in question occurs left of the column where $g^\lambda_v$ inserts $z_v$, hence $g^\lambda_v(w^\lambda_{E2})=g^\lambda_v(w^\lambda)$ is immediate.
		
Finally, using the basis $u_1\dots,u_n$ and \eqref{eqn:sections}, we get a basis $s_{u_\rho}, 1\leq \rho\leq n$ of $\Hom(\cO_{\Gr(n,2)},\cW)$ and in turn a collection of linearly independent global sections $z_{u_\rho}$. Then for any given $\lambda$ the maps $f^\lambda_{u_\rho}$ (or $g^\lambda_{u_\rho}$) are also linearly independent: for $1\leq \rho\leq n$ the images of $w^\lambda$ under each $f^\lambda_{u_\rho}$ (or $g^\lambda_{u_\rho}$) are pairwise distinct, the sections $z_{u_\rho}$ are linearly independent and no sequence of exchanges will produce a linear dependence relation since exchanges never introduce new variables, only move around the existing ones.
\end{proof}
\end{proposition}

We will hereafter refer to the maps constructed in \Cref{prop:schurmaps} as `\emph{$f$-type}' and `\emph{$g$-type}' maps.

\subsection{Structure of the tilting quiver for $\Gr(n,2)$}
\label{sec:tilt_quiver_structure}

Recall that a \textit{quiver} $Q=(Q_0,Q_1)$ is a directed graph with vertex set $Q_0$ and arrow set $Q_1$. For each arrow $a\in Q_1$ we denote by $\hd(a),\tl(a)\in Q_0$ the vertices at the head and tail of $a$ respectively. Write $\kk Q$ for the path algebra of $Q$.

We now construct the \emph{tilting quiver} $\TQ$ for $\Gr(n,2)$, a quiver such that we may write down a surjective $\kk$-algebra homomorphism $\Phi\colon\kk\TQ\rightarrow A=\End(E)$, where $E$ is the tilting bundle given in \eqref{eqn:tilt_bundle_grass}. The vertex set $\TQ_0$ will be given by the irreducible summands of $E$, namely $\SS^\lambda \cW$ for all $\lambda\in \Young(n-2,2)$. Note that we will sometimes directly refer to vertices by $\lambda$ rather than $\SS^\lambda\cW$. Recall that $A=\End(E)$ may be decomposed as the collection of spaces $\Hom(\SS^\lambda \cW,\SS^\mu \cW)$ for all pairs $\lambda,\mu \in \Young(n-2,2)$. The arrow set $\TQ_1$ will be given by a minimal set of generators for the spaces satisfying $\lambda<\mu$. By \Cref{cor:length1arrows_and_syms}, for adjacent vertices we have $\Hom(\SS^\lambda \cW, \SS^{\lambda+e_i}\cW)\cong V$ and, depending on $i$, these spaces are spanned by a collection of $f$-type or $g$-type maps as defined in \Cref{prop:schurmaps}. Hence, for all pairs $\lambda,\lambda+e_i\in\Young(n-2,2)$ we will have $n$ arrows in $\TQ_1$ from $\SS^\lambda \cW\rightarrow\SS^{\lambda+e_i}\cW$ corresponding to the $f$-type or $g$-type basis of $\Hom(\SS^\lambda \cW, \SS^{\lambda+e_i}\cW)$ as appropriate.

\smallskip
\noindent \textbf{Claim:} For any pair $\lambda<\mu\in\Young(n-2,2)$, every map in $\Hom(\SS^\lambda \cW, \SS^\mu\cW)$ may be written as a linear combination of compositions of $f$-type and $g$-type maps.

\smallskip
A proof of this claim implies that the collection of $f$-type and $g$-type maps constitutes a minimal set of generators for the spaces $\Hom(\SS^\lambda \cW,\SS^\mu \cW)$ with $\lambda<\mu$. Therefore, the arrows between adjacent summands described above form the complete arrow set $\TQ_1$. We actually prove a stronger statement than in the claim, which is that the compositions formed strictly by a sequence of $f$-type maps followed by a sequence of $g$-type maps is enough. In other words, given $\lambda<\mu$ define $m_1=\mu_1-\lambda_1$ and $m_2=\mu_2-\lambda_2$, and for $0\leq k\leq m_1+m_2$ define the sequence of partitions
\begin{align}
\label{eqn:nu_sequence}
\tau_k:=
\begin{cases}
\lambda+ke_1 & \text{if} \ 0\leq k\leq m_1, \\
\lambda+m_1e_1+(k-m_1)e_2 & \text{if} \ m_1\leq k\leq m_1+m_2 .
\end{cases}
\end{align}
Then the claim follows from the following proposition.

\begin{theorem}
\label{prop:surj_onto_homs}
Let $\lambda<\mu\in\Young(n-2,2)$ and $\tau_k$ be the sequence of partitions defined in \eqref{eqn:nu_sequence}. Then the composition map
\[
\Theta_{\lambda,\mu}\colon\bigotimes_{k=1}^{m_1+m_2}\Hom(\SS^{\tau_{k-1}}\cW, \SS^{\tau_{k}}\cW) \longrightarrow \Hom(\SS^\lambda \cW, \SS^\mu\cW),
\]
where
\[
f^{\tau_0}_{v_1} \otimes \cdots \otimes f^{\tau_{m_1-1}}_{v_{m_1}} \otimes g^{\tau_{m_1}}_{v_{m_1+1}}\otimes\cdots\otimes g^{\tau_{m_1+m_2-1}}_{v_{m_1+m_2}} \mapsto g^{\tau_{m_1+m_2-1}}_{v_{m_1+m_2}}\circ \cdots \circ g^{\tau_{m_1+m_2-1}}_{v_{m_1+m_2}} \circ f^{\tau_{m_1-1}}_{v_{m_1}} \circ \cdots \circ f^{\tau_0}_{v_1},
\]
is surjective.
\end{theorem}
The proof is technical and we postpone it until \Cref{sec:surj_proof}. Assuming \Cref{prop:surj_onto_homs}, we can now properly define the tilting quiver $\TQ$ and establish the main result of this chapter.

\begin{definition}
\label{defn:tilting_quiver_gr_n2}
Define the tilting quiver $\TQ$ of ${\Gr(n,2)}$ by
\begin{align*}
\TQ_0&=\left\{\lambda \in \ZZ^2 \mid n-2\geq\lambda_1\geq\lambda_2\geq 0\right\},\\
\TQ_1&=\left\{a_\rho^{\lambda,i} \,\middle\vert\, 
\begin{array}{l}
1\leq \rho\leq n\\
i\in\{1,2\}, \ \lambda,\lambda+e_i\in \TQ_0\\
\tl(a_\rho^{\lambda,i})=\lambda, \  \hd(a_\rho^{\lambda,i})=\lambda+e_i
\end{array}\right\}.
\end{align*}
See \Cref{fig:grn2tiltQ}. Note that when drawing the tilting quiver we will use the presentation of $\SS^\lambda\cW$ given in \eqref{eqn:schur_2partsform} to label the vertices, though we will omit writing $E_\lambda$ on every vertex.
\begin{figure}[!ht]
	\begin{center}
		\begin{tikzpicture}
		\tikzset{edge/.style = {->,> = latex'}}
		\tikzset{font={\fontsize{8pt}{12}\selectfont}}
		
		\node[thick] (A) at  (0,0) {$\cO_{\Gr(n,2)}$};
		\node[thick] (B) at  (2,0) {$\cW$};
		\node[thick] (C) at  (4,0) {$\Sym^2 \cW$};
		\node[thick] (D) at  (2,2) {$\bw^2 \cW$};
		\node[thick] (E) at  (4,2) {$\bw^2 \cW \otimes \cW$};
		\node[thick] (F) at  (4,4) {$(\bw^2 \cW)^{\otimes 2}$};
		\node[thick] (G) at  (9,0) {$\Sym^{n-3}\cW$};
		\node[thick] (H) at  (12,0) {$\Sym^{n-2}\cW$};
		\node[thick] (I) at  (9,2) {\begin{tabular}{c} $\bw^2 \cW$ \\ $\otimes \Sym^{n-4}\cW$ \end{tabular}};
		\node[thick] (J) at  (12,2) {\begin{tabular}{c} $\bw^2 \cW$ \\ $\otimes \Sym^{n-3}\cW$ \end{tabular}};
		\node[thick] (K) at  (12,7) {\begin{tabular}{c} $(\bw^2 \cW)^{\otimes (n-4)}$ \\ $\otimes \Sym^2 \cW$ \end{tabular}};
		\node[thick] (L) at  (9,9) {$(\bw^2 \cW)^{\otimes(n-3)}$};		
		\node[thick] (M) at  (12,9) {\begin{tabular}{c} $(\bw^2 \cW)^{\otimes (n-3)}$ \\ $\otimes \cW$ \end{tabular}};
		\node[thick] (N) at  (12,11) {$(\bw^2 \cW)^{\otimes(n-2)}$};
		\node[thick] (CC) at  (6,0) {};
		\node[thick] (EE) at  (6,2) {};		
		\node[thick] (FF) at  (6,4) {};
		\node[thick] (GG) at  (7,0) {};
		\node[thick] (II1) at  (7,2) {};
		\node[thick] (II2) at  (9,4) {};		
		\node[thick] (JJ) at  (12,4) {};		
		\node[thick] (KK1) at  (12,5) {};	
		\node[thick] (KK2) at  (9.5,7) {};		
		\node[thick] (LL) at  (9,7) {\begin{tabular}{c} $(\bw^2 \cW)^{\otimes (n-4)}$ \\ $\otimes \cW$ \end{tabular}};
		\node[thick] (LL1) at  (9,5) {};
		\node[thick] (LL2) at  (6.8,7) {};
		
		\draw[edge] (A) -- (B);
		\draw[edge] (B) -- (C);
		\draw[edge] (B) -- (D);
		\draw[edge] (C) -- (E);
		\draw[edge] (D) -- (E);
		\draw[edge] (E) -- (F);
		\draw[edge] (G) -- (H);
		\draw[edge] (G) -- (I);
		\draw[edge] (H) -- (J);
		\draw[edge] (I) -- (J);
		\draw[edge] (K) -- (M);
		\draw[edge] (L) -- (M);
		\draw[edge] (M) -- (N);
		
		\draw[dashed, ->] (C) -- (CC);
		\draw[dashed, ->] (E) -- (EE);
		\draw[dashed, ->] (F) -- (FF);
		\draw[dashed, ->] (GG) -- (G);		
		\draw[dashed, ->] (II1) -- (I);		
		\draw[dashed, ->] (I) -- (II2);		
		\draw[dashed, ->] (J) -- (JJ);		
		\draw[dashed, ->] (KK1) -- (K);	
		\draw[dashed, ->] (KK2) -- (K);
		\draw[dashed, ->] (LL) -- (L);
		\draw[dashed, ->] (LL1) -- (LL);
		\draw[dashed, ->] (LL2) -- (LL);
		
		\draw[dotted] (5.5,1) -- (7.5,1);
		\draw[dotted] (5.5,3) -- (7.5,3);
		\draw[dotted] (10.5,3.5) -- (10.5,5.5);
		\draw[dotted] (5.5,4.7) -- (7,6.5);
		
		\end{tikzpicture}
	\end{center}
	\caption[The tilting quiver for $\Gr(n,2)$]{The tilting quiver $\TQ$ for $\Gr(n,2)$. Each arrow in the figure represents $n$ arrows in the quiver corresponding to a basis of $V$.}
	\label{fig:grn2tiltQ}	
\end{figure}
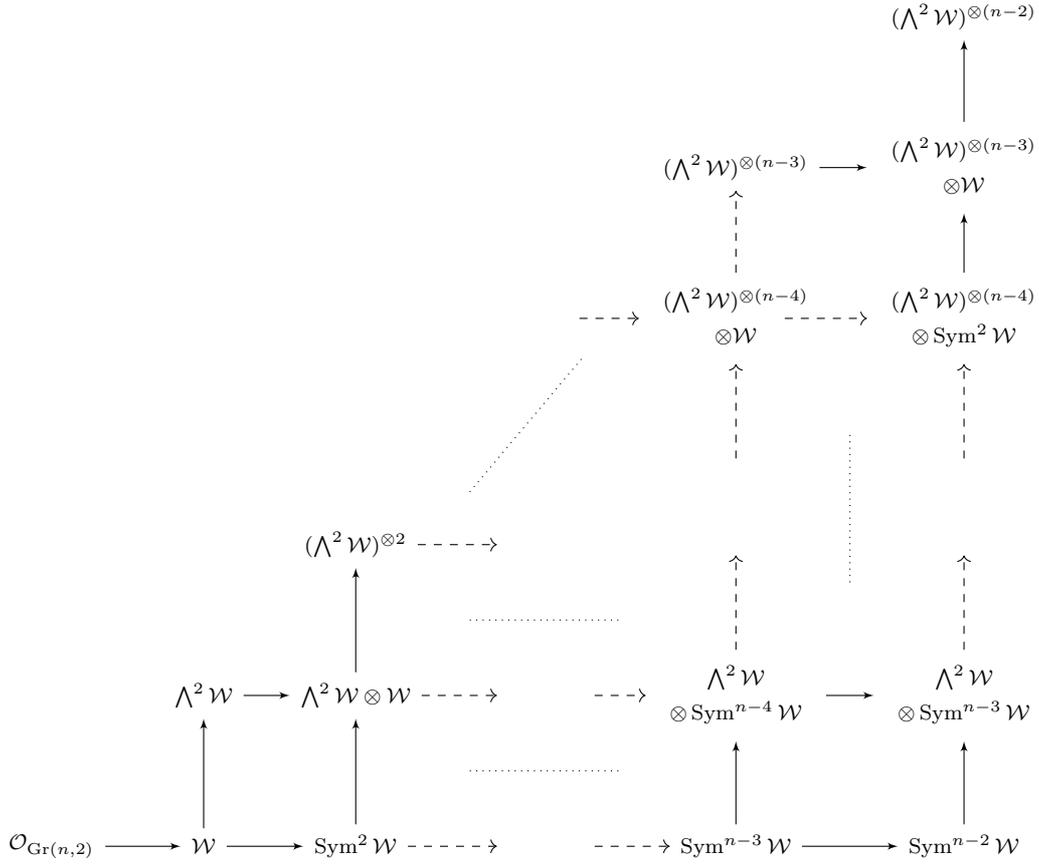
\end{definition}

For each $\lambda\in\TQ_0$ let $e_\lambda\in\kk\TQ$ denote the idempotent corresponding to the path of length zero at that vertex. For all $\lambda,\mu\in\TQ_0$ satisfying $\lambda<\mu$, with the convention that we traverse paths from right to left (the same way that composition of maps is performed), $e_\mu\kk\TQ e_\lambda$ denotes the space of paths $\lambda\rightarrow\mu$. Fix a basis $u_1,\dots,u_n$ of V.

\begin{definition}
\label{defn:k-alg-hom}
By definition of $\TQ_1$ and \Cref{cor:length1arrows_and_syms}, when $\mu=\lambda+e_i$ for $i\in\{1,2\}$, we have $e_{\lambda+e_i}\kk\TQ e_\lambda\cong V \cong \Hom(\SS^\lambda\cW,\SS^{\lambda+e_i}\cW)$, and by \Cref{prop:schurmaps} the latter space has a basis given by $f^\lambda_{u_1},\dots,f^\lambda_{u_n}$ if $i=1$ or $g^\lambda_{u_1},\dots,g^\lambda_{u_n}$ if $i=2$. Thus, for all $\lambda\in\TQ, 1\leq\rho\leq n$ and $i\in\{1,2\}$ as appropriate, we define a $\kk$-algebra homomorphism
\begin{align}
\label{eqn:k-alg_surj}
\Phi\colon\kk\TQ\longrightarrow A
\end{align}
as follows:
\begin{align*}
\Phi(e_\lambda)&=\text{id}_\lambda\in\Hom(\SS^\lambda\cW,\SS^\lambda\cW)\cong\kk, \\
\Phi(a_\rho^{\lambda,i})&=
\begin{cases}
f^\lambda_{u_\rho}\in \Hom(\SS^\lambda\cW,\SS^{\lambda+e_1}\cW) & \text{if} \ i=1, \\
g^\lambda_{u_\rho}\in \Hom(\SS^\lambda\cW,\SS^{\lambda+e_2}\cW) & \text{if} \ i=2.
\end{cases}
\end{align*}
The images of the horizontal and vertical arrows in $\TQ_1$ are therefore exactly the $f$-type and $g$-type maps defined in \Cref{prop:schurmaps} respectively. We extend $\Phi$ to any path in $\TQ$ by mapping the concatenation of arrows $a_\rho^{\lambda,i}$ to the composition of maps in $A$ as appropriate and finally, we extend $\Phi$ linearly over $\kk$ to all linear combinations of paths in $\kk\TQ$.
\end{definition}

\begin{corollary}
\label{thm:tilting_quiver_rank2}
The $\kk$-algebra homomorphism $\Phi\colon\kk\TQ\rightarrow A$ is surjective.
\begin{proof}
Firstly, whenever $\lambda$ is not contained in $\mu$ we have $\Hom(\SS^\lambda\cW,\SS^\mu\cW)=0$, thus surjectivity is trivial in these cases; indeed, we defined no arrows in $\TQ$ for such pairs $\lambda,\mu$. Consequently, $\TQ$ is acyclic and $(0,0)$ is the unique source vertex. Now suppose $\lambda<\mu\in\TQ_0$ and $h\in\Hom(\SS^\lambda\cW,\SS^\mu\cW)\subset A$. \Cref{prop:surj_onto_homs} implies that $h$ may be factorised as a linear combination of compositions of $f$-type and $g$-type maps. Therefore the corresponding linear combination of paths given by concatenating arrows of the form $a_\rho^{\lambda,1},a_\rho^{\lambda,2}$ maps to $h$ under $\Phi$. Thus, $\Phi$ is surjective.
\end{proof}
\end{corollary}

\begin{remark}
\label{rem:BLV_newproof}
\Cref{thm:tilting_quiver_rank2} thus provides a new proof of \cite[Theorem~6.9]{BLV16} for $\Gr(n,2)$.
\end{remark}

\subsection{Proof of \Cref{prop:surj_onto_homs}}
\label{sec:surj_proof}

First of all, due to the invariance result \Cref{cor:invariant_homs}, it is enough to consider the special case where $\lambda_2=0$; the general case follows immediately since $\Hom(\SS^{(\lambda_1,\lambda_2)} \cW, \SS^{(\mu_1,\mu_2)} \cW)\cong\Hom(\SS^{(\lambda_1-\lambda_2,0)}\cW,\SS^{(\mu_1-\lambda_2,\mu_2-\lambda_2)}\cW)$, and we adjust the sequence $\tau_k$ accordingly. We therefore assume $\lambda_2=0$ throughout the entirety of the proof, and in turn we have $m_1=\mu_1-\lambda_1$, $m_2=\mu_2$.

We first consider the domain and codomain of $\Theta_{\lambda,\mu}$. By \Cref{cor:length1arrows_and_syms}, the domain of $\Theta_{\lambda,\mu}$ is isomorphic to $V^{\otimes (m_1+m_2)}$ and by \Cref{prop:betterkap}, the codomain is isomorphic to $\SS^{\mu/\lambda}V$, which by construction is a $\GL(V)$-submodule of $V^{\otimes (m_1+m_2)}$.

\begin{lemma}
\label{lem:skew_hom_decomp}
The irreducible decomposition of $\Hom(\SS^\lambda \cW, \SS^\mu\cW)$ is given by
\begin{align}
\label{eqn:hom_decomp}
\SS^{\mu/\lambda}V\cong \bigoplus_{\gamma\in\Gamma_{\mu/\lambda}} \SS^\gamma V,
\end{align}
where
\begin{itemize}
\item if $\mu_2\leq \lambda_1$, $\Gamma_{\mu/\lambda}$ consists of the partitions $(\max\{m_1,m_2\},\min\{m_1,m_2\}),$\\ $(\max\{m_1,m_2\}+1,\min\{m_1,m_2\}-1),\dots,(m_1+m_2,0)$.
\item if $\mu_2>\lambda_1$, $\Gamma_{\mu/\lambda}$ consists of the partitions $(\max\{m_1,m_2\},\min\{m_1,m_2\}),$\\ $(\max\{m_1,m_2\}+1,\min\{m_1,m_2\}-1),\dots, (\mu_1,\mu_2-\lambda_1)$.
\end{itemize}
\begin{proof}
The main tool for this is \Cref{prop:skewschurdecomp}, which tells us that
\[
\SS^{\mu/\lambda} V \cong \bigoplus_{\gamma}(\SS^\gamma V)^{\oplus \lrn}
\]
where $\gamma$ ranges over all partitions satisfying $\modgamma = \modmu - \modlam=m_1+m_2$. By \Cref{rem:basiclrn}\three \ we have $\lrn\neq0\implies \gamma\leq \mu$, hence we only need to consider $\gamma$ with at most two parts that satisfy $\gamma_1\leq\mu_1$ and $\gamma_2\leq\mu_2$. Additionally, $\lrn$ is equal to either $0$ or $1$ by \Cref{rem:basiclrn}\one.

The set $\Gamma_{\mu/\lambda}$ is given by the collection of such $\gamma$ satisfying $\lrn=1$, i.e.\ those such that the skew diagram $\mu/\lambda$ filled with content $\gamma$ is a Littlewood-Richardson tableau. First suppose that $\mu_2\leq\lambda_1$, which means the two rows of $\mu/\lambda$ do not overlap. Since there are no columns of height two, the strictly increasing columns condition cannot be broken and so a filling consisting of only $1$'s, i.e.\ $\gamma=(m_1+m_2,0)$, is permissible. By \Cref{rem:basiclrn}\one, $2$'s may only be placed in the right of the bottom row; therefore, to avoid breaking the reverse lattice word condition, the number of $2$'s we may insert is bounded above by the length of the top row, thus $\gamma_2\leq m_1$. Since we must also have $\gamma_2\leq \mu_2=m_2$, we have $\min\{m_1,m_2\}\geq\gamma_2\geq0$ as required. 

The argument for the case when $\mu_2>\lambda_1$ carries over from above, except we now also have a positive lower bound for $\gamma_2$ because columns of height two exist in $\mu/\lambda$. We therefore require at least $\mu_2-\lambda_1$ many $2$'s to place in the height two columns, and this gives the second case in the statement of the lemma.
\end{proof}
\end{lemma}

We now outline the strategy of proof for \Cref{prop:surj_onto_homs}. Fix $\lambda<\mu$. We will show that the restriction of $\Theta_{\lambda,\mu}$ to the submodule $\SS^{\mu/\lambda}V$ of $V^{\otimes (m_1+m_2)}$ is an isomorphism, and it will therefore follows that $\Theta_{\lambda,\mu}$ is surjective.  To do this we take advantage of \text{Schur's lemma} (see \cite[Lemma~1.7]{FH91}): if $\varphi\colon W_1\rightarrow W_2$ is a $G$-module homomorphism of irreducible $G$-modules, then $\varphi$ is either an isomorphism or zero. Recall $\Gamma_{\mu/\lambda}$ from \Cref{lem:skew_hom_decomp}. We have that $\SS^\gamma V$ is irreducible over $\GL(V)$ for all $\gamma\in\Gamma_{\mu/\lambda}$ and $\Theta_{\lambda,\mu}$ is a $\GL(V)$-module homomorphism. Since every summand in \eqref{eqn:hom_decomp} appears with multiplicity one, Schur's lemma implies that it is enough to show $\Theta_{\lambda,\mu}$ is non-zero when restricted to  $\SS^\gamma V\subset V^{\otimes (m_1+m_2)}$ for each $\gamma\in\Gamma_{\mu/\lambda}$. Hence, the proof of \Cref{prop:surj_onto_homs} is completed by the following lemma.

\begin{lemma}
\label{lem:gamma_map_nonzero}
Let $\gamma\in\Gamma_{\mu/\lambda}$. Then $\Theta_{\lambda,\mu}\big|_{\SS^\gamma V} \neq 0$.
\end{lemma}

The strategy for the proof of \Cref{lem:gamma_map_nonzero} is as follows. For each $\gamma\in\Gamma_{\mu/\lambda}$ we first write down an element $h_\gamma$ in $V^{\otimes (m_1+m_2)}$, the domain of $\Theta_{\lambda,\mu}$. We then get an element $c_\gamma(h_\gamma)$ of $\SS^\gamma V$ by applying the Young symmetrizer $c_\gamma\colon V^{\otimes m_1+m_2}\twoheadrightarrow \SS^\gamma V$; see \Cref{exm:young_syms}. Then, by formulating the evaluation of a section $w^\lambda$ under the map  $\Theta_{\lambda,\mu}(c_\gamma(h_\gamma))$, we show this is non-zero to complete the proof. In order to implement this strategy we need some new notation.

\begin{notation}
\label{not:superscript_drop}
\noindent \one \ Fix a basis $u_1,\dots,u_n$ of $V$ and let  $\lambda<\mu\in\Young(n-2,2)$. Recall $m_1=\mu_1-\lambda_1$ and $m_2=\mu_2$ (we are assuming $\lambda_2=0$), so that $\mu=\lambda+m_1e_1+m_2e_2$. In this section we always consider compositions of maps given by the image of $\Theta_{\lambda,\mu}$, i.e.\ those of the form
\[
g^{\tau_{m_1+m_2-1}}_{v_{m_1+m_2}}\circ \cdots \circ g^{\tau_{m_1+m_2-1}}_{v_{m_1+1}} \circ f^{\tau_{m_1-1}}_{v_{m_1}} \circ \cdots \circ f^{\tau_0}_{v_1}
\]
where each $v_i$ is one of the basis vectors $u_1,\dots,u_n$. Because the domain of each map in this composition may be derived from whether the previous map is $f$-type or $g$-type, we will suppress the superscript of all maps except the first and instead write the above as
\[
g_{v_{m_1+m_2}}\circ\cdots\circ g_{v_{m_1+1}}\circ f_{v_{m_1}}\circ\cdots\circ f^\lambda_{v_1}.
\]

\noindent \two \ We require the notion of \emph{multisets}; these are sets with possible multiplicities of elements, and we distinguish multisets from sets by using square brackets $[ \ ]$. The cardinality of a multiset counts these multiplicities, e.g. $[1,1,2]$ has cardinality $3$. Denote \emph{ordered multisets} using $[ \ ]^o$, and given a multiset $M$ define $\text{O}_k(M)$ to be the collection of all ordered sub-multisets of $M$ with cardinality $k$. If the cardinality of $M$ is $m\geq0$, then the size of the collection $\text{O}_k(M)$ is $m!/(m-k)!$. For example, if $M=[1,1,2]$ then $\text{O}_2(M)$ is the collection  $[1,1]^o,[1,1]^o,[1,2]^o,[1,2]^o,[2,1]^o,[2,1]^o$.
\end{notation}

For every $\gamma\in\Gamma_{\mu/\lambda}$, we now define our candidates $h_\gamma\in V^{\otimes (m_1+m_2)}$ for use in the proof of \Cref{lem:gamma_map_nonzero}. Recall that the Young symmetrizer $c_\gamma$ is defined by $b_\gamma a_\gamma$ (see the discussion prior to \Cref{defn:young_sym_skewSchur}), so $c_\gamma(h_\gamma)$ will be a double sum taken over all ways of first symmetrizing the rows of $\gamma$, followed by anti-symmetrizing the columns. We will simplify matters by defining $h_\gamma$ such that $a_\gamma$ is trivial. The simplest such map is the basis vector of $V^{\otimes (m_1+m_2)}$ corresponding to the skew tableau of shape $\gamma$ with the top row filled with $1$'s and the bottom row filled with $2$'s, i.e.\
\begin{align*}
&\overbrace{\young(1111111,2222)}^{\gamma_1} \ \longleftrightarrow \ \underbrace{f^\lambda_{u_1} \otimes \cdots \otimes f_{u_1}}_{m_1} \otimes \underbrace{g_{u_1}\otimes\cdots\otimes g_{u_1}\otimes g_{u_2}\otimes\cdots\otimes g_{u_2}}_{m_2}=:h_\gamma.\\[-2em]
&\underbrace{\phantom{\young(2222)}}_{\gamma_2} \phantom{\young(111) \ \longleftrightarrow \ }  \underbrace{\phantom{f^\lambda_{u_1} \otimes \cdots \otimes f_{u_1} \otimes g_{u_1}\otimes\cdots\otimes g_{u}}}_{\gamma_1}\phantom{\otimes \ } \underbrace{\phantom{g_{u_2}\otimes\cdots\otimes g_{u_2}}}_{\gamma_2}
\end{align*}
Thus, there are $m_1$ $f$-type maps followed by $m_2$ $g$-type maps, the first $\gamma_1$ of which are defined using the basis vector $u_1$ and the last $\gamma_2$ defined using the basis vector $u_2$. By the conditions on $\gamma\in\Gamma_{\mu/\lambda}$, we always have $\gamma_2\leq\min\{m_1,m_2\}$.

We now write down $c_\gamma(h_\gamma)=b_\gamma a_\gamma(h_\gamma)$. By construction, $a_\gamma$ acts trivially on $h_\gamma$ so it remains to sum over all ways of anti-symmetrizing the height two columns of $\gamma$. Therefore
\[
c_\gamma(h_\gamma)=b_\gamma(h_\gamma)=\sum_{\sigma\in P_{\text{col}}}\sign(\sigma)h_\gamma\cdot \sigma,
\]
which, since $\gamma$ has $\gamma_2$ columns of height two, is a sum of $2^{\gamma_2}$ terms. Each $\sigma$ defines a unique ordered multiset
\[
\uve=[v_1,\dots,v_{\gamma_2}]^o,
\]
where each $v_i$ is equal to either $u_1$ or $u_2$ as given by the top row of the height two columns in $h_\gamma\cdot \sigma$. Define $\alpha_\sigma$ and $\beta_\sigma$ to be the number of $v_1$'s and $v_2$'s appearing in $\uve$ respectively; then $\alpha_\sigma+\beta_\sigma=\gamma_2$ and we have $\sign(\sigma)=(-1)^{\beta_\sigma}$. Given $\uve$ define $v'_i$ for $1\leq i\leq \gamma_2$ to be the vectors on the bottom row of the height two columns in $h\cdot \sigma$, i.e.\ if $v_i=u_1$ then $v'_i=u_2$ and vice versa. This yields $c_\gamma(h_\gamma)$ equal to
\[
\sum_{\sigma\in P_{\text{col}}} (-1)^{\beta_\sigma} f^\lambda_{v_1}\otimes\cdots\otimes f_{v_{\gamma_2}}\otimes f_{u_1}\otimes\cdots\otimes f_{u_1}\otimes g_{u_1}\otimes\cdots\otimes g_{u_1} \otimes g_{v'_1}\otimes\cdots\otimes g_{v'_{\gamma_2}}
\]
and therefore $\Theta_{\lambda,\mu}(c_\gamma(h_\gamma))$ is given by
\begin{align}
\label{eqn:c_gamma_h}
\sum_{\sigma\in P_{\text{col}}} (-1)^{\beta_\sigma} g_{v'_{\gamma_2}}\circ\cdots\circ g_{v'_1}\circ g_{u_1}\circ\cdots\circ g_{u_1}\circ f_{u_1}\circ\cdots\circ f_{u_1}\circ  f_{v_{\gamma_2}}\circ\cdots\circ f^\lambda_{v_1}.
\end{align}

\begin{notation}
\label{not:fg_maps_comp}
Before detailing the evaluation of a section $w^\lambda=y_1\cdots y_{\lambda_1}\in\SS^\lambda\cW\cong\Sym^{\lambda_1}\cW$ under the above sum, we will first describe the evaluation of $w^\lambda$ under a single composition $g_{v_{m_1+m_2}}\circ\cdots\circ g_{v_{m_1+1}}\circ f_{v_{m_1}}\circ\cdots\circ f^\lambda_{v_1}\in \Hom(\SS^\lambda\cW,\SS^\mu\cW)$ where each $v_i$ is one of the basis vectors $u_1,\dots,u_n$; this will also be useful in \Cref{chap:relnsgrn2}. Using \Cref{prop:schurmaps}, evaluating the composition of the $f$-type maps is easy: we simply have
\[
f_{v_{m_1}}\circ\cdots\circ f^\lambda_{v_1}(w^\lambda)=y_1\cdots y_{\lambda_1} z_{v_1} \cdots z_{v_{m_1}}.
\]
The evaluation of a $g$-type map is the sum over each way of pairing a variable in the symmetric part with the new variable introduced by the map, hence the evaluation of a succession of $g$-type maps is given by summing over all ordered ways of doing this. Thus, define the multiset $M=[y_1,\dots, y_{\lambda_1},z_{v_1},\dots,z_{v_{m_1}}]$, and recall that $O_{m_2}(M)$ is the collection of ordered sub-multisets of $M$ with cardinality $m_2$. Then the image of $w^\lambda$ under the composition above is given by
\begin{align}
\label{eqn:point_image}
w^\lambda\mapsto \sum_{\substack{X=[x_1,\dots,x_{m_2}]^o \\ \in \text{O}_{m_2}(M)}} x_1 \ww z_{v_{m_1+1}} \otimes \cdots \otimes x_{m_2} \ww z_{v_{m_1+m_2}} \otimes \prod_{z\in M\setminus X} z.
\end{align}
\end{notation}

\noindent \textbf{Proof of \Cref{lem:gamma_map_nonzero}.} Now recall \eqref{eqn:c_gamma_h}; the goal is to find a section $w^\lambda$ such that $\Theta_{\lambda,\mu}(c_\gamma(h_\gamma))(w^\lambda)\neq0$. Define $\delta$ to be the number of $g$-type maps in each term of $c_\gamma(h_\gamma)$ with fixed defining basis vector $u_1$, i.e.\ $\delta=m_2-\gamma_2$. Set
\[
w^\lambda:=z_{u_1}^{\lambda_1-\delta}z_{u_2}^{\delta}.
\]
Note that this choice of section is possible because $\gamma\in\Gamma_{\mu/\lambda}\implies \gamma_2\geq \mu_2-\lambda_1=m_2-\lambda_1\implies \lambda_1\geq m_2-\gamma_2=\delta$.

Since $\Theta_{\lambda,\mu}$ is linear we will analyse the evaluation of $w^\lambda$ under each term $\Theta_{\lambda,\mu}(h_\gamma\cdot\sigma)$ in the sum \eqref{eqn:c_gamma_h} separately; we will deal with the sign $(-1)^{\beta_\sigma}$ at the end. Thus, fix $\sigma\in P_{\text{col}}$ and consider $\uve,\alpha_\sigma,\beta_\sigma$ as defined above. Following \Cref{not:fg_maps_comp}, the composition of the $f$-type maps in $\Theta_{\lambda,\mu}(h_\gamma\cdot\sigma)$ contributes $z_{u_1}^{m_1-\beta_\sigma}$ and $z_{u_2}^{\beta_\sigma}$ to $w^\lambda$, bringing the total exponent of $z_{u_1}$ to $\lambda_1-\delta+m_1-\beta_\sigma=\mu_1-\mu_2+\alpha_\sigma$ and the total exponent of $z_{u_2}$ to $\delta+\beta_\sigma=\mu_2-\gamma_2+\beta_\sigma$. Hence, the multiset $M$ in \eqref{eqn:point_image} is given by
\[
M=[\underbrace{z_{u_1},\dots,z_{u_1}}_{\mu_1-\mu_2+\alpha_\sigma},\underbrace{z_{u_2},\dots,z_{u_2}}_{\mu_2-\gamma_2+\beta_\sigma}].
\]
After composing the remaining $g$-type maps, $\Theta_{\lambda,\mu}(h_\gamma\cdot\sigma)(w^\lambda)$ is equal to
\begin{align}
\label{eqn:h_sigma_eval}
\text{\footnotesize $ \sum_{\substack{X=[x_1,\dots,x_{m_2}]^o \\ \in \text{O}_{m_2}(M)}} x_1 \ww z_{u_1} \otimes \cdots \otimes x_{m_2-\gamma_2} \ww z_{u_1} \otimes x_{m_2-\gamma_2+1} \ww z_{v'_1} \otimes\cdots\otimes x_{m_2} \ww z_{v'_{\gamma_2}} \otimes \prod_{z\in M\setminus X} z. $}
\end{align}
Every variable in this sum is either $z_{u_1}$ or $z_{u_2}$, thus the only $X\in\text{O}_{m_2}(M)$ that produce a non-zero term are those with $x_1,\dots,x_{m_2-\gamma_2}$ equal to $z_{u_2}$ and for $m_2-\gamma_2+1\leq j \leq m_2$, $x_j=z_{u_1}$ if $v'_{j-m_2+\gamma_2}=z_{u_2}$ and vice versa. Hence every ordered sub-multiset $X$ that produces a non-zero term is identical, and consists of $\alpha_\sigma$ many $z_{u_1}$'s and $\mu_2-\gamma_2+\beta_\sigma$ many $z_{u_2}$'s. Define the total number of such $X$ to be $\eta_\sigma$; this is given by the number of ordered ways of choosing $\alpha_\sigma$ many $z_{u_1}$'s from $M$, multiplied by the number of ordered ways of choosing $\mu_2-\gamma_2+\beta_\sigma$ many $z_{u_2}$'s from $M$, i.e.\
\[
\eta_\sigma=\frac{(\mu_2-\gamma_2+\beta_\sigma)! (\mu_1-\mu_2+\alpha_\sigma)!}{(\mu_1-\mu_2)!}.
\]
The sum \eqref{eqn:h_sigma_eval} therefore simplifies to
\[
\eta_\sigma z_{u_2}\ww z_{u_1} \otimes\cdots\otimes z_{u_2}\ww z_{u_1} \otimes z_{v_1} \ww z_{v'_1} \otimes\cdots\otimes z_{v_{\gamma_2}} \ww z_{v'_{\gamma_2}} \otimes z_{u_1}^{\mu_1-\mu_2}.
\]
We now use anti-symmetrization in the columns with content $z_{v_i} \ww z_{v'_i}$ to ensure that the first entry is $z_{u_2}$ while the second is $z_{u_1}$. There are $\alpha_\sigma$ many columns that are not in this order, since this is the number of $z_{v'_i}$ terms that are equal to $z_{u_2}$. Hence, rearranging each column so that the content reads $z_{u_2}\ww z_{u_1}$ means we must multiply by $(-1)^{\alpha_\sigma}$. Therefore we simplify the above once more, yielding
\[
\Theta_{\lambda,\mu}(h_\gamma\cdot\sigma)(w^\lambda)=(-1)^{\alpha_\sigma} \eta_\sigma z_{u_2}\ww z_{u_1} \otimes\cdots\otimes z_{u_2}\ww z_{u_1} \otimes z_{u_1}^{\mu_1-\mu_2}.
\]
In conclusion, using \eqref{eqn:c_gamma_h} we have
\begin{align*}
\Theta_{\lambda,\mu}(c_\gamma (h_\gamma))(w^\lambda)&=\sum_{\sigma\in P_{\text{col}}} (-1)^{\beta_\sigma}(-1)^{\alpha_\sigma} \eta_\sigma z_{u_2}\ww z_{u_1} \otimes\cdots\otimes z_{u_2}\ww z_{u_1} \otimes z_{u_1}^{\mu_1-\mu_2}\\
&=(-1)^{\gamma_2} \sum_{\sigma\in P_{\text{col}}} \eta_\sigma z_{u_2}\ww z_{u_1} \otimes\cdots\otimes z_{u_2}\ww z_{u_1} \otimes z_{u_1}^{\mu_1-\mu_2}\\
&\neq 0,
\end{align*}
since $z_{u_2}\ww z_{u_1} \otimes\cdots\otimes z_{u_2}\ww z_{u_1} \otimes z_{u_1}^{\mu_1-\mu_2}\neq0$, $(-1)^{\gamma_2}$ is constant, and $\eta_\sigma>0$ for all $\sigma$. This completes the proof of \Cref{lem:gamma_map_nonzero}, and hence \Cref{prop:surj_onto_homs}. \qed

\section{The ideal of relations for the tilting algebra of $\Gr(n,2)$}
\label{chap:relnsgrn2}

In this section we identify the kernel of the $\kk$-algebra homomorphism $\Phi:\kk\TQ\twoheadrightarrow A$ from \Cref{thm:tilting_quiver_rank2}. This ideal describes the relations on the tilting quiver $\TQ$ from \Cref{defn:tilting_quiver_gr_n2}. Throughout, fix a basis $u_1,\dots,u_n$ of $V$. We will write angle brackets $\langle \ \rangle$ for linear subspaces and round brackets $( \ )$ for ideals.

\subsection{Strategy for finding $\ker(\Phi)$} Recall that for each $\lambda\in\TQ_0$, $e_\lambda\in\kk\TQ$ denotes the idempotent corresponding to the path of length zero at that vertex. Then for all pairs $\lambda< \mu\in\TQ_0$, we denote by $\Phi_{\lambda,\mu}$ the induced $\kk$-linear map obtained by restricting $\Phi$ to the subspace spanned by paths with tail at $\lambda$ and head at $\mu$, i.e.\
\[
\Phi_{\lambda,\mu} \colon e_\mu\kk\TQ e_\lambda \xrightarrowdbl{\phantom{aaa}} \Hom(\SS^\lambda\cW,\SS^\mu \cW)\cong \SS^{\mu/\lambda}V.
\]
Note that $\Phi_{\lambda,\mu}$ is surjective by \Cref{prop:surj_onto_homs}. Now, $\TQ$ is acyclic and  there are no relations involving paths of length one; indeed, relations only arise between paths that  share the same head and tail, and for all $a\in\TQ_1$ the only paths $p$ in $\TQ$ satisfying $\tl(a)=\tl(p)$ and $\hd(a)=\hd(p)$ are the arrows between the same two vertices, and these have linearly independent images under the map $\Phi$. It therefore suffices to find $\ker(\Phi_{\lambda,\mu})$ for every pair $(\lambda,\mu)$ in the set
\[
P:=\{(\lambda,\mu)\in{\TQ_0}^2 \mid \ \lambda<\mu, \modmu\geq\modlam+2\}.
\]
Denote by $K_{\lambda,\mu}$ a set of basis vectors for the subspace $\ker(\Phi_{\lambda,\mu})$. Then $\ker(\Phi)$ is the ideal generated by the union of these bases:
\begin{align}
\label{eqn:ker_phi_decomp}
\ker(\Phi)=\left(\bigcup_{(\lambda,\mu)\in P} K_{\lambda,\mu}\right).
\end{align}

We divide this section into two main steps. Define
\begin{align}
\label{eqn:P_2}
P_2:=\left\{(\lambda,\mu)\in{\TQ_0}^2 \mid \lambda<\mu,   \modmu=\modlam+2\right\}\subset P,
\end{align}
the pairs of vertices separated by paths of length two. The first step is to find $\ker(\Phi_{\lambda,\mu})$, and therefore $K_{\lambda,\mu}$, for all $(\lambda,\mu)\in P_2$. We find elements of $\ker(\Phi_{\lambda,\mu})$ by studying various compositions of the $f$-type and $g$-type maps defined in \Cref{prop:schurmaps} evaluated on an arbitrary section $w^\lambda$ of $\SS^\lambda\cW$. We then define the ideal
\begin{align}
\label{eqn:ideal_I_defn}
I:=\left( \bigcup_{(\lambda,\mu)\in P_2}K_{\lambda,\mu}\right)\subset \kk\TQ
\end{align}
generated by all relations of length two. Hence $I\subseteq \ker(\Phi)$, and by considering longer paths in $\kk\TQ$, in the second step we prove that $I=\ker(\Phi)$.

\begin{remark}
Note that in the case of $\Gr(4,2)$, the ideal $I$ was written down by Buchweitz, Leuschke and Van den Bergh in \cite[Example~8.4]{BLV15}. We recover this example in our discussion of $\Gr(5,2)$ in \Cref{exa:gr52}.
\end{remark}

\subsection{Relations between paths of length two}
\label{sec:gr2_relations}

In this section we identify the vector spaces $\ker(\Phi_{\lambda,\mu})$ for all $(\lambda,\mu)\in P_2$, and then extract bases $K_{\lambda,\mu}$ in order to define the ideal $I\subseteq\ker(\Phi)$ from \eqref{eqn:ideal_I_defn}. For each $(\lambda,\mu)\in P_2$, after finding a certain collection of relations we will perform a dimension count to prove these relations span $\ker(\Phi_{\lambda,\mu})$. To do this, we first need to decompose the codomain $\Hom(\SS^\lambda \cW, \SS^\mu \cW)\cong\SS^{\mu/\lambda}V$ into a sum of irreducibles.

It will be useful to recall the construction of $\Sym ^2 V$ and $\bw^2 V$ as subspaces embedded in $V^{\otimes 2}$. We have the quotients
\[
\Sym^2 V:= \frac{V^{\otimes 2}}{\langle v_1\otimes v_2 - v_2\otimes v_1 \rangle}, \quad \bw^2 V:= \frac{V^{\otimes 2}}{\langle v_1\otimes v_2 + v_2\otimes v_1 \rangle},
\]
and there is a natural isomorphism $V^{\otimes 2} \longrightarrow \Sym^2 V \oplus \bw^2 V$ given by
\begin{align}
\label{eqn:sym_wedge}
v_1\otimes v_2 \mapsto \frac{1}{2}(v_1\otimes v_2+v_2\otimes v_1, \ v_1\otimes v_2 - v_2\otimes v_1)=:(v_1v_2,v_1\ww v_2),
\end{align}
which enables us to identify $\Sym^2 V $ and $\bw^2 V$ as subspaces of $V^{\otimes 2}$.

\begin{proposition}
\label{prop:length2_maps}
Let $(\lambda,\mu)\in P_2$.
\begin{itemize}
	\setlength\itemsep{0em}
	\item[\one] If $\mu=(\lambda_1+2,\lambda_2)$ then $\Hom(\SS^\lambda \cW, \SS^\mu \cW)\cong\Sym^2 V$.
	\item[\two] If $\mu=(\lambda_1,\lambda_2+2)$ then $\Hom(\SS^\lambda \cW, \SS^\mu \cW)\cong\Sym^2 V$.
	\item[\three] If $\lambda_1=\lambda_2$ and $\mu=(\lambda_1+1,\lambda_1+1)$ then $\Hom(\SS^\lambda \cW, \SS^\mu \cW)\cong\bw^2 V $.
	\item[\four] If $\lambda_1>\lambda_2$ and $\mu=(\lambda_1+1,\lambda_2+1)$ then $\Hom(\SS^\lambda \cW, \SS^\mu \cW)\cong V^{\otimes 2}$.
\end{itemize}
\begin{proof}
For \one \ and \two \ this is just \Cref{cor:length1arrows_and_syms}. For \three \ and \four \ fix $\mu=(\lambda_1+1,\lambda_2+1)$, and using the notation of \Cref{lem:skew_hom_decomp}, let $\Gamma_{\mu/\lambda}$ be the set of partitions $\gamma$ corresponding to the irreducible summands of $\Hom(\SS^\lambda \cW, \SS^\mu \cW)\cong\SS^{\mu/\lambda}V$. By \Cref{lem:skew_hom_decomp}, we have $\gamma\in\Gamma_{\mu/\lambda}$ if and only if $\modgamma=2$ and $\lambda_2+1-\lambda_1\leq \gamma_2\leq 1$. If $\lambda_1=\lambda_2$ then the only $\gamma\in\Gamma_{\mu/\lambda}$ is $(1,1)$, giving  $\Hom(\SS^\lambda \cW, \SS^\mu \cW)\cong \SS^{(1,1)}V=\bw^2 V$ as required. If $\lambda_1>\lambda_2$ then both $(2,0),(1,1)\in\Gamma_{\mu/\lambda}$, hence $\Hom(\SS^\lambda \cW, \SS^\mu \cW)\cong\Sym^2 V\oplus\bw^2 V\cong V^{\otimes 2}$.
\end{proof}
\end{proposition}

We now describe $\ker(\Phi_{\lambda,\mu})$ in each of the four cases of \Cref{prop:length2_maps} respectively. To do this, we will find relations by evaluating elements of $\im(\Phi_{\lambda,\mu})$ on an arbitrary section $w^\lambda$ of $\SS^\lambda \cW$. Such elements are given by compositions of the $f$-type and $g$-type maps defined in \Cref{prop:schurmaps}.

\begin{remark}
\label{rem:first_iso}
In the first three cases there is only one route from $\lambda$ to $\mu$, so the domain of $\Phi_{\lambda,\mu}$ satisfies $e_\mu\kk\TQ e_\lambda\cong V^{\otimes 2}$. Since $\Phi_{\lambda,\mu}$ is surjective and $V^{\otimes 2}\cong\Sym^2 V\oplus\bw^2 V$, we just need to find relations that span a space isomorphic to the irreducible summand in the decomposition of $V^{\otimes 2}$ that is complement to the summand given by $\Hom(\SS^\lambda \cW, \SS^\mu \cW)$ in \Cref{prop:length2_maps}; these relations must then span $\ker(\Phi_{\lambda,\mu})$. Case \four \ is slightly different and we deal with that in \Cref{subsec:square}.
\end{remark}

\begin{notation}
\label{not:abuse}
\begin{itemize}	\setlength\itemsep{0em}
\item[\one] Let $w^\lambda$ be a section of $\SS^\lambda \cW$ as in \eqref{eqn:w_lambda}. In the following calculations, the $x_{1,1} \ww x_{1,2} \otimes \cdots \otimes x_{\lambda_2,1} \ww x_{\lambda_2,2}$ part of $w^\lambda$ is never altered so we simplify the notation by denoting $\underline{x}:=x_{1,1} \ww x_{1,2} \otimes \cdots \otimes x_{\lambda_2,1} \ww x_{\lambda_2,2}$. Thus we write
\[
w^\lambda=\underline{x}\otimes y_1\cdots y_{\lambda_1-\lambda_2}.
\]
\item[\two] When writing down the subspaces $\ker(\Phi_{\lambda,\mu})$, we will abuse notation slightly by writing $f^\lambda_{u_\rho}$ and $g^\lambda_{u_\rho}$ for the arrows $a^{\lambda,1}_\rho$ and $a^{\lambda,2}_\rho$ respectively. Henceforth the symbols $f^\lambda_{u_\rho},g^\lambda_{u_\rho}$ are juxtaposed to denote paths in $\kk\TQ$, but separated by $\circ$ for homomorphisms in $A$.
\end{itemize}
\end{notation}

\subsubsection{Paths of two horizontal arrows}
\label{subsec:horiz}
Here $\mu=(\lambda_1+2,\lambda_2)$ we consider paths of the form $$\SS^{(\lambda_1,\lambda_2)}\cW \longrightarrow \SS^{(\lambda_1+1,\lambda_2)}\cW \longrightarrow \SS^{(\lambda_1+2,\lambda_2)}\cW.$$ Denote $\nu=(\lambda_1+1,\lambda_2)$. Then for all $v_1,v_2 \in \{u_1,\dots,u_n\}$ we have
\begin{align*}
f^\nu_{v_2} \circ f^\lambda_{v_1}(w^\lambda)&=f^\nu_{v_2}\left( \underline{x} \otimes y_1\cdots y_{\lambda_1-\lambda_2} z_{v_1} \right ) \\
&=\underline{x} \otimes y_1\cdots y_{\lambda_1-\lambda_2} z_{v_1} z_{v_2} \\ &=\underline{x} \otimes y_1\cdots y_{\lambda_1-\lambda_2} z_{v_2} z_{v_1} \\&=f^\nu_{v_1}\left( \underline{x} \otimes y_1\cdots y_{\lambda_1-\lambda_2} z_{v_2} \right) \\ &=f^\nu_{v_1} \circ f^\lambda_{v_2}(w^\lambda)
\end{align*}
and so $f^\nu_{v_2}\circ f^\lambda_{v_1}-f^\nu_{v_1}\circ f^\lambda_{v_2}=0$. Using \eqref{eqn:sym_wedge} and identifying the tensor product with composition of maps, we may write down an isomorphism
\[
\bw^2 V \xrightarrow{\phantom{aa}\large{\cong}\phantom{aa}} \left\langle f^\nu_{v_2}\circ f^\lambda_{v_1}-f^\nu_{v_1}\circ f^\lambda_{v_2} \mid v_1,v_2 \in \{u_1,\dots,u_n\} \right\rangle
\]
where $v_1\ww v_2\mapsto f^\nu_{v_2}\circ f^\lambda_{v_1}-f^\nu_{v_1}\circ f^\lambda_{v_2}$. We have $\Hom(\SS^\lambda \cW,  \SS^{(\lambda_1+2,\lambda_2)} \cW)\cong\Sym^2 V$ by \Cref{prop:length2_maps}\one; hence, it follows from \Cref{rem:first_iso} that this is the entire subspace of relations because the domain of $\Phi_{\lambda,\mu}$ satisfies $e_\mu\kk\TQ e_\lambda \cong V^{\otimes 2}\cong\Sym^2 V\oplus\bw^2 V$. Thus we may conclude
\[
\ker(\Phi_{\lambda,\mu})=\left\langle f^\nu_{v_2} f^\lambda_{v_1}-f^\nu_{v_1} f^\lambda_{v_2} \mid v_1,v_2\in \{u_1,\dots,u_n\} \right\rangle.
\]
Note that in the original notation used in \Cref{defn:tilting_quiver_gr_n2}, this subspace is given by  $\ker(\Phi_{\lambda,\mu})=\left\langle a^{\nu,1}_{\rho_1} a^{\lambda,1}_{\rho_2}-a^{\nu,1}_{\rho_2} a^{\lambda,1}_{\rho_1} \mid 1\leq \rho_1,\rho_2\leq n \right\rangle$.

\subsubsection{Paths of two vertical arrows and paths between vertices on the diagonal}
\label{subsec:vert}
These two cases are similar to the first case; one must simply check that the desired relations between the $f$-type and $g$-type maps hold. We merely state the results here; for a full proof, see \cite[Subsections~5.1.2,~5.1.3]{Gre18}

As for paths of two vertical arrows, let $\mu=(\lambda_1,\lambda_2+2)$ and consider paths of the form $\SS^{(\lambda_1,\lambda_2)}\cW \rightarrow \SS^{(\lambda_1,\lambda_2+1)}\cW \rightarrow \SS^{(\lambda_1,\lambda_2+2)}\cW$. Denote $\nu=(\lambda_1,\lambda_2+1)$. Then we have
\[
\ker(\Phi_{\lambda,\mu})=\left\langle g^\nu_{v_2} g^\lambda_{v_1}-g^\nu_{v_1} g^\lambda_{v_2} \mid   v_1,v_2\in \{u_1,\dots,u_n\} \right\rangle.
\]

As for paths between vertices on the diagonal, suppose $\lambda_1=\lambda_2$ and let $\mu=(\lambda_1+1,\lambda_1+1)$. Consider paths of the form $\SS^{(\lambda_1,\lambda_1)}\cW \rightarrow \SS^{(\lambda_1+1,\lambda_1)}\cW \rightarrow \SS^{(\lambda_1+1,\lambda_1+1)}\cW$ and denote $\nu=(\lambda_1+1,\lambda_1)$. Then we have 
\[
\ker(\Phi_{\lambda,\mu})=\left\langle g^\nu_{v_2} f^\lambda_{v_1}+g^\nu_{v_1} f^\lambda_{v_2} \mid   v_1,v_2\in \{u_1,\dots,u_n\} \right\rangle.
\]

\subsubsection{Paths around a square}
\label{subsec:square}
Now suppose $\lambda_1>\lambda_2$ and $\mu=(\lambda_1+1,\lambda_2+1)$, and consider paths around a square as in \Cref{fig:maps_square_full}. In this case there are two routes from $\lambda$ to $\mu$, so we have $e_\mu \kk\TQ e_\lambda\cong V^{\otimes 2}\oplus V^{\otimes 2}$. Akin to \Cref{rem:first_iso}, surjectivity of $\Phi_{\lambda,\mu}$ and counting dimensions implies that $\ker(\Phi_{\lambda,\mu})\cong V^{\otimes 2}$ since $\Hom(\SS^\lambda\cW,\SS^\mu\cW)\cong V^{\otimes 2}$ by \Cref{prop:length2_maps}\four. Thus, we are looking for relations that span a space isomorphic to $ V^{\otimes 2}$.
\begin{figure}[!ht]
	\begin{center}
		\begin{tikzpicture}[xscale=0.9,yscale=0.9]
		\tikzset{edge/.style = {->,> = latex'}}
		
		\node[thick] (A) at  (0,0) {\begin{tabular}{c}$(\bw^2 \cW)^{\otimes \lambda_2}$ \\ $\otimes \Sym^{\lambda_1-\lambda_2} \cW$ \end{tabular}};
		\node[thick] (B) at  (5,0) {\begin{tabular}{c}$(\bw^2 \cW)^{\otimes \lambda_2}$ \\ $\otimes \Sym^{\lambda_1-\lambda_2+1} \cW$ \end{tabular}};
		\node[thick] (C) at  (0,4) {\begin{tabular}{c}$(\bw^2 \cW)^{\otimes \lambda_2+1}$ \\ $\otimes \Sym^{\lambda_1-\lambda_2-1} \cW$ \end{tabular}};
		\node[thick] (D) at  (5,4) {\begin{tabular}{c}$(\bw^2 \cW)^{\otimes \lambda_2+1}$ \\ $\otimes \Sym^{\lambda_1-\lambda_2} \cW$ \end{tabular}};
		
		\draw[edge] (A) -- node[above]{\textcolor{black}{$f^\lambda_\bullet$}} (B);
		\draw[edge] (A) -- node[right]{\textcolor{black}{$g^\lambda_\bullet$}} (C);
		\draw[edge] (B) -- node[right]{\textcolor{black}{$g^\nu_\bullet$}} (D);
		\draw[edge] (C) -- 
		node[above]{\textcolor{black}{$f^\delta_\bullet$}} (D);
		\draw [dashed, ->] (A) -- (D) node[midway,fill=white]{$V^{\otimes 2}$};
		
		\end{tikzpicture}
	\end{center}
	\caption[$\Gr(n,2)$ case: relations of paths around a square]{Paths around a square, with space of relations isomorphic to $V^{\otimes 2}$, generated by \eqref{eqn:square_relns}. Each $\bullet$ represents a basis vector $u_i, 1\leq i \leq n$, that varies independently.}
	\label{fig:maps_square_full}	
\end{figure}

\begin{lemma}
\label{lem:square_relns}
Denote $\nu=(\lambda_1+1,\lambda_2)$ and $\delta=(\lambda_1,\lambda_2+1)$. Then for all $v_1,v_2\in \{u_1,\dots,u_n\}$, we have
\begin{align}
\label{eqn:square_relns}
\left( \lambda_1-\lambda_2 \right)g^\nu_{v_2}\circ f^\lambda_{v_1}=\left( \lambda_1-\lambda_2+1\right)f^\delta_{v_1}\circ g^\lambda_{v_2} - f^\delta_{v_2}\circ g^\lambda_{v_1}.
\end{align}
\begin{proof}
Starting with the right hand side, we have
\begin{align*}
&\left( \lambda_1-\lambda_2+1\right)  f^\delta_{v_1}\circ g^\lambda_{v_2}(w^\lambda) - f^\delta_{v_2}\circ g^\lambda_{v_1}(w^\lambda)\\&=\left( \lambda_1-\lambda_2+1\right)f^\delta_{v_1}\left( \sum_{k=1}^{\lambda_1-\lambda_2} \underline{x} \otimes y_k \ww z_{v_2} \otimes \prod_{j\neq k} y_j \right)  -f^\delta_{v_2}\left( \sum_{k=1}^{\lambda_1-\lambda_2} \underline{x} \otimes y_k \ww z_{v_1} \otimes \prod_{j\neq k} y_j \right) \\&=\left( \lambda_1-\lambda_2+1\right) \sum_{k=1}^{\lambda_1-\lambda_2} \underline{x} \otimes y_k \ww z_{v_2} \otimes z_{v_1} \prod_{j\neq k} y_j - \sum_{k=1}^{\lambda_1-\lambda_2} \underline{x} \otimes y_k \ww z_{v_1} \otimes z_{v_2} \prod_{j\neq k} y_j.
\end{align*}
The left hand side becomes
\begin{align*}
&\left( \lambda_1-\lambda_2\right)g^\nu_{v_2} \circ f^\lambda_{v_1}(w^\lambda)\\&=\left( \lambda_1-\lambda_2\right)g^\nu_{v_2}\left( \underline{x} \otimes y_1\cdots y_{\lambda_1-\lambda_2} z_{v_1} \right) \\&=\left( \lambda_1-\lambda_2\right)\sum_{k=1}^{\lambda_1-\lambda_2} \left( \underline{x} \otimes y_k \ww z_{v_2} \otimes z_{v_1} \prod_{j\neq k} y_j\right)+\left( \lambda_1-\lambda_2\right) \underline{x} \otimes z_{v_1} \ww z_{v_2} \otimes  \prod_{j=1}^{\lambda_1-\lambda_2} y_j.
\end{align*}
Now subtract the left hand side from the right hand side to give
\begin{align*}
&\left( \lambda_1-\lambda_2+1\right)f^\delta_{v_1}\circ g^\lambda_{v_2} - f^\delta_{v_2}\circ g^\lambda_{v_1} - \left( \lambda_1-\lambda_2\right)g^\nu_{v_2}\circ f^\lambda_{v_1}\\
&=\sum_{k=1}^{\lambda_1-\lambda_2} \underline{x} \otimes y_k \ww z_{v_2} \otimes z_{v_1} \prod_{j\neq k} y_j- \sum_{k=1}^{\lambda_1-\lambda_2} \underline{x} \otimes y_k \ww z_{v_1} \otimes z_{v_2} \prod_{j\neq k} y_j\\&\quad -\left( \lambda_1-\lambda_2\right) \underline{x} \otimes z_{v_1} \ww z_{v_2} \otimes  \prod_{j=1}^{\lambda_1-\lambda_2} y_j \quad=:(\dagger).
\end{align*}
We now perform slightly different exchanges with each of the $\lambda_1-\lambda_2$ copies of $\underline{x} \otimes z_{v_1} \ww z_{v_2} \otimes  \prod_{j=1}^{\lambda_1-\lambda_2} y_j$ defining the last term of $(\dagger)$. With the first copy, perform an (E1)-type exchange by moving into $y_1$ the column containing $z_{v_1} \ww z_{v_2}$, yielding
\begin{align*}
\underline{x} \otimes z_{v_1} \ww z_{v_2} \otimes  \prod_{j=1}^{\lambda_1-\lambda_2} y_j&=\underline{x} \otimes y_1 \ww z_{v_2} \otimes z_{v_1} \prod_{j\neq 1} y_j+\underline{x} \otimes z_{v_1} \ww y_1 \otimes z_{v_2} \prod_{i\neq 1} y_i\\&=-\underline{x} \otimes z_{v_2} \ww y_1 \otimes z_{v_1} \prod_{j\neq 1} y_j+\underline{x} \otimes z_{v_1} \ww y_1 \otimes z_{v_2} \prod_{j\neq 1} y_j
\end{align*}
Now perform a similar exchange with the second copy using $y_2$, and in general with the $k$-th copy using $y_k$. Adding these all together, we get
\begin{align*}
&\left( \lambda_1-\lambda_2\right) \underline{x} \otimes z_{v_1} \ww z_{v_2} \otimes  \prod_{j=1}^{\lambda_1-\lambda_2} y_j\\&=- \sum_{k=1}^{\lambda_1-\lambda_2} \underline{x} \otimes z_{v_2} \ww y_k \otimes z_{v_1} \prod_{j\neq k} y_j+\sum_{k=1}^{\lambda_1-\lambda_2} \underline{x} \otimes z_{v_1} \ww y_k \otimes z_{v_2} \prod_{j\neq k} y_j.
\end{align*}
Finish by substituting the right-hand side of this identity into the last term of $(\dagger)$ to get zero.
\end{proof}
\end{lemma}

Since
\[
V^{\otimes 2}\cong \left\langle \left( \lambda_1-\lambda_2 \right)g^\nu_{v_2}\circ f^\lambda_{v_1}-\left( \lambda_1-\lambda_2+1\right)f^\delta_{v_1}\circ g^\lambda_{v_2} + f^\delta_{v_2}\circ g^\lambda_{v_1} \mid v_1,v_2\in \cB \right\rangle,
\]
the discussion prior to \Cref{lem:square_relns} implies 
\[
\ker(\Phi_{\lambda,\mu})=\left\langle \left( \lambda_1-\lambda_2 \right)g^\nu_{v_2} f^\lambda_{v_1}-\left( \lambda_1-\lambda_2+1\right)f^\delta_{v_1} g^\lambda_{v_2} + f^\delta_{v_2} g^\lambda_{v_1} \mid v_1,v_2\in \cB \right\rangle.
\]
\begin{remark}
\label{rem:square_either_way}
Consequently, because each path going in one direction around the square may be written as a linear combination of paths going the opposite way, there are no relations amongst paths traversing in the same direction.
\end{remark}

\begin{proposition}
\label{prop:relns_2_ideal}
For each $(\lambda,\mu)\in P_2$, let $\nu,\delta$ be the vertices that lie on paths $\lambda$ and $\mu$ as defined in Sections~\ref{subsec:horiz}-\ref{subsec:square}. Define the sets $K_{\lambda,\mu}$ as follows:
\begin{itemize} \setlength\itemsep{0em}
	\item[\one] if $\mu=(\lambda_1+2,\lambda_2)$,   $K_{\lambda,\mu}=\left\{ f^\nu_{u_j} f^\lambda_{u_i}-f^\nu_{u_i} f^\lambda_{u_j} \mid   1\leq i,j\leq n \right\}$.
	\item[\two] if $\mu=(\lambda_1,\lambda_2+2)$,  $K_{\lambda,\mu}=\left\{ g^\nu_{u_j} g^\lambda_{u_i}-g^\nu_{u_i} g^\lambda_{u_j} \mid 1\leq i,j\leq n \right\}$.
	\item[\three] if $\lambda_1=\lambda_2$ and $\mu=(\lambda_1+1,\lambda_1+1)$,  $K_{\lambda,\mu}=\left\{ g^\nu_{u_j} f^\lambda_{u_i}+g^\nu_{u_i} f^\lambda_{u_j} \mid 1\leq i,j\leq n \right\}$.
	\item[\four] if $\lambda_1>\lambda_2$ and  $\mu=(\lambda_1+1,\lambda_2+1)$,  
	\[
	K_{\lambda,\mu}=\left\{ \left( \lambda_1-\lambda_2 \right)g^\nu_{u_j} f^\lambda_{u_i}-\left( \lambda_1-\lambda_2+1\right)f^\delta_{u_i} g^\lambda_{u_j} + f^\delta_{u_j} g^\lambda_{u_i} \mid 1\leq i,j\leq n \right\}.
	\] 
\end{itemize}
Then each $K_{\lambda,\mu}$ is a basis of $\ker(\Phi_{\lambda,\mu})$ and the ideal
\[
I:=\left( \bigcup_{(\lambda,\mu)\in P_2}K_{\lambda,\mu}\right)\subseteq\ker(\Phi)
\]
contains all of the relations in $\TQ$ generated by paths of length two.

\begin{proof}
Subsections~\ref{subsec:horiz}-\ref{subsec:square}.
\end{proof}
\end{proposition}

\subsection{Relations between longer paths}
\label{sec:relations_longer}

Having described the ideal $I\subseteq\ker(\Phi)$ in \Cref{prop:relns_2_ideal}, we now prove that $I=\ker(\Phi)$. Recall from equation \eqref{eqn:ker_phi_decomp} that we have $\ker(\Phi)=(\cup_{(\lambda,\mu)\in P} K_{\lambda,\mu})$. Since the ideal $I$ is taken over $P_2$, a subset of $P$,  we must show that $\ker(\Phi_{\lambda,\mu})=e_\mu I e_\lambda$ for all pairs $(\lambda,\mu)$ in the complement of $P_2$, i.e.\ the set
\[
P_l:=\left\{(\lambda,\mu)\in{\TQ_0}^2 \mid \lambda<\mu,   \modmu>\modlam+2\right\}
\]
where we write $l$ simply to mean `longer'.

We do this in two propositions. First we compute $\ker(\Phi_{\lambda,\mu})$ in the special cases that all paths $\lambda\rightarrow\mu$ are straight lines in $\TQ$, i.e.\ $(\lambda,\mu)\in P_l$ where $\mu$ is of the form either $(\mu_1,\lambda_2)$ or $(\lambda_1,\mu_2)$. Afterwards we deal with the remaining cases where both $\lambda_1<\mu_1$ and $\lambda_2<\mu_2$.

First recall that the $k$-th symmetric power of $V$ is given by
\begin{align}
\label{eqn:sym_power_defn}
\Sym^k V:= \frac{V^{\otimes k}}{\langle v_1\otimes\cdots \otimes v_{k}-v_{\sigma(1)}\otimes\cdots \otimes v_{\sigma(k)} \mid \sigma\in S_k \rangle},
\end{align}
where $S_k$ is the permutation group on $\{1,\dots,k\}$. Following \cite[$\S1.1.1$~p.3]{Wey03}, there is a natural embedding given by
\begin{align}
\label{eqn:sym_embedding}
\begin{split}
\Delta_k\colon \Sym^k V &\lhook\joinrel\longrightarrow V^{\otimes k} \\ v_1\cdots v_{k} &\xmapsto{\phantom{aaa}}  \frac{1}{k!}\sum_{\sigma \in S_{k}} v_{\sigma(1)}\otimes\cdots \otimes v_{\sigma(k)}.
\end{split}
\end{align}

\begin{notation}
\label{not:abuse_again}
We will make use of \Cref{not:abuse}\two \ again: when writing down the subspaces $\ker(\Phi_{\lambda,\mu})$, we will abuse notation slightly by writing $f^\lambda_{u_\rho}$ and $g^\lambda_{u_\rho}$ for the arrows $a^{\lambda,1}_\rho$ and $a^{\lambda,2}_\rho$ respectively. Henceforth the symbols $f^\lambda_{u_\rho},g^\lambda_{u_\rho}$ are juxtaposed to denote paths in $\kk\TQ$, but separated by $\circ$ for homomorphisms in $A$.

In addition, we will reuse some notation from \Cref{sec:surj_proof}; denote $m_1=\mu_1-\lambda_1$ and $m_2=\mu_2-\lambda_2$, and since we are using the $f$ and $g$ notation for arrows in $\TQ_1$, without ambiguity we may drop the superscript on all arrows in a path except the first as in \Cref{not:superscript_drop}\one.
\end{notation}

\begin{proposition}
\label{prop:symmetric_relations}
Let $(\lambda,\mu)\in P_l$ and suppose $\mu$ is of the form either $(\mu_1,\lambda_2)$ or $(\lambda_1,\mu_2)$. Then
\begin{align}
\label{eqn:main_sym_relns}
\begin{split}
\ker(\Phi_{\lambda,\mu})=
\begin{cases}
\left\langle f_{v_{m_1}}\cdots f^\lambda_{v_{1}}-f_{v_{\sigma(m_1)}}\cdots f^\lambda_{v_{\sigma(1)}} \mid \sigma\in S_{m_1}, v_i\in \{u_1,\dots,u_n\}\right\rangle & \ \text{if} \  \mu=(\mu_1,\lambda_2), \\
\left\langle g_{v_{m_2}}\cdots g^\lambda_{v_{1}}-g_{v_{\sigma(m_2)}}\cdots g^\lambda_{v_{\sigma(1)}} \mid  \sigma\in S_{m_2}, v_i\in \{u_1,\dots,u_n\} \right\rangle & \ \text{if} \  \mu=(\lambda_1,\mu_2).
\end{cases}
\end{split}
\end{align}
In particular, $\ker(\Phi_{\lambda,\mu})=e_\mu I e_\lambda$.
\begin{proof}
The fact that the subspaces \eqref{eqn:main_sym_relns} are contained in $\ker(\Phi_{\lambda,\mu})$ follows from straightforward induction arguments on the results of Sections~\ref{subsec:horiz}-\ref{subsec:vert}. We claim that these relations span $\ker(\Phi_{\lambda,\mu})$ by dimension count. The domain of the surjective map $\Phi_{\lambda,\mu}$ is $e_\mu \kk\TQ e_\lambda$, which is isomorphic to either $V^{\otimes m_1}$ or $V^{\otimes m_2}$ when $\mu$ is equal to $(\mu_1,\lambda_2)$ or $(\lambda_1,\mu_2)$ respectively. Taking the quotient of these spaces by the appropriate subspace in \eqref{eqn:main_sym_relns} gives $\Sym^{m_1} V$ or $\Sym^{m_2} V$ respectively by \eqref{eqn:sym_power_defn}. Since the codomain of $\Phi_{\lambda,\mu}$ is $\Hom(\SS^\lambda \cW,\SS^\mu\cW)$, which is isomorphic to $\Sym^{m_1}V$ or $\Sym^{m_2}V$ respectively by \Cref{cor:length1arrows_and_syms}, the claim follows from the first isomorphism theorem.

For the final statement, first suppose $\mu=(\mu_1,\lambda_2)$. The subspace $e_\mu I e_\lambda$ consists only of the relations amongst straight line paths $\lambda\rightarrow\mu$, and all of these are generated by those in \Cref{prop:relns_2_ideal}\one, specifically those of the form $\{f_{v_2} f^{\gamma_{i}}_{v_1} - f_{v_1} f^{\gamma_{i}}_{v_2} \mid v_1,v_2\in\{u_1,\dots,u_n\} \}$ where $\gamma_i=\lambda+ie_1$ for $0\leq i\leq m_1-2$. Hence, define $S^T_{m_1}\subseteq S_{m_1}$ to be the subset of adjacent transpositions, i.e.\ $\sigma\in S^T_{m_1}$ if for some $1\leq k\leq m_1-1$ we have $\sigma(k)=k+1, \ \sigma(k+1)=k$, and $\sigma(j)=j$ for all $j\neq k,k+1$. Then
\[
e_\mu I e_\lambda =\left\langle f_{v_{m_1}}\cdots f^\lambda_{v_{1}}-f_{v_{\sigma(m_1)}}\cdots f^\lambda_{v_{\sigma(1)}} \mid \sigma\in S^T_{m_1}, v_i\in \{u_1,\dots,u_n\} \right\rangle.
\]
Hence we have $e_\mu I e_\lambda\subseteq \ker(\Phi_{\lambda,\mu})$, but since $S_{m_1}$ is generated by the elements of $S^T_{m_1}$ the reverse inclusion follows simply by performing a sequence of permutations in $S^T_{m_1}$. The proof is similar for $\mu=(\lambda_1,\mu_2)$. 
\end{proof}
\end{proposition}

\begin{proposition}
\label{prop:no_more_relns}
Let $(\lambda,\mu)\in P_l$ and suppose $\lambda_1<\mu_1$ and $\lambda_2<\mu_2$. Then $\ker(\Phi_{\lambda,\mu})=e_\mu I e_\lambda$.
\begin{proof}
First of all, as in \Cref{sec:surj_proof}, it is enough to consider the special case that $\lambda_2=0$ as each of the spaces $e_\mu \kk\TQ e_\lambda, e_\mu I e_\lambda$ and $\Hom(\SS^\lambda \cW, \SS^\mu \cW)$ is unchanged if we replace $(\lambda_1,\lambda_2)$ and $(\mu_1,\mu_2)$ by $(\lambda_1-\lambda_2,\lambda_2-\lambda_2)$ and $(\mu_1-\lambda_2,\mu_2-\lambda_2)$ respectively; see \Cref{cor:invariant_homs}.

There is a commutative diagram
\[
\begin{tikzpicture}[xscale=1,yscale=1]
\node (A) at (0,2.5) {$e_\mu\kk\TQ e_\lambda$};
\node (B) at (6,2.5) {$\Hom(\SS^\lambda\cW,\SS^\mu\cW)$};
\node (C) at (0,0) {$\displaystyle\frac{e_\mu\kk\TQ e_\lambda}{e_\mu I e_\lambda}$};
\draw [->>] (A) -- node [above]{$\Phi_{\lambda,\mu}$} (B);
\draw [->>] (A) -- node [left]{$\pi$} (C);
\draw [->>] (C) -- node [below right]{$\Psi_{\lambda,\mu}$} (B);
\end{tikzpicture}
\]
where $\pi$ is the quotient map. The goal is to show that $\Psi_{\lambda,\mu}$ is injective: then $\ker(\Psi_{\lambda,\mu})=0$ and so $\ker(\Phi_{\lambda,\mu})=\pi^{-1}(0)=e_\mu I e_\lambda$ as required.

Consider a path (or more generally, a linear combination of paths) $p\in e_\mu\kk\TQ e_\lambda$ and let $\nu$ be the vertex $(\mu_1,\lambda_2)$. If at any point on the path(s) $p$ there is a vertical arrow immediately before a horizontal arrow, it is possible to use relations from $I$, namely those of \Cref{prop:relns_2_ideal}\four, to rewrite those two arrows as a linear combination of arrows around the same square in $\TQ$ that instead go horizontally before vertically. We can repeat this process until $p$ has been  rewritten completely as linear combination of paths that all go strictly horizontally before vertically; in other words, there exists an element $p_2\otimes_{\kk} p_1\in e_\mu\kk\TQ e_\nu\otimes_{\kk} e_\nu\kk\TQ e_\lambda$ such that $[p]=[p_2p_1]\in e_\mu\kk\TQ e_\lambda/e_\mu I e_\lambda$.

Now, in $\kk\TQ/I$ we have $[p_2p_1]=[p_2][p_1]$ where
\[
[p_2]\in \frac{e_\mu \kk\TQ e_\nu}{e_\mu I e_\nu}\cong \Sym^{m_2} V, \quad [p_1]\in \frac{e_\nu \kk\TQ e_\lambda}{e_\nu I e_\lambda}\cong \Sym^{m_1} V,
\]
and we have used the isomorphisms from  \Cref{prop:symmetric_relations} with $m_1=\mu_1-\lambda_1$ and $m_2=\mu_2$. Since every  $[p]\in e_\mu\kk\TQ e_\lambda/e_\mu I e_\lambda$ can be written in the form $[p_2p_1]$, there exists a surjective homomorphism
\[
\xi_1\colon \Sym^{m_2}V\otimes_{\kk}\Sym^{m_1}V \xrightarrowdbl{\phantom{aaa}} \frac{e_\mu \kk\TQ e_\lambda}{e_\mu I e_\lambda}
\]
where $[p_2]\otimes_{\kk}[p_1]\mapsto [p_2p_1]=[p]$. We now split into two subcases according to the decomposition of $\Hom(\SS^\lambda \cW, \SS^\mu \cW)$ into irreducibles.

\smallskip
\noindent\textbf{Case (i): $\mu_2\leq \lambda_1$}. In this case $\Hom(\SS^\lambda \cW, \SS^\mu \cW)$ is isomorphic to $\Sym^{m_2}V\otimes_{\kk}\Sym^{m_1}V$; this follows from \Cref{lem:skew_hom_decomp} and the first Pieri rule (\Cref{prop:pieri}). We have the diagram
\[
\begin{tikzpicture}[xscale=1,yscale=1]
\node (A) at (0,2) {$\Sym^{m_2}V\otimes_{\kk}\Sym^{m_1}V$};
\node (B) at (0,0) {$\displaystyle\frac{e_\mu\kk\TQ e_\lambda}{e_\mu I e_\lambda}$};
\node (C) at (5,2) {$\Hom(\SS^\lambda \cW, \SS^\mu \cW)$};
\node (D) at (5,0) {$\Hom(\SS^\lambda \cW, \SS^\mu \cW)$};
\draw [->>] (A) -- node [left]{$\xi_1$} (B);
\draw [->>] (B) -- node [above]{$\Psi_{\lambda,\mu}$} (D);
\draw [->] (C) -- node [above]{$\cong$} (A);
\draw [->] (C) -- node [right]{$\text{id}$} (D);
\end{tikzpicture}
\]
and therefore $\xi_1$ and $\Psi_{\lambda,\mu}$ must be bijections. In particular, $\Psi_{\lambda,\mu}$ is injective as required.

\smallskip
Note that the relations between paths on the main diagonal of $\TQ$, i.e.\ those of \Cref{prop:relns_2_ideal}\three, are absent in case \one. Indeed, any vertex $\gamma$ on a path $\lambda\rightarrow\mu$ satisfies $0\leq \gamma_2\leq \mu_2$ and $\lambda_1\leq \gamma_1\leq \mu_1$. Vertices on the diagonal also satisfy $\gamma_1=\gamma_2$, and when $\mu_2\leq \lambda_1$ the only such possible vertex is $(\lambda_1,\mu_2)$. Thus, with at most one vertex on the diagonal on any path $\lambda\rightarrow\mu$,  \Cref{prop:relns_2_ideal}\three \ plays no role in $e_\mu\kk\TQ e_\lambda/e_\mu I e_\lambda$ in this case. The next case is different.

\smallskip
\noindent\textbf{Case (ii): $\mu_2>\lambda_1$}. The irreducible summands of $\Hom(\SS^\lambda \cW, \SS^\mu \cW)$ form a proper subset of those in the irreducible decomposition of $\Sym^{m_2}V\otimes_{\kk}\Sym^{m_1}V$; again this follows from \Cref{lem:skew_hom_decomp} and the first Pieri rule. Indeed, as hinted above we must now also consider the possibility of relations between vertices along the diagonal. Let $d=\mu_2-\lambda_1>0$. Then the vertices $(\lambda_1+k,\lambda_1+k)$ for all $0\leq k\leq d$ may appear on paths $\lambda\rightarrow \mu$. Previously, we used the relations around squares to rewrite $p$ as a linear combination of paths going strictly horizontally before vertically. While we may also do that here and the surjective map $\xi_1$ still applies, we can also use the relations around squares to rewrite $p$ as a linear combination of paths that take the route
\[
(\lambda_1,0)\rightarrow (\lambda_1,\lambda_1)\rightarrow (\lambda_1+1,\lambda_1+1)\rightarrow \cdots \rightarrow (\lambda_1+d,\lambda_1+d)=(\mu_2,\mu_2)\rightarrow(\mu_1,\mu_2),
\]
in other words, paths in $p$ travel vertically from $\lambda$ to the diagonal and then staircase along it as much as possible, exiting horizontally towards $\mu$ at height $\mu_2$. Define the sequence $\nu_0=\lambda, \nu_i=(\lambda_1+i-1,\lambda_1+i-1)$ for all $1\leq i \leq d+1$, and $\nu_{d+2}=\mu$. Then there exists an element $q_2 \otimes_{\kk} z_d \otimes_{\kk} \cdots \otimes_{\kk} z_1 \otimes_{\kk} q_1$ in
\begin{align}
\label{eqn:diagonal_space}
\bigotimes_{i=0}^{d+1} e_{\nu_{i+1}} \kk\TQ e_{\nu_i}
\end{align}
such that $[p]=[q_2z_d\cdots z_1 q_1]\in e_\mu\kk\TQ e_\lambda/e_\mu I e_\lambda$.

Following the strategy above, we take the quotient of each subspace in \eqref{eqn:diagonal_space} by the appropriate graded slice of $I$. For $i=0$ we have straight vertical paths $(\lambda_1,0)\rightarrow (\lambda_1,\lambda_1)$ and for $i=d+1$ we have straight horizontal paths $(\mu_2,\mu_2)\rightarrow (\mu_1,\mu_2)$. Hence, using \Cref{prop:symmetric_relations} we have
\[
[q_2]\in \frac{e_\mu \kk\TQ e_{\nu_{d+1}}}{e_\mu I e_{\nu_{d+1}}}\cong \Sym^{\mu_1-\mu_2} V, \quad [q_1]\in \frac{e_{\nu_1} \kk\TQ e_\lambda}{e_{\nu_1} I e_\lambda}\cong \Sym^{\lambda_1} V.
\]
Each of the remaining subspaces, $e_{\nu_{i+1}} \kk\TQ e_{\nu_i}$ for $1\leq i \leq d$, is spanned by paths starting at a vertex on the diagonal and going horizontally then vertically to the next vertex on the diagonal. These paths were considered in \Cref{subsec:vert} and therefore we have
\[
[z_i]\in \frac{e_{\nu_{i+1}} \kk\TQ e_{\nu_i}}{e_{\nu_{i+1}} I e_{\nu_i}}\cong \bw^2 V, \quad 1\leq i \leq d.
\]
Consider the quotient
\[
D:=\bigotimes_{i=0}^{d+1} \frac{e_{\nu_{i+1}} \kk\TQ e_{\nu_i}}{e_{\nu_{i+1}} I e_{\nu_i}}.
\]
Then by the prior discussion, there is a surjective homomorphism
\[
\xi_2\colon D\cong \Sym^{\mu_1-\mu_2}V \otimes_{\kk} \left(\bw^2 V\right)^{\otimes d} \otimes_{\kk} \Sym^{\lambda_1} V \xrightarrowdbl{\phantom{aaa}} \frac{e_\mu \kk\TQ e_\lambda}{e_\mu I e_\lambda}
\]
where $[q_2]\otimes_{\kk} [z_d] \otimes_{\kk} \cdots \otimes_{\kk} [z_1] \otimes_{\kk} [q_1] \mapsto [q_2z_d\cdots z_1q_1]=[p]$.

We must now find the irreducible decomposition of $D$, which we accomplish using the Pieri rules; see \Cref{prop:pieri}. We first decompose the central collection of terms $(\bw^2 V)^{\otimes d}$. The second Pieri rule states that tensoring a Schur power $\SS^\gamma V$ by $\bw^2 V$ yields a direct sum taken over the ways of adding two new boxes to distinct rows of $\gamma$. Starting with $\SS^{(1,1)}V=\bw^2 V$ and applying this rule $d-1$ times, we have
\[
\left(\bw^2 V\right)^{\otimes d} \cong \SS^{(d,d)}V \oplus X,
\]
where $X$ is a direct sum of Schur powers of $V$ defined by Young diagrams with at least three rows.

Next we tensor $\SS^{(d,d)}V \oplus X$ by the $\Sym^{\lambda_1}V$ term. By the first Pieri rule this has decomposition given over the ways of adding $\lambda_1$ boxes to each partition in the sum $\SS^{(d,d)}V \oplus X$ with no two in the same column. By adding $\lambda_1$ boxes to the top row of $(d,d)$ we have the term $\SS^{(\lambda_1+d,d)}V=\SS^{(\mu_2,\mu_2-\lambda_1)}V$, but again, every other term has at least three rows in this second decomposition. Amending $X$ to $X'$ (we don't care exactly what these terms with at least three rows are), we have
\[
\left(\bw^2 V\right)^{\otimes d}\otimes_{\kk}\Sym^{\lambda_1}V \cong \SS^{(\mu_2,\mu_2-\lambda_1)}V\oplus X'.
\]

Lastly, we must tensor this decomposition by $\Sym^{\mu_1-\mu_2}V$. Using the first Pieri rule again and focusing only on the terms that will produce Young diagrams with at most two rows, we have $D\cong (\bigoplus_\gamma \SS^\gamma V ) \oplus X''$ where $\gamma$ ranges over the partitions $(\max\{m_1,m_2\},\min\{m_1,m_2\}),(\max\{m_1,m_2\}+1,\min\{m_1,m_2\}-1),\dots, (\mu_1,\mu_2-\lambda_1)$. By \Cref{lem:skew_hom_decomp} these partitions are precisely those that describe the irreducible decomposition of $\Hom(\SS^\lambda \cW, \SS^\mu \cW)$. Therefore,
\[
D\cong \Hom(\SS^\lambda \cW, \SS^\mu \cW) \oplus X''.
\]

Hence, we have a diagram of surjective maps
\[
\begin{tikzpicture}[xscale=1,yscale=1]
\node (A) at (0,3) {$\Sym^{m_2}V\otimes_{\kk} \Sym^{m_1}V$};
\node (B) at (0,0) {$D\cong \Hom(\SS^\lambda \cW, \SS^\mu \cW) \oplus X''$};
\node (C) at (4,1.5) {$\displaystyle\frac{e_\mu\kk\TQ e_\lambda}{e_\mu I e_\lambda}$};
\node (D) at (9,1.5) {$\Hom(\SS^\lambda \cW, \SS^\mu \cW)$};
\draw [->>] (A) -- node [above]{$\xi_1$} (C);
\draw [->>] (B) -- node [above]{$\xi_2$} (C);
\draw [->>] (C) -- node [above]{$\Psi_{\lambda,\mu}$} (D);
\end{tikzpicture}
\]
Since $\xi_1$ and $\xi_2$ are surjective, $e_\mu\kk\TQ e_\lambda/e_\mu I e_\lambda$ must be isomorphic to a subspace of the direct sum of the summands that appear in both the irreducible decompositions of $\Sym^{m_2}V\otimes_{\kk} \Sym^{m_1}V$ and $D$. Since the decomposition of $\Sym^{m_2}V\otimes_{\kk} \Sym^{m_1}V$ consists of only partitions with at most two rows, of which those comprising $\Hom(\SS^\lambda \cW, \SS^\mu \cW)$ form a proper subset, we conclude that $e_\mu\kk\TQ e_\lambda/e_\mu I e_\lambda$ is isomorphic to a subspace of $\Hom(\SS^\lambda \cW, \SS^\mu \cW)$. This forces the surjective map $\Psi_{\lambda,\mu}$ to be an isomorphism, and in particular is injective as required.
\end{proof}
\end{proposition}

We now have a full presentation of Kapranov's tilting algebra for $\Gr(n,2)$.

\begin{theorem}
\label{thm:complete_relns_rank2}
For $\Gr(n,2)$ let $E$ be the tilting bundle  \eqref{eqn:tilt_bundle_grass} and let $A=\End(E)$. Let $\TQ$ be the quiver defined in \Cref{defn:tilting_quiver_gr_n2}. Then the $\kk$-algebra $A$ is isomorphic to $\kk\TQ/I$, where
\[
I=\left(\bigcup_{(\lambda,\mu)\in P_2} K_{\lambda,\mu}\right)
\]
and
\begin{itemize}
	\item[\one] if $\mu=(\lambda_1+2,\lambda_2)$,   $K_{\lambda,\mu}=\left\{ f_{u_j} f^\lambda_{u_i}-f_{u_i} f^\lambda_{u_j} \mid   1\leq i,j\leq n \right\}$.
	\item[\two] if $\mu=(\lambda_1,\lambda_2+2)$,  $K_{\lambda,\mu}=\left\{ g_{u_j} g^\lambda_{u_i}-g_{u_i} g^\lambda_{u_j} \mid 1\leq i,j\leq n \right\}$.
	\item[\three] if $\lambda_1=\lambda_2$ and $\mu=(\lambda_1+1,\lambda_1+1)$,  $K_{\lambda,\mu}=\left\{ g_{u_j} f^\lambda_{u_i}+g_{u_i} f^\lambda_{u_j} \mid 1\leq i,j\leq n \right\}$.
	\item[\four] if $\lambda_1>\lambda_2$ and  $\mu=(\lambda_1+1,\lambda_2+1)$,  
	\[
	K_{\lambda,\mu}=\left\{ \left( \lambda_1-\lambda_2 \right)g_{u_j} f^\lambda_{u_i}-\left( \lambda_1-\lambda_2+1\right)f_{u_i} g^\lambda_{u_j} + f_{u_j} g^\lambda_{u_i} \mid 1\leq i,j\leq n \right\}.
	\] 
\end{itemize}
\begin{proof}
In \Cref{chap:tiltingalgebra} we defined a $\kk$-algebra homomorphism $\Phi\colon \kk\TQ\rightarrow A$ and proved it is surjective. After establishing that $\ker(\Phi)=(\bigcup_{(\lambda,\mu)\in P} K_{\lambda,\mu})$, we presented the ideal  $I=(\bigcup_{(\lambda,\mu)\in P_2} K_{\lambda,\mu})\subseteq \ker(\Phi)$ in \Cref{prop:relns_2_ideal}. Propositions~\ref{prop:symmetric_relations}~and~\ref{prop:no_more_relns} then prove that $\langle K_{\lambda,\mu}\rangle = \ker(\Phi_{\lambda,\mu})= e_\mu I e_\lambda\subseteq I$ for all $(\lambda,\mu)\in P \setminus P_2$. This completes the proof that $\ker(\Phi)=I$.
\end{proof}
\end{theorem}

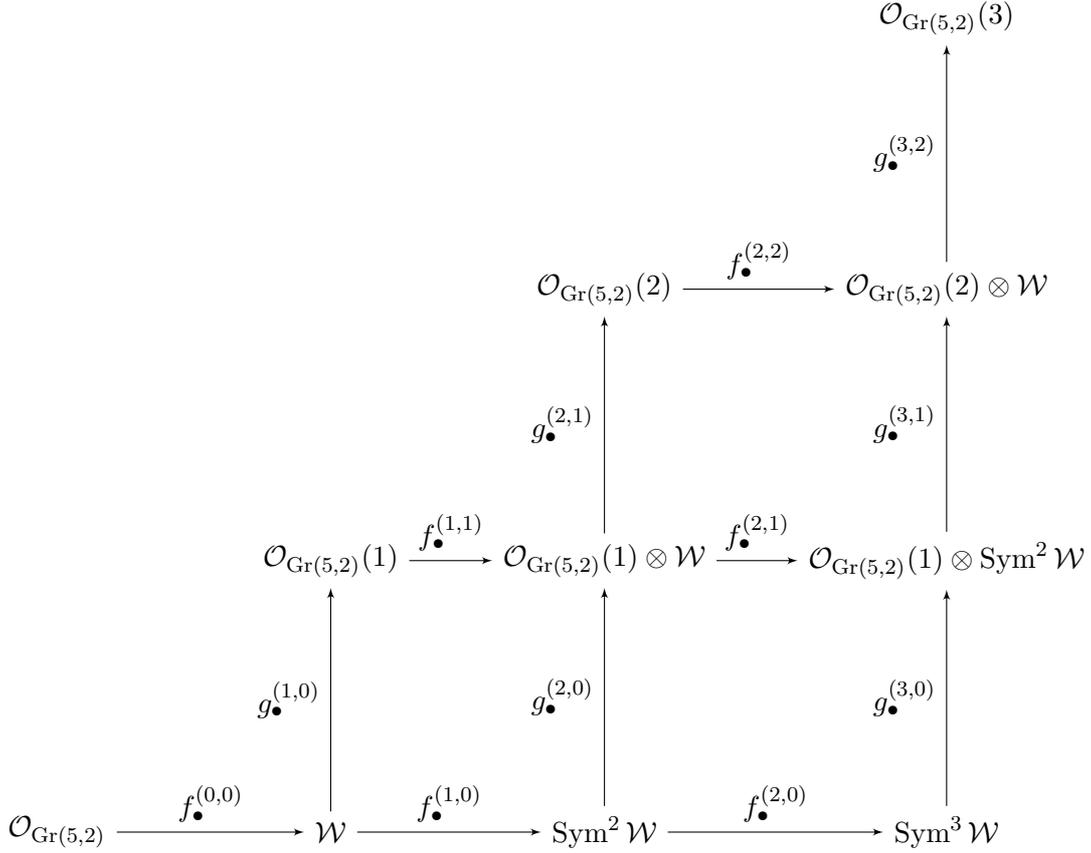
\begin{figure}[!ht]
	\begin{center}
		\begin{tikzpicture}[xscale=0.9,yscale=0.9]
		\tikzset{edge/.style = {->,> = latex'}}
		\node[thick] (A) at  (0,0) {$\cO_{\Gr(5,2)}$};
		\node[thick] (B) at  (4,0) {$\cW$};
		\node[thick] (C) at  (4,4) {$\cO_{\Gr(5,2)}(1)$};
		\node[thick] (D) at  (8,0) {$\Sym^2 \cW$};
		\node[thick] (E) at  (8,4) {$\cO_{\Gr(5,2)}(1) \otimes \cW$};
		\node[thick] (F) at  (8,8) {$\cO_{\Gr(5,2)}(2)$};
		\node[thick] (G) at  (13,0) {$\Sym^3 \cW$};
		\node[thick] (H) at  (13,4) {$\cO_{\Gr(5,2)}(1)\otimes \Sym^2 \cW$};
		\node[thick] (I) at  (13,8) {$\cO_{\Gr(5,2)}(2)\otimes \cW$};
		\node[thick] (J) at  (13,12) {$\cO_{\Gr(5,2)}(3)$};
		
		\draw[edge] (A) -- node[above]{\textcolor{black}{$f^{(0,0)}_\bullet$}} (B);
		\draw[edge] (B) -- node[above]{\textcolor{black}{$f^{(1,0)}_\bullet$}} (D);
		\draw[edge] (C) -- node[above]{\textcolor{black}{$f^{(1,1)}_\bullet$}} (E);
		\draw[edge] (B) -- node[left]{\textcolor{black}{$g^{(1,0)}_\bullet$}} (C);
		\draw[edge] (D) -- node[left]{\textcolor{black}{$g^{(2,0)}_\bullet$}} (E);
		\draw[edge] (E) -- node[left]{\textcolor{black}{$g^{(2,1)}_\bullet$}} (F);

		\draw[edge] (D) -- node[above]{\textcolor{black}{$f^{(2,0)}_\bullet$}} (G);
		\draw[edge] (E) -- node[above]{\textcolor{black}{$f^{(2,1)}_\bullet$}} (H);
		\draw[edge] (F) -- node[above]{\textcolor{black}{$f^{(2,2)}_\bullet$}} (I);
		\draw[edge] (G) -- node[left]{\textcolor{black}{$g^{(3,0)}_\bullet$}} (H);
		\draw[edge] (H) -- node[left]{\textcolor{black}{$g^{(3,1)}_\bullet$}} (I);
		\draw[edge] (I) -- node[left]{\textcolor{black}{$g^{(3,2)}_\bullet$}} (J);
		\end{tikzpicture}
	\end{center}
	\caption[The tilting quiver for $\Gr(5,2)$]{The tilting quiver for $\Gr(5,2)$. Each arrow represents $5$ arrows corresponding to the basis $u_1,\dots,u_5$ of $V$.}
	\label{fig:gr4-2_more_relns}
\end{figure}

\begin{example}
\label{exa:gr52}
Consider the case where $n=5$. Then for $\Gr(5,2)$, the tilting quiver is given by \Cref{fig:gr4-2_more_relns}, and below we list the relations that span $I=\ker(\Phi)$. Following \Cref{not:abuse_again}, we only require a superscript for the first arrow in a path since the $f$ and $g$-type notation determines the remaining arrows. For all $1\leq i,j\leq 5$, we have the following.
\begin{itemize}
\item Horizontal paths: for $\lambda=(0,0),(1,0),(1,1)$ we have
\[
f_{u_j} f^\lambda_{u_i}=f_{u_i} f^\lambda_{u_j}.
\]
\item Vertical paths: for $\lambda=(2,0),(3,0),(3,1)$ we have
\[
g_{u_j} g^\lambda_{u_i}=g_{u_i} g^\lambda_{u_j}.
\]
\item Paths on the main diagonal: for $\lambda=(0,0),(1,1),(2,2)$ we have
\[
g_{u_j} f^\lambda_{u_i}=-g_{u_i} f^\lambda_{u_j}.
\]
\item Lower-left and upper squares: for $\lambda=(1,0),(2,1)$ we have
\[
g_{u_j} f^\lambda_{u_i}=2f_{u_i} g^\lambda_{u_j}-f_{u_j} g^\lambda_{u_i}.
\]
\item Lower-right square: for $\lambda=(2,0)$ we have
\[
2g_{u_j} f^\lambda_{u_i}=3f_{u_i} g^\lambda_{u_j}-f_{u_j} g^\lambda_{u_i}.
\]
\end{itemize}
\end{example}

\begin{remark}
\label{rem:BLV_example}
Following \Cref{exa:gr52}, consider the full sub-quiver of $\TQ$ for $\Gr(5,2)$ defined by deleting the vertices $(3,0),\dots,(3,3)\in\TQ_0$ and any arrows with head or tail at those vertices. By also removing all arrows associated to $u_5$, we recover the tilting quiver for $\Gr(4,2)$; see \Cref{fig:gr4-2}. In particular, the relations defining $\ker(\Phi)$ for $\Gr(4,2)$ form a sublist of those in the above example; these were calculated by Buchweitz, Leuschke and Van den Bergh in \cite[Example~8.4]{BLV15}.
\end{remark}

\section{Reconstructing $\Gr(n,2)$ from Kapranov's tilting bundle}
\label{chap:moduli_grn2}

Recall that Kapranov's tilting bundle for $\Gr(n,2)$ is given by
\begin{align}
\label{eqn:tilt_bundle_grass2}
E=\bigoplus_{\lambda\in \Young(n-2,2)} \SS^\lambda \cW
\end{align}
where $\cW$ is the tautological quotient bundle.
Having given an explicit presentation of $A=\End(E)$ in \Cref{thm:complete_relns_rank2}, we now apply this result to reconstruct $\Gr(n,2)$ from $E$.

We begin with the construction of the \emph{multigraded linear series} from \cite{Cra11}, an example of the moduli spaces originally constructed by King~\cite{Kin94}. Let $E_0,\dots,E_m$ be the indecomposable summands of $E$ from \eqref{eqn:tilt_bundle_grass2} with $E_0=\cO_{\Gr(n,2)}$. Denote by $\vv:=(v_j)\in \NN^{m+1}$ the dimension vector satisfying $v_j:= \rank(E_j)$ for all $0\leq j\leq m$. To introduce our choice of stability condition $\theta$, first set $\theta^\prime=(-\sum_{i=1}^m v_i,1,1,\dots,1)\in \Hom(\ZZ^{m+1},\QQ)$. An $A$-module $M=\bigoplus_{0\leq j\leq m} M_j$ of dimension vector $\vv$ is $\theta^\prime$-stable if and only if $M$ is generated as an $A$-module by any nonzero element of $M_0$; any such stability condition is called \emph{0-generated}. Since $\vv$ is indivisible, King~\cite[Proposition~5.3]{Kin94} constructs the fine moduli space $\mathcal{M}(A,\vv,\theta^\prime)$ of isomorphism classes of $\theta^\prime$-stable $A$-modules of dimension vector $\vv$ as a GIT quotient. In particular, $\mathcal{M}(A,\vv,\theta^\prime)$ comes with an ample bundle $\mathcal{O}(1)$. Let $k\geq 1$ denote the smallest positive integer such that $\mathcal{O}(k)$ is very ample. Then $\theta:=k\theta'$ is also a 0-generated stability condition, and we write
\[
\modulie:=\moduli
\]
for the fine moduli space of $\theta$-stable $A$-modules of dimension vector $\vv$. The universal property of $\modulie$ determines a morphism
\begin{equation}
\label{eqn:closed_immersion}
f_E\colon \Gr(n,2)\longrightarrow \modulie
\end{equation}
and moreover we have the following.

\begin{proposition}
\label{thm:closed_immersion}
The morphism $f_E\colon \Gr(n,2)\rightarrow \modulie$ is a closed immersion.
\begin{proof}
The tautological quotient bundle $\cW$ is globally generated by \cite[Corollary~3.5(i)]{Cra11}, and hence so is every indecomposable summand $E_i$ of $E$ from \eqref{eqn:tilt_bundle_grass2}. Therefore $\det(E_i)$ is also globally generated for each $E_i$. By \cite[Lemma~3.7]{Cra11} the line bundle $\det(\cW)$ is ample, hence the line bundle $\bigotimes_{0\leq j\leq m} \det(E_j)$, which is a tensor product of $\det(\cW)$ and other globally generated bundles, is also ample by \cite[Exercise~II.7.5(d)]{Har77}. The result now follows from \cite[Theorem~2.6]{CIK17}.
\end{proof}
\end{proposition}

\begin{theorem}
\label{thm:grn2_isomorphism}
The closed immersion $f_E:\Gr(n,2)\rightarrow \modulie$ from \eqref{eqn:closed_immersion} is an isomorphism.
\begin{proof}
We follow a similar strategy used by Craw and the author in the main result of \cite{CG18}. First, in \Cref{subsec:part1} we give an explicit description of points in the image of $f_E$. Then, in \Cref{subsec:part2} we take an arbitrary point $w\in\modulie$ and find a point $y\in\Gr(n,2)$ such that, modulo the change of basis group action, $f_E(y)$ is equivalent to $w$. This implies $f_E$ is surjective on closed points, and so by \Cref{thm:closed_immersion}, $f_E$ is an isomorphism of $\Gr(n,2)$ onto the underlying reduced scheme of $\modulie$. The proof of Bergman-Proudfoot~\cite[Theorem~1.13]{BP08} then implies that $f_E$ identifies the tangent space at each closed point $y$ of $\Gr(n,2)$ with the tangent space of $\modulie$ at $f_E(y)$. Since $f_E$ is surjective on closed points, it follows that $\modulie$ is reduced and therefore $f_E$ is an isomorphism of schemes. 
\end{proof}
\end{theorem}

\subsection{Describing the image of $f_E$.}
\label{subsec:part1}
Fix a closed point $y\in\Gr(n,2)$ and a basis $u_i$, $1\leq i \leq n$, of $V$. Since closed points of $\Gr(n,2)$ are one-to-one with surjective linear maps $\kk^n\rightarrow \kk^2$, we may represent $y$ as a full rank $2\times n$ matrix. Without loss of generality we assume that the first two columns are linearly independent; then we can use the group action on the quiver flag variety to change basis such that
\[
y=\begin{pmatrix}
1 & 0 & x_1 & x_3 & \cdots & x_{2n-7} & x_{2n-5}\\ 
0 & 1 & x_2 & x_4 & \cdots &x_{2n-6} & x_{2n-4}
\end{pmatrix}
\]
for some scalars $x_1,\dots,x_{2n-4}\in\kk$. We will now see how the closed point $f_E(y)\in\modulie$ is determined by only these scalars and the choice of maps given in \Cref{prop:schurmaps}, but first we will need to introduce some notation to describe closed points of $\modulie$. Since $\modulie$ is also a quiver flag variety, we have a matrix for each arrow in the tilting quiver; see \Cref{fig:grn2tiltQ}, remembering that each arrow in the figure represents $n$ arrows corresponding to a basis of $V$. Since each of these arrows corresponds to a map defined in \Cref{prop:schurmaps}, we denote the matrix of the arrow $f^\lambda_{u_i}$ by $F^\lambda_i$ and  $g^\lambda_{u_i}$ by $G^\lambda_i$. Observe that these matrices must satisfy the matrix relations corresponding to those in \Cref{thm:complete_relns_rank2}.

As seen in \eqref{eqn:sections}, the basis $u_i$ gives us a basis of sections $z_{u_i}$ of $H^0(\Gr(n,2),\cW)$ which hereafter we simply denote $z_i$. Then for any point $v\in \Gr(n,2)$ there exists an open set $U_{i,j}\subset \Gr(n,2)$, $1\leq i<j\leq n$, such that $z_i(v)$ and $z_j(v)$ form a basis of the fibre $\cW_v$. For our fixed point $y$ with matrix above we will reorder the basis elements if necessary so that we have $y\in U_{1,2}$. We then write $b_1:=z_1(y)$ and $b_2:=z_2(y)$ for the basis of $\cW_y$ and for $3\leq i \leq n$ we have $z_{i}(y)=x_{2i-5}b_1+x_{2i-4}b_2$. Now, the basis $b_1,b_2$ of $\cW_y$ induces a basis on each fibre $(\SS^\lambda \cW)_y=\SS^\lambda \cW_y$ for all $\lambda\in \Young(n-2,2)$ in the canonical way; for example, $\bw^2 \cW_y$ has basis $b_2\ww b_1$ and $\Sym^2 \cW_y$ has basis $b_1b_1, b_1b_2, b_2b_2$. In general, we have $\dim(\SS^\lambda\cW_y)=\lambda_1-\lambda_2+1$ and the induced basis on the fibre $\SS^\lambda\cW_y$ is given by
\[
p^\lambda_j:=\left(b_2\ww b_1\right)^{\otimes \lambda_2} \otimes b_1^{\lambda_1-\lambda_2-j} b_2^j \ , \ 0\leq j\leq \lambda_1-\lambda_2.
\]

To describe $f_E(y)$ we therefore must find all of the $F^\lambda_i$ and $G^\lambda_i$ matrices. We deal with the $F^\lambda_i$ matrices first. Fix $\lambda$ such that $\lambda+e_1\in \Young(n-2,2)$. We will only calculate $F^\lambda_3$ here; for the other $F^\lambda_i$ simply substitute scalars from the $i$-th column of $y$ and then the calculation is the same. The induced basis of $\SS^{\lambda+e_1}\cW_y$ is given by
\[
q^\lambda_j:=\left(b_2\ww b_1\right)^{\otimes \lambda_2} \otimes b_1^{\lambda_1-\lambda_2+1-j} b_2^j \ , \ 0\leq j\leq \lambda_1-\lambda_2+1,
\]
and for all $0\leq j \leq \lambda_1-\lambda_2$ we have
\begin{align*}
(f_{u_3}^\lambda)_y(p^\lambda_j)&=x_1 \left(b_2\ww b_1\right)^{\otimes \lambda_2} \otimes b_1^{\lambda_1-\lambda_2-j+1} b_2^j + x_2\left(b_2\ww b_1\right)^{\otimes \lambda_2} \otimes b_1^{\lambda_1-\lambda_2-j} b_2^{j+1}\\
&=x_1q^\lambda_j+x_2q^\lambda_{j+1}
\end{align*}
which yields the following $(\lambda_1-\lambda_2+2)\times (\lambda_1-\lambda_2+1)$ matrix for $F_3^{\lambda}$,
\begin{align}
\label{eqn:Fmatrix}
F_3^{\lambda}=
\begin{pmatrix}
x_1 & 0 & 0 & \cdots & \cdots & 0\\ 
x_2 & x_1 & 0 & \cdots & \cdots & 0 \\
0 & x_2 & x_1 & \cdots & \cdots & 0 \\
\vdots & \ddots & \ddots & \ddots & \ddots & \vdots \\
0 & 0 & \cdots & \cdots & x_2 & x_1 \\
0 & 0 & \cdots & \cdots & 0 & x_2 \\ 
\end{pmatrix}.
\end{align}

Next we find the $G^\lambda_i$ matrices. Fix $\lambda$ such that $\lambda+e_2\in \Young(n-2,2)$. Like above we will calculate $G^\lambda_3$ only. The induced basis of $\SS^{\lambda+e_2}\cW_y$ is given by
\[
s^\lambda_j:=\left(b_2\ww b_1\right)^{\otimes \lambda_2+1} \otimes b_1^{\lambda_1-\lambda_2-1-j} b_2^j \ , \ 0\leq j\leq \lambda_1-\lambda_2-1,
\]
and for all $0\leq j \leq \lambda_1-\lambda_2$ we have
\begin{align*}
(g_{u_3}^\lambda)_y(p^\lambda_j)&=jx_1(b_2\ww b_1)^{\otimes \lambda_2+1}\otimes b_1^{\lambda_1-\lambda_2-j}b_2^{j-1}\\&\phantom{=}-(\lambda_1-\lambda_2-j)x_2 (b_2\ww b_1)^{\otimes \lambda_2+1}\otimes b_1^{\lambda_1-\lambda_2-j-1}b_2^j\\
&=\begin{cases}
-x_2(k-1)s^\lambda_0 & \text{if} \ j=0\\
jx_1s_{j-1}-(k-1-j)x_2s^\lambda_j&\text{if} \ 1\leq j \leq \lambda_1-\lambda_2-1\\
x_1(k-1)s^\lambda_{\lambda_1-\lambda_2-1} & \text{if} \ j=\lambda_1-\lambda_2\\
\end{cases}
\end{align*}
which yields the following $(\lambda_1-\lambda_2)\times (\lambda_1-\lambda_2+1)$ matrix for $G_3^{\lambda}$,
\begin{align}
\label{eqn:Gmatrix}
G_3^{\lambda}=
\begin{pmatrix}
-(k-1)x_2 & x_1 & 0 & \cdots & \cdots & \cdots & 0\\ 
0 & -(k-2)x_2 & 2x_1 & \cdots & \cdots & \cdots & 0 \\
\vdots & \ddots & \ddots & \ddots & \ddots & \ddots & \vdots  \\
0 & \cdots & \cdots & \cdots & -2x_2 & (k-2)x_1 & 0 \\
0 & \cdots & \cdots & \cdots & 0 & -x_2 & (k-1)x_1 \\ 
\end{pmatrix}.
\end{align}
where $k:=\lambda_1-\lambda_2+1$. As mentioned previously, for the other $F^\lambda_i$ and $G^\lambda_i$ matrices simply take $F^\lambda_3$ or $G^\lambda_3$ respectively and substitute $x_1,x_2$ for the entries in the $i$-th column of $y$.

\begin{example}
\label{exa:newgr42}
In the case $n=4$, which serves as our base case in the upcoming induction argument, we have
\[
y=\begin{pmatrix}
1&0&x_1&x_3\\
0&1&x_2&x_4
\end{pmatrix}
\]
and $f_E(y)$ is given by the system of matrices in \Cref{fig:gr4-2mats}. Observe that the four $2\times1$ matrices corresponding to the arrows $\cO_{\Gr(4,2)}\rightarrow \cW$ together produce naturally the $2\times4$ matrix $y$.

\begin{figure}[!ht]
	\begin{center}
		\begin{tikzpicture}[xscale=0.9,yscale=0.9]
		\tikzset{edge/.style = {->,> = latex'}}
		
		\node[thick] (A) at  (0,0) {$\cO_{\Gr(4,2)}$};
		\node[thick] (B) at  (5,0) {$\cW$};
		\node[thick] (C) at  (5,6) {$\bw^2 \cW$};
		\node[thick] (D) at  (13,0) {$\Sym^2 \cW$};
		\node[thick] (E) at  (13,6) {$\bw^2 \cW \otimes \cW$};
		\node[thick] (F) at  (13,12) {$(\bw^2 \cW)^{\otimes 2}$};
		
		\draw[edge] (A) -- node[below,yshift=-0.5em]{\textcolor{black}{\small{$\begin{pmatrix}
					1\\
					0\\
					\end{pmatrix}
					\begin{pmatrix}
					0\\
					1\\
					\end{pmatrix}
					\begin{pmatrix}
					x_1\\
					x_2\\
					\end{pmatrix}
					\begin{pmatrix}
					x_3\\
					x_4\\
					\end{pmatrix}$}}}(B);
		\draw[edge] (B) -- node[left,text width=2cm,align=right]{\textcolor{black}{\small{$\begin{pmatrix}
					0 & 1\\
					\end{pmatrix}$\\
					$\begin{pmatrix}
					-1 & 0\\
					\end{pmatrix}$\\
					$\begin{pmatrix}
					-x_2 & x_1\\
					\end{pmatrix}$\\
					$\begin{pmatrix}
					-x_4 & x_3\\
					\end{pmatrix}$}}}(C);			
		\draw[edge] (B) -- node[below,yshift=-0.5em]{\textcolor{black}{\small{$\begin{pmatrix}
					1 & 0\\
					0 & 1\\
					0 & 0\\
					\end{pmatrix}
					\begin{pmatrix}
					0 & 0\\
					1 & 0\\
					0 & 1\\
					\end{pmatrix}
					\begin{pmatrix}
					x_1 & 0\\
					x_2 & x_1\\
					0 & x_2\\
					\end{pmatrix}
					\begin{pmatrix}
					x_3 & 0\\
					x_4 & x_3\\
					0 & x_4\\
					\end{pmatrix}$}}}(D);		
		\draw[edge] (C) -- node[above]{\textcolor{black}{\small{$\begin{pmatrix}
					1\\
					0\\
					\end{pmatrix}
					\begin{pmatrix}
					0\\
					1\\
					\end{pmatrix}
					\begin{pmatrix}
					x_1\\
					x_2\\
					\end{pmatrix}
					\begin{pmatrix}
					x_3\\
					x_4\\
					\end{pmatrix}$}}}(E);
		\draw[edge] (D) -- node[left,text width=4.5cm,align=right]{\textcolor{black}{\small{$\begin{pmatrix}
					0 & 1 & 0\\
					0 & 0 & 2\\
					\end{pmatrix}$\\
					$\begin{pmatrix}
					-2 & 0 & 0\\
					0 & -1 & 0\\
					\end{pmatrix}$\\
					$\begin{pmatrix}
					-2x_2 & x_1 & 0\\
					0 & -x_2 & 2x_1\\
					\end{pmatrix}$\\
					$\begin{pmatrix}
					-2x_4 & x_3 & 0\\
					0 & -x_4 & 2x_3\\
					\end{pmatrix}$}}}(E);
		\draw[edge] (E) -- node[left,text width=2cm,align=right]{\textcolor{black}{\small{$\begin{pmatrix}
					0 & 1\\
					\end{pmatrix}$\\
					$\begin{pmatrix}
					-1 & 0\\
					\end{pmatrix}$\\
					$\begin{pmatrix}
					-x_2 & x_1\\
					\end{pmatrix}$\\
					$\begin{pmatrix}
					-x_4 & x_3\\
					\end{pmatrix}$}}}(F);	
		\end{tikzpicture}
	\end{center}
	\caption[Image of a point $y\in\Gr(4,2)$ under $f_E$]{The image of $y\in\Gr(4,2)$ under $f_E$. Matrices on horizontal or vertical arrows are of the form $F^\lambda_i$ or $G^\lambda_i$ respectively.}
	\label{fig:gr4-2mats}
\end{figure}
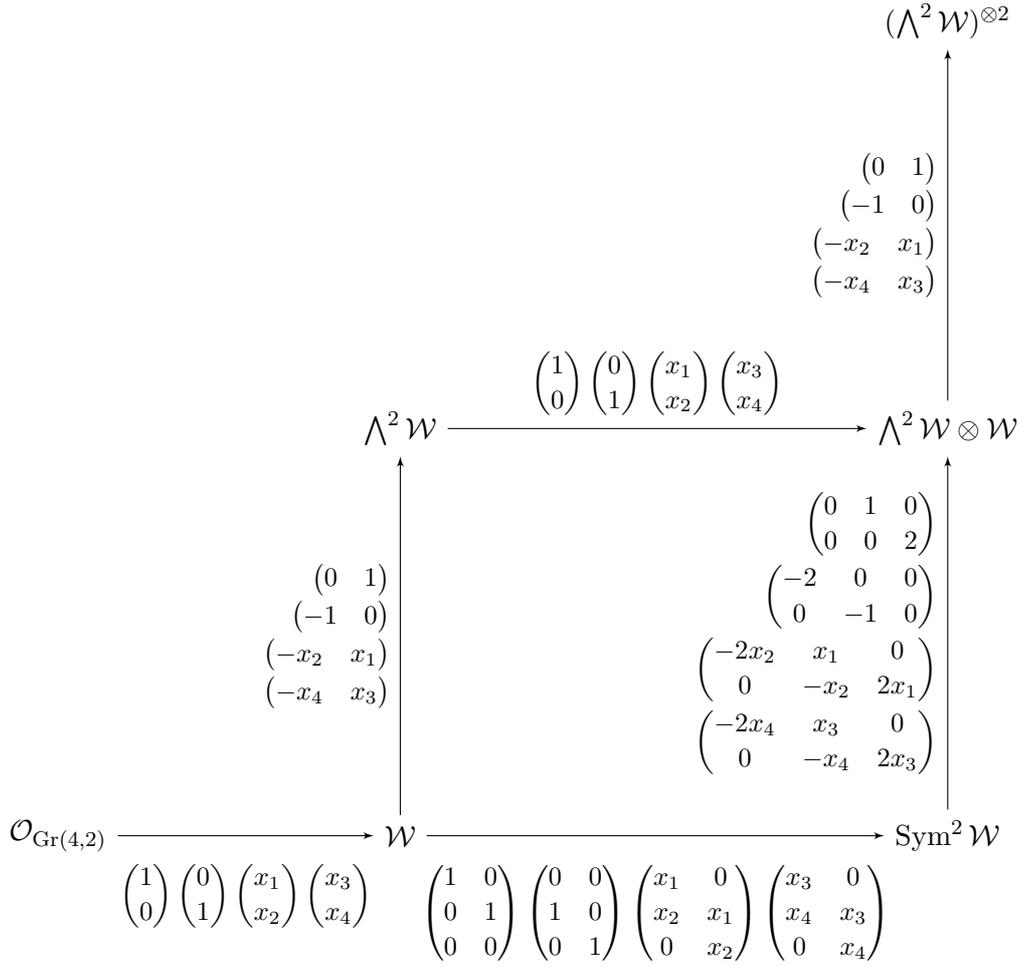
\end{example}

\subsection{Surjectivity of $f_E$ on closed points.}
\label{subsec:part2}

We prove surjectivity using \Cref{thm:complete_relns_rank2} and a technical induction argument. It will be easier to prove the base case, $\Gr(4,2)$, before actually stating the induction hypothesis.

\begin{figure}[!ht]
	\begin{center}
		\begin{tikzpicture}[xscale=0.9,yscale=0.9]
		\tikzset{edge/.style = {->,> = latex'}}
		\node[thick] (A) at  (0,0) {$\cO_{\Gr(4,2)}$};
		\node[thick] (B) at  (4,0) {$\cW$};
		\node[thick] (C) at  (4,4) {$\bw^2 \cW$};
		\node[thick] (D) at  (8,0) {$\Sym^2 \cW$};
		\node[thick] (E) at  (8,4) {$\bw^2 \cW \otimes \cW$};
		\node[thick] (F) at  (8,8) {$(\bw^2 \cW)^{\otimes 2}$};
		
		\draw[edge] (0.7,0.24) -- node[above]{\textcolor{black}{$F^{(0,0)}_\bullet$}} (3.5,0.24);
		\draw[edge] (0.7,0.08) -- (3.5,0.08);
		\draw[edge] (0.7,-0.08) -- (3.5,-0.08);
		\draw[edge] (0.7,-0.24) -- (3.5,-0.24);
		
		\draw[edge] (4.5,0.24) -- node[above]{\textcolor{black}{$F^{(1,0)}_\bullet$}} (7.1,0.24);
		\draw[edge] (4.5,0.08) -- (7.1,0.08);
		\draw[edge] (4.5,-0.08) -- (7.1,-0.08);
		\draw[edge] (4.5,-0.24) -- (7.1,-0.24);
		
		\draw[edge] (4.7,4.24) -- node[above]{\textcolor{black}{$F^{(1,1)}_\bullet$}} (6.8,4.24);
		\draw[edge] (4.7,4.08) -- (6.8,4.08);
		\draw[edge] (4.7,3.92) -- (6.8,3.92);
		\draw[edge] (4.7,3.76) -- (6.8,3.76);
		
		\draw[edge] (3.76,0.5) -- node[left]{\textcolor{black}{$G^{(1,0)}_\bullet$}} (3.76,3.5);
		\draw[edge] (3.92,0.5) -- (3.92,3.5);
		\draw[edge] (4.08,0.5) -- (4.08,3.5);
		\draw[edge] (4.24,0.5) -- (4.24,3.5);
		
		\draw[edge] (7.76,0.5) -- node[left]{\textcolor{black}{$G^{(2,0)}_\bullet$}} (7.76,3.5);
		\draw[edge] (7.92,0.5) -- (7.92,3.5);
		\draw[edge] (8.08,0.5) -- (8.08,3.5);
		\draw[edge] (8.24,0.5) -- (8.24,3.5);
		
		\draw[edge] (7.76,4.5) -- node[left]{\textcolor{black}{$G^{(2,1)}_\bullet$}} (7.76,7.5);
		\draw[edge] (7.92,4.5) -- (7.92,7.5);
		\draw[edge] (8.08,4.5) -- (8.08,7.5);
		\draw[edge] (8.24,4.5) -- (8.24,7.5);	
		\end{tikzpicture}
	\end{center}
	\caption[The tilting quiver for $\Gr(4,2)$]{The tilting quiver for $\Gr(4,2)$ with matrices representing the point  $w\in \modulie$. Each {\small$\bullet$} varies independently from $1$ to $4$.}
	\label{fig:gr4-2}
\end{figure}
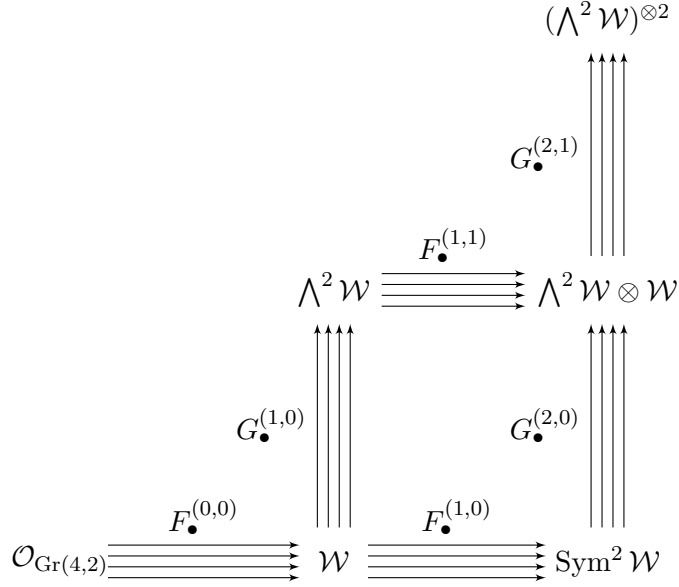

\smallskip
\noindent\textbf{Base case: $\Gr(4,2)$.}
Let $n=4$ and fix $w\in\modulie$.  Begin by fixing a basis $u_1,u_2,u_3,u_4$ of $V=\kk^4$. By \Cref{thm:complete_relns_rank2} we have the tilting quiver given in \Cref{fig:gr4-2} and for all $1\leq i,j \leq 4$ the  relations are:
\begin{align}
\begin{split}
\label{eqn:base_case_relns}
\bullet& \ G^{(1,0)}_{i} F^{(0,0)}_{j}+G^{(1,0)}_{j} F^{(0,0)}_{i},\\
\bullet& \ F^{(1,0)}_{i} F^{(0,0)}_{j}-F^{(1,0)}_{j} F^{(0,0)}_{i},\\
\bullet& \ G^{(2,0)}_{i} F^{(1,0)}_{j}-2F^{(1,1)}_{j} G^{(1,0)}_{i}+F^{(1,1)}_{i} G^{(1,0)}_{j},\\
\bullet& \ G^{(2,1)}_{i} F^{(1,1)}_{j}+G^{(2,1)}_{j} F^{(1,1)}_{i},\\
\bullet& \ G^{(2,1)}_{i} G^{(2,0)}_{j}-G^{(2,1)}_{j} G^{(2,0)}_{i}.
\end{split}
\end{align}

Besides the $F^\lambda_i,G^\lambda_i$ matrices satisfying these relations, $w$ must also be $\theta$-stable. We capture this notion as follows: as per \cite[Equation~2.2]{Cra11}, for any vertex $\lambda$ juxtapose all matrices whose corresponding arrows have head at $\lambda$ into a single matrix denoted $W_\lambda$. For example, $W_{(2,1)}$ is the $2\times 16$ matrix formed by concatenating the four $2\times1$ $F^{(1,1)}_i$ matrices together with the four $2\times 3$ $G^{(2,0)}_i$ matrices. Stability of $w$ is equivalent to each of the matrices $W_\lambda$ being full rank. Like in the previous section, we will relabel the basis $u_1,\dots,u_4$ if necessary to ensure that the first two columns of $W_{(1,1)}$ are linearly independent. Then, using the group action at the vertex $(1,1)$, $w$ takes the form given in \Cref{fig:gr4-2mats_general} for some unknown $y_1,\dots,y_{72}\in\kk$.

\begin{figure}[!ht]
	\begin{center}
		\begin{tikzpicture}[xscale=0.9,yscale=0.9]
		\tikzset{edge/.style = {->,> = latex'}}
		
		\node[thick] (A) at  (0,0) {$\cO_{\Gr(4,2)}$};
		\node[thick] (B) at  (5,0) {$\cW$};
		\node[thick] (C) at  (5,6) {$\bw^2 \cW$};
		\node[thick] (D) at  (13,0) {$\Sym^2 \cW$};
		\node[thick] (E) at  (13,6) {$\bw^2 \cW \otimes \cW$};
		\node[thick] (F) at  (13,12) {$(\bw^2 \cW)^{\otimes 2}$};
		
		\draw[edge] (A) -- node[below,yshift=-0.5em]{\textcolor{black}{\small{$\begin{pmatrix}
					1\\
					0\\
					\end{pmatrix}
					\begin{pmatrix}
					0\\
					1\\
					\end{pmatrix}
					\begin{pmatrix}
					x_1\\
					x_2\\
					\end{pmatrix}
					\begin{pmatrix}
					x_3\\
					x_4\\
					\end{pmatrix}$}}}(B);
		\draw[edge] (B) -- node[left,text width=1.7cm,align=left]{\textcolor{black}{\small{$\begin{pmatrix}
					y_1 & y_2\\
					\end{pmatrix}$\\
					$\begin{pmatrix}
					y_3 & y_4\\
					\end{pmatrix}$\\
					$\begin{pmatrix}
					y_5 & y_6\\
					\end{pmatrix}$\\
					$\begin{pmatrix}
					y_7 & y_8\\
					\end{pmatrix}$}}}(C);			
		\draw[edge] (B) -- node[below,xshift=0.5em,yshift=-0.5em]{\textcolor{black}{\small{$\begin{pmatrix}
					y_9 & y_{10}\\
					y_{11} & y_{12}\\
					y_{13} & y_{14}\\
					\end{pmatrix}
					\begin{pmatrix}
					y_{15} & y_{16}\\
					y_{17} & y_{18}\\
					y_{19} & y_{20}\\
					\end{pmatrix}
					\begin{pmatrix}
					y_{21} & y_{22}\\
					y_{23} & y_{24}\\
					y_{25} & y_{26}\\
					\end{pmatrix}
					\begin{pmatrix}
					y_{27} & y_{28}\\
					y_{29} & y_{30}\\
					y_{31} & y_{32}\\
					\end{pmatrix}$}}}(D);		
		\draw[edge] (C) -- node[above]{\textcolor{black}{\small{$\begin{pmatrix}
					y_{33}\\
					y_{34}\\
					\end{pmatrix}
					\begin{pmatrix}
					y_{35}\\
					y_{36}\\
					\end{pmatrix}
					\begin{pmatrix}
					y_{37}\\
					y_{38}\\
					\end{pmatrix}
					\begin{pmatrix}
					y_{39}\\
					y_{40}\\
					\end{pmatrix}$}}}(E);
		\draw[edge] (D) -- node[left,text width=3cm,align=left]{\textcolor{black}{\small{$\begin{pmatrix}
					y_{41} & y_{42} & y_{43}\\
					y_{44} & y_{45} & y_{46}\\
					\end{pmatrix}$\\
					$\begin{pmatrix}
					y_{47} & y_{48} & y_{49}\\
					y_{50} & y_{51} & y_{52}\\
					\end{pmatrix}$\\
					$\begin{pmatrix}
					y_{53} & y_{54} & y_{55}\\
					y_{56} & y_{57} & y_{58}\\
					\end{pmatrix}$\\
					$\begin{pmatrix}
					y_{59} & y_{60} & y_{61}\\
					y_{62} & y_{63} & y_{64}\\
					\end{pmatrix}$}}}(E);
		\draw[edge] (E) -- node[left,text width=2cm,align=left]{\textcolor{black}{\small{$\begin{pmatrix}
					y_{65} & y_{66}\\
					\end{pmatrix}$\\
					$\begin{pmatrix}
					y_{67} & y_{68}\\
					\end{pmatrix}$\\
					$\begin{pmatrix}
					y_{69} & y_{70}\\
					\end{pmatrix}$\\
					$\begin{pmatrix}
					y_{71} & y_{72}\\
					\end{pmatrix}$}}}(F);	
		\end{tikzpicture}
	\end{center}
	\caption[General form of a point $w\in U_{1,2}\subset\modulie$ for $\Gr(4,2)$]{General form of a point $w\in U_{1,2}\subset\modulie$ for $\Gr(4,2)$.}
	\label{fig:gr4-2mats_general}
\end{figure}

\textbf{Aim:} Working through the quiver step by step, the aim is to show that $w$ is equivalent to the point given in \Cref{fig:gr4-2mats}, i.e.\ the point $y\in\Gr(4,2)$ defined using the data from the matrix
\[
W_{(1,0)}=\begin{pmatrix}
1&0&x_1&x_3\\
0&1&x_2&x_4
\end{pmatrix}
\]
satisfies $f_E(y)$ is equivalent to $w$; hence $f_E$ is surjective for $n=4$ as required.

\smallskip
\textsc{Step 1:} \emph{The maps $\cO_{\Gr(4,2)}\rightarrow\cW\rightarrow\bw^2\cW$.}

In this first part of the quiver we have the matrices
\begin{align}
\label{eqn:first_Fs}
F^{(0,0)}_1=\begin{pmatrix}
1\\
0
\end{pmatrix},
F^{(0,0)}_2=\begin{pmatrix}
0\\
1
\end{pmatrix},
F^{(0,0)}_3=\begin{pmatrix}
x_1\\
x_2
\end{pmatrix},
F^{(0,0)}_4=\begin{pmatrix}
x_3\\
x_4
\end{pmatrix}
\end{align}
and
\[
G^{(1,0)}_1=\begin{pmatrix}
y_1&y_2
\end{pmatrix},
G^{(1,0)}_2=\begin{pmatrix}
y_3&y_4
\end{pmatrix},
G^{(1,0)}_3=\begin{pmatrix}
y_5&y_6
\end{pmatrix},
G^{(1,0)}_4=\begin{pmatrix}
y_7&y_8
\end{pmatrix},
\]
which are subject to the relations $G^{(1,0)}_iF^{(0,0)}_j+G^{(1,0)}_jF^{(0,0)}_i=0$ for $1\leq i,j \leq 4$. By first considering when $i=1=j$ and $i=2=j$, we find that $y_1=0=y_4$. When $i=1,j=2$ we have
\begin{align*}
\begin{pmatrix}
0&y_2
\end{pmatrix}
\begin{pmatrix}
0\\
1
\end{pmatrix}+
\begin{pmatrix}
y_3&0
\end{pmatrix}
\begin{pmatrix}
1\\
0
\end{pmatrix}=0\implies y_3=-y_2.
\end{align*}
Next, setting $i=1,j=3$ yields
\begin{align*}
\begin{pmatrix}
0&y_2
\end{pmatrix}
\begin{pmatrix}
x_1\\
x_2
\end{pmatrix}+
\begin{pmatrix}
y_5&y_6
\end{pmatrix}
\begin{pmatrix}
1\\
0
\end{pmatrix}=0\implies y_5=-y_2x_2,
\end{align*}
and setting $i=2,j=3$ yields
\begin{align*}
\begin{pmatrix}
-y_2&0
\end{pmatrix}
\begin{pmatrix}
x_1\\
x_2
\end{pmatrix}+
\begin{pmatrix}
y_5&y_6
\end{pmatrix}
\begin{pmatrix}
0\\
1
\end{pmatrix}=0\implies y_6=y_2x_1.
\end{align*}
Repeating the above two substitutions with $j=4$ yields $y_7=-y_2x_4$ and $y_8=y_2x_3$. So far we have
\[
W_{(1,1)}=\begin{pmatrix}
0&y_2&-y_2&0&-y_2x_2&y_2x_1&-y_2x_4&y_2x_3
\end{pmatrix}.
\]
Since $W_{(1,1)}$ must be full rank we must have $y_2\neq0$. We now use the $\GL(1)$ action at the vertex $\bw^2\cW$ to multiply $W_{(1,1)}$ by $y_2^{-1}$, giving
\begin{align}
\label{eqn:first_Gs}
G^{(1,0)}_1=\begin{pmatrix}
0&1
\end{pmatrix},
G^{(1,0)}_2=\begin{pmatrix}
-1&0
\end{pmatrix},
G^{(1,0)}_3=\begin{pmatrix}
-x_2&x_1
\end{pmatrix},
G^{(1,0)}_4=\begin{pmatrix}
-x_4&x_3
\end{pmatrix}
\end{align}
as required.

\smallskip
\textsc{Step 2:} \emph{The maps $\cO_{\Gr(4,2)}\rightarrow\cW\rightarrow\Sym^2\cW$.}

Here we have the matrices
\[
F^{(1,0)}_1=\begin{pmatrix}
y_9 & y_{10}\\
y_{11} & y_{12}\\
y_{13} & y_{14}\\
\end{pmatrix}, 
F^{(1,0)}_2=\begin{pmatrix}
y_{15} & y_{16}\\
y_{17} & y_{18}\\
y_{19} & y_{20}\\
\end{pmatrix}, 
F^{(1,0)}_3=\begin{pmatrix}
y_{21} & y_{22}\\
y_{23} & y_{24}\\
y_{25} & y_{26}\\
\end{pmatrix}, 
F^{(1,0)}_4=\begin{pmatrix}
y_{27} & y_{28}\\
y_{29} & y_{30}\\
y_{31} & y_{32}\\
\end{pmatrix}
\]
subject to the relations $F^{(1,0)}_i F^{(0,0)}_j-F^{(1,0)}_j F^{(0,0)}_i=0$. First set $i=1,j=2$. Then
\[
\begin{pmatrix}
y_9 & y_{10}\\
y_{11} & y_{12}\\
y_{13} & y_{14}\\
\end{pmatrix}
\begin{pmatrix}
0\\
1
\end{pmatrix}-
\begin{pmatrix}
y_{15} & y_{16}\\
y_{17} & y_{18}\\
y_{19} & y_{20}\\
\end{pmatrix}
\begin{pmatrix}
1\\
0
\end{pmatrix}=
\begin{pmatrix}
0\\
0\\
0
\end{pmatrix}\implies
\begin{cases}
y_{10}=y_{15} \\ y_{12}=y_{17} \\ y_{14}=y_{19}\\
\end{cases},
\]
so the second column of $F^{(1,0)}_1$ equals the first column of $F^{(1,0)}_2$.
Next fix $j=3$ and in turn substitute $i=1$ then $i=2$, yielding
\[
\begin{pmatrix}
y_9 & y_{10}\\
y_{11} & y_{12}\\
y_{13} & y_{14}\\
\end{pmatrix}
\begin{pmatrix}
x_1\\
x_2
\end{pmatrix}-
\begin{pmatrix}
y_{21} & y_{22}\\
y_{23} & y_{24}\\
y_{25} & y_{26}\\
\end{pmatrix}
\begin{pmatrix}
1\\
0
\end{pmatrix}=
\begin{pmatrix}
0\\
0\\
0
\end{pmatrix}\implies
\begin{cases}
y_{21}=y_{9}x_1+y_{10}x_2 \\ y_{23}=y_{11}x_1+y_{12}x_2 \\ y_{25}=y_{13}x_1+y_{14}x_2\\
\end{cases},
\]
\[
\begin{pmatrix}
y_{10} & y_{16}\\
y_{12} & y_{18}\\
y_{14} & y_{20}\\
\end{pmatrix}
\begin{pmatrix}
x_1\\
x_2
\end{pmatrix}-
\begin{pmatrix}
y_{21} & y_{22}\\
y_{23} & y_{24}\\
y_{25} & y_{26}\\
\end{pmatrix}
\begin{pmatrix}
0\\
1
\end{pmatrix}=
\begin{pmatrix}
0\\
0\\
0
\end{pmatrix}\implies
\begin{cases}
y_{22}=y_{10}x_1+y_{16}x_2 \\ y_{24}=y_{12}x_1+y_{18}x_2 \\ y_{26}=y_{14}x_1+y_{20}x_2\\
\end{cases}.
\]
Similarly, fix $j=4$ while substituting $i=1$ and $i=2$ as above to get
\[\begin{array}{ccc}
y_{27}=y_{9}x_3+y_{10}x_4 &\phantom{aa} & y_{28}=y_{10}x_3+y_{16}x_4 \\ y_{29}=y_{11}x_3+y_{12}x_4 & \phantom{aa} & y_{30}=y_{12}x_3+y_{18}x_4 \\ y_{31}=y_{13}x_3+y_{14}x_4 & \phantom{aa} & y_{32}=y_{14}x_3+y_{20}x_4 \\
\end{array}.\]
Substituting all of the above into $W_{(2,0)}$, we have
{ \fontsize{9.5}{5}
\[
W_{(2,0)}=
\begin{pmatrix}
y_9&y_{10}&y_{10}&y_{16}&y_{9}x_1+y_{10}x_2&y_{10}x_1+y_{16}x_2&y_{9}x_3+y_{10}x_4&y_{10}x_3+y_{16}x_4 \\
y_{11}&y_{12}&y_{12}&y_{18}&y_{11}x_1+y_{12}x_2&y_{12}x_1+y_{18}x_2&y_{11}x_3+y_{12}x_4&y_{12}x_3+y_{18}x_4 \\
y_{13}&y_{14}&y_{14}&y_{20}&y_{13}x_1+y_{14}x_2&y_{14}x_1+y_{20}x_2&y_{13}x_3+y_{14}x_4&y_{14}x_3+y_{20}x_4 \\
\end{pmatrix}.
\]
}
We now show that the minor of $W_{(2,0)}$ given by the first, second and fourth columns must be full rank. Suppose for a contradiction that it is not full rank. Then we may use the group action (specifically, $\GL(3)$ acting at the vertex $\Sym^2 \cW$) to produce a row of zeros in this minor. Suppose this is the top row, i.e.\ $y_9,y_{10},y_{16}$ become zero in the new basis (the argument is similar for the other rows). The effect this has on the rest of $W_{(2,0)}$ is that now the entire top row is zero. This contradicts the condition that $W_{(2,0)}$ must be full rank, thus we conclude that the chosen minor must be full rank. Consequently, we may use the group action to turn this minor into the identity matrix. The matrix with respect to this new basis is
\[
W_{(2,0)}=
\begin{pmatrix}
1&0&0&0&x_1&0&x_3&0 \\
0&1&1&0&x_2&x_1&x_4&x_3 \\
0&0&0&1&0&x_2&0&x_4 \\
\end{pmatrix},
\]
thereby yielding the matrices $F^{(1,0)}_1,\dots,F^{(1,0)}_4$ in \Cref{fig:gr4-2mats} as required.

\smallskip
\textsc{Step 3:} \emph{The central square, including the maps  $\cW\rightarrow\bw^2\cW\rightarrow \bw^2\otimes\cW$ and $\cW\rightarrow\Sym^2\cW\rightarrow \bw^2\otimes\cW$.}

Around the bottom right square of $\TQ$ we have the relations 
\begin{align}
\label{eqn:base_square_relns}
G^{(2,0)}_i F^{(1,0)}_j=2F^{(1,1)}_j G^{(1,0)}_i-F^{(1,1)}_i G^{(1,0)}_j, \ 1\leq i,j\leq 4.
\end{align}
When $i=j$ this simplifies to $G^{(2,0)}_i F^{(1,0)}_i=F^{(1,1)}_i G^{(1,0)}_i$. Additionally, the matrix $W_{(2,1)}$, which is formed by concatenating the eight matrices $F^{(1,1)}_1,\dots,F^{(1,1)}_4$, $G^{(2,0)}_1,\dots,G^{(2,0)}_4$, must be full rank.

\smallskip
\textsc{Step 3A:} \emph{Write $y_{41},\dots,y_{52}$ (the entries of $G^{(2,0)}_1,G^{(2,0)}_2$) in terms of $y_{33},\dots,y_{36}$ (the entries of $F^{(1,1)}_1,F^{(1,1)}_2$).}

First, consider \eqref{eqn:base_square_relns} when $i,j=1$ and $i,j=2$ in turn. We have
\[
\begin{pmatrix}
y_{41} & y_{42} & y_{43}\\
y_{44} & y_{45} & y_{46}\\
\end{pmatrix}
\begin{pmatrix}
1 & 0\\
0 & 1\\
0 & 0\\
\end{pmatrix}=
\begin{pmatrix}
y_{33}\\
y_{34}
\end{pmatrix}
\begin{pmatrix}
0&1
\end{pmatrix}\implies
\begin{pmatrix}
y_{41} & y_{42}\\
y_{44} & y_{45}\\
\end{pmatrix}=
\begin{pmatrix}
0 & y_{33}\\
0 & y_{34}\\
\end{pmatrix},
\]
\[
\begin{pmatrix}
y_{47} & y_{48} & y_{49}\\
y_{50} & y_{51} & y_{52}\\
\end{pmatrix}
\begin{pmatrix}
0 & 0\\
1 & 0\\
0 & 1\\
\end{pmatrix}=
\begin{pmatrix}
y_{35}\\
y_{36}
\end{pmatrix}
\begin{pmatrix}
-1&0
\end{pmatrix}\implies
\begin{pmatrix}
y_{48} & y_{49}\\
y_{51} & y_{52}\\
\end{pmatrix}=
\begin{pmatrix}
-y_{35} & 0\\
-y_{36} & 0\\
\end{pmatrix}.
\]
Next, for $i=1$ and $j=2$ we have $G^{(2,0)}_1F^{(1,0)}_2=2F^{(1,1)}_2G^{(1,0)}_1-F^{(1,1)}_1G^{(1,0)}_2$, giving
\begin{align*}
&\begin{pmatrix}
y_{41} & y_{42} & y_{43}\\
y_{44} & y_{45} & y_{46}\\
\end{pmatrix}
\begin{pmatrix}
0 & 0\\
1 & 0\\
0 & 1\\
\end{pmatrix}=
2\begin{pmatrix}
y_{35}\\y_{36}
\end{pmatrix}
\begin{pmatrix}
0&1
\end{pmatrix}-
\begin{pmatrix}
y_{33}\\y_{34}
\end{pmatrix}
\begin{pmatrix}
-1&0
\end{pmatrix}\\
&\implies
\begin{pmatrix}
y_{42} & y_{43}\\
y_{45} & y_{46}\\
\end{pmatrix}=
\begin{pmatrix}
y_{33} & 2y_{35}\\
y_{34} & 2y_{36}\\
\end{pmatrix},
\end{align*}
and when $i=2$ and $j=1$, we have  $G^{(2,0)}_2F^{(1,0)}_1=2F^{(1,1)}_1G^{(1,0)}_2-F^{(1,1)}_2G^{(1,0)}_1$ so that
\begin{align*}
&\begin{pmatrix}
y_{47} & y_{48} & y_{49}\\
y_{50} & y_{51} & y_{52}\\
\end{pmatrix}
\begin{pmatrix}
1 & 0\\
0 & 1\\
0 & 0\\
\end{pmatrix}=
2\begin{pmatrix}
y_{33}\\y_{34}
\end{pmatrix}
\begin{pmatrix}
-1&0
\end{pmatrix}-
\begin{pmatrix}
y_{35}\\y_{36}
\end{pmatrix}
\begin{pmatrix}
0&1
\end{pmatrix}\\
&\implies
\begin{pmatrix}
y_{47} & y_{48}\\
y_{50} & y_{51}\\
\end{pmatrix}=
\begin{pmatrix}
-2y_{33} & -y_{35}\\
-2y_{34} & -y_{36}\\
\end{pmatrix}.
\end{align*}
Combining all of the above, we are able to write the entries of $G^{(2,0)}_1, G^{(2,0)}_2$ in terms of those in $F^{(1,1)}_1,F^{(1,1)}_2$ as follows:
\[
G^{(2,0)}_1=\begin{pmatrix}
0 & y_{33} & 2y_{35}\\
0 & y_{34} & 2y_{36}\\
\end{pmatrix}, \ 
G^{(2,0)}_2=\begin{pmatrix}
-2y_{33} & -y_{35} & 0\\
-2y_{34} & -y_{36} & 0\\
\end{pmatrix}.
\]

\smallskip
\textsc{Step 3B:} \emph{Write $y_{37},y_{38}$ (the entries of $F^{(1,1)}_3$) in terms of $y_{33},\dots,y_{36},x_1,x_2$.}

Taking $i=3$ and $j=1$ we have
\begin{align*}
&\begin{pmatrix}
0 & y_{33} & 2y_{35}\\
0 & y_{34} & 2y_{36}\\
\end{pmatrix}
\begin{pmatrix}
x_1 & 0\\
x_2 & x_1\\
0 & x_2\\
\end{pmatrix}=
2\begin{pmatrix}
y_{37}\\y_{38}
\end{pmatrix}
\begin{pmatrix}
0&1
\end{pmatrix}-
\begin{pmatrix}
y_{33}\\y_{34}
\end{pmatrix}
\begin{pmatrix}
-x_2&x_1
\end{pmatrix}\\
&\implies
\begin{pmatrix}
y_{33}x_2 & y_{33}x_1+2y_{35}x_2\\
y_{34}x_2 & y_{34}x_1+2y_{36}x_2\\
\end{pmatrix}=
\begin{pmatrix}
y_{33}x_2 & 2y_{37}-y_{33}x_1\\
y_{34}x_2 & 2y_{38}-y_{34}x_1\\
\end{pmatrix}\\
&\implies
\begin{cases}
y_{37}=y_{33}x_1+y_{35}x_2\\
y_{38}=y_{34}x_1+y_{36}x_2
\end{cases}
\end{align*}
and so
\[
F^{(1,1)}_3=
\begin{pmatrix}
y_{33}x_1+y_{35}x_2\\y_{34}x_1+y_{36}x_2
\end{pmatrix}.
\]

\smallskip
\textsc{Step 3C:} \emph{Write $y_{53},\dots,y_{58}$ (the entries of $G^{(2,0)}_3$) in terms of $y_{33},\dots,y_{36},x_1,x_2$.}

Taking $i=1$ and $j=3$ we have
\begin{align*}
&\begin{pmatrix}
y_{53} & y_{54} & y_{55}\\
y_{56} & y_{57} & y_{58}\\
\end{pmatrix}
\begin{pmatrix}
1 & 0\\
0 & 1\\
0 & 0\\
\end{pmatrix}=
2\begin{pmatrix}
y_{33}\\y_{34}
\end{pmatrix}
\begin{pmatrix}
-x_2&x_1
\end{pmatrix}-
\begin{pmatrix}
y_{33}x_1+y_{35}x_2\\y_{34}x_1+y_{36}x_2
\end{pmatrix}
\begin{pmatrix}
0&1
\end{pmatrix}\\
&\implies
\begin{pmatrix}
y_{53} & y_{54}\\
y_{56} & y_{57}\\
\end{pmatrix}=
\begin{pmatrix}
-2y_{33}x_2 & y_{33}x_1-y_{35}x_2\\
-2y_{34}x_2 & y_{34}x_1-y_{36}x_2\\
\end{pmatrix}.
\end{align*}
Now take $i=2$ and $j=3$ to get
\begin{align*}
&\begin{pmatrix}
y_{53} & y_{54} & y_{55}\\
y_{56} & y_{57} & y_{58}\\
\end{pmatrix}
\begin{pmatrix}
0 & 0\\
1 & 0\\
0 & 1\\
\end{pmatrix}=
2\begin{pmatrix}
y_{35}\\y_{36}
\end{pmatrix}
\begin{pmatrix}
-x_2&x_1
\end{pmatrix}-
\begin{pmatrix}
y_{33}x_1+y_{35}x_2\\y_{34}x_1+y_{36}x_2
\end{pmatrix}
\begin{pmatrix}
-1&0
\end{pmatrix}\\
&\implies
\begin{pmatrix}
y_{54} & y_{55}\\
y_{57} & y_{58}\\
\end{pmatrix}=
\begin{pmatrix}
y_{33}x_1-y_{35}x_2 & 2y_{35}x_1\\
y_{34}x_1-y_{36}x_2 & 2y_{36}x_1\\
\end{pmatrix}.
\end{align*}
Combining the above and summarising steps \textsc{3B} and \textsc{3C}, we have
\[
F^{(1,1)}_3=\begin{pmatrix}
y_{33}x_1+y_{35}x_2\\
y_{34}x_1+y_{36}x_2\\
\end{pmatrix}, \ 
G^{(2,0)}_3=\begin{pmatrix}
-2y_{33}x_2&y_{33}x_1-y_{35}x_2 & 2y_{35}x_1 \\
-2y_{34}x_2&y_{34}x_1-y_{36}x_2 & 2y_{36}x_1 \\
\end{pmatrix}.
\]

\smallskip
\textsc{Step 3D:} \emph{Write $y_{39},y_{40},y_{59},\dots,y_{64}$ (the entries of $F^{(1,1)}_4$ and $G^{(2,0)}_4$) in terms of $y_{33},\dots,y_{36},x_3,x_4$.}

This step is identical to \textsc{Steps 3B} and \textsc{3C}, only whenever $i$ or $j$ equals $3$, we instead substitute $4$. This gives
\[
F^{(1,1)}_4=\begin{pmatrix}
y_{33}x_3+y_{35}x_4\\
y_{34}x_3+y_{36}x_4\\
\end{pmatrix}, \ 
G^{(2,0)}_4=\begin{pmatrix}
-2y_{33}x_4&y_{33}x_3-y_{35}x_4 & 2y_{35}x_3 \\
-2y_{34}x_4&y_{34}x_3-y_{36}x_4 & 2y_{36}x_3 \\
\end{pmatrix}.
\]

\smallskip
\textsc{Step 3E:} We complete \textsc{Step 3} by repeating the same argument used to conclude \textsc{Step 2}. The matrix $W_{(2,1)}$, formed by concatenating $F^{(1,1)}_1,\dots,F^{(1,1)}_4$, $G^{(2,0)}_1,\dots,G^{(2,0)}_4$, is a $2\times16$ matrix where, due to \textsc{Steps 3A-3D}, every term on the top row is a multiple of either $y_{33}$ or $y_{35}$ and every term on the bottom row is a multiple of either $y_{34}$ or $y_{36}$. The argument at the end of \textsc{Step 2} now applies: the minor formed by the first two columns of $W_{(2,1)}$ must be full rank, otherwise it is possible to use the group action in such a way that an entire row of $W_{(2,1)}$ becomes zero, which contradicts the condition that $W_{(2,1)}$ must be full rank. Therefore, the group action allows us to change basis such that $(F^{(1,1)}_1 \ F^{(1,1)}_2)$ becomes the identity matrix. This forces  $y_{33},y_{36}\mapsto 1$ and $y_{34},y_{35}\mapsto 0$, and the resulting change to the rest of $W_{(2,1)}$ yields the matrices $F^{(1,1)}_1,\dots,F^{(1,1)}_4,G^{(2,0)}_1,\dots,G^{(2,0)}_4$ in \Cref{fig:gr4-2mats} as required.

\smallskip
\textsc{Step 4:} \emph{The maps $\bw^2\cW\rightarrow \bw^2\otimes\cW\rightarrow (\bw^2\cW)^{\otimes 2}$.}

The final step is identical to \textsc{Step 1} because the matrices for the maps $\bw^2\cW\rightarrow \bw^2\otimes\cW$ coincide with those for $\cO_{\Gr(4,2)}\rightarrow \cW$, and the relations involved are identical. This implies $G^{(2,1)}_i=G^{(1,0)}_i$ for $1\leq i\leq 4$ as required (see \Cref{fig:gr4-2mats}) and ensures that $W_{(2,2)}$ is full rank. It is routine to check that the relations $G^{(2,1)}_i G^{(2,0)}_j-G^{(2,1)}_j G^{(2,0)}_i=0$ hold. This completes the proof in the case $n=4$.

\smallskip
\noindent\textbf{Induction step.} Now let $n>4$ and as before, for some fixed point $w\in\modulie$ fix a basis $u_1,\dots,u_n$ of $V$ such that without loss of generality we may assume $W_{(1,0)}$ takes the form
\[
W_{(1,0)}=
\begin{pmatrix}
1 & 0 & x_1 & x_3 & \cdots & x_{2n-7} & x_{2n-5}\\ 
0 & 1 & x_2 & x_4 & \cdots &x_{2n-6} & x_{2n-4}
\end{pmatrix}.
\]
Recall the tilting quiver $\TQ$ of $\Gr(n,2)$ (see \Cref{fig:grn2tiltQ}) and notice that the tilting quiver for $\Gr(n-1,2)$ is the following sub-quiver $S$ of $\TQ$:
\begin{align}
\label{eqn:quiver_S}
\begin{split}
S_0&:=\left\{\lambda \in \ZZ^2 \mid n-3\geq\lambda_1\geq\lambda_2\geq 0\right\}, \\
S_1&:=\left\{a_\rho^{\lambda,i} \,\middle\vert\, 
\begin{array}{l}
1\leq \rho\leq n-1\\
i\in\{1,2\}, \ \lambda,\lambda+e_i\in S_0\\
\tl(a_\rho^{\lambda,i})=\lambda, \  \hd(a_\rho^{\lambda,i})=\lambda+e_i
\end{array}\right\}.
\end{split}
\end{align}
Observe that $S$ is not a full sub-quiver of $\TQ$ since the arrows corresponding to $u_n$ are missing.

The induction hypothesis is as follows: suppose that all matrices corresponding to the arrows in $S_1$ are of the forms defined in \eqref{eqn:Fmatrix} and \eqref{eqn:Gmatrix}. The base case for this argument is the previous section of calculations proving this for $\Gr(4,2)$.

Now, we begin by making the important observation that the matrices defined in \eqref{eqn:Fmatrix} and \eqref{eqn:Gmatrix} do not actually depend on the precise vertex $\lambda$ but only the value $\lambda_1-\lambda_2+1$, i.e.\ the rank of the bundle $\SS^\lambda\cW$. For example, in \Cref{fig:gr4-2mats} compare the matrices for $\cO_{\Gr(4,2)} \rightarrow \cW$ with $\bw^2 \cW \rightarrow \bw^2 \cW \otimes \cW$ ($2\times1$ matrices), and the matrices for $\cW \rightarrow \bw^2 \cW$ with $\bw^2 \cW \otimes \cW \rightarrow (\bw^2 \cW)^{\otimes 2}$ ($1\times2$ matrices). Hence, the following notation is well-defined for all $1\leq i \leq n-1$:
\begin{align}
\label{eqn:matrix_defns}
\begin{split}
\text{For all} \ 1\leq k\leq n-3, \ F_i^{(k)}:=F^\lambda_i \phantom{a} \text{for any $\lambda\in S_0$ with $\rank(\SS^\lambda \cW)=k$}, \\
\text{For all} \ 2\leq k\leq n-2, \ G_i^{(k)}:=G^\lambda_i \phantom{a} \text{for any $\lambda\in S_0$ with $\rank(\SS^\lambda \cW)=k$}.
\end{split}
\end{align}
\textbf{Aim:} Prove that the remaining matrices corresponding to the arrows in $\TQ_1\setminus S_1$ are also of this form. This implies $w$ is equivalent to a point in the image of $f_E$ and therefore that $f_E$ is an isomorphism; this completes the proof of \Cref{thm:grn2_isomorphism}.

\smallskip
We first deal with arrows in $\TQ_1$ corresponding to the basis vector $u_n$ that have head at a vertex in $S_0$. For these matrices, recall that certain $\lambda$ in the base case, the work done to find $F_4^\lambda$ or $G_4^\lambda$ is identical to the work required to find $F_3^\lambda,G_3^\lambda$. This generalises to all $i>3$, and we may therefore extend our above definition of the matrices $F_i^{(k)}$ and $G_i^{(k)}$ to include $i=n$. It is routine to check that these matrices satisfy all of the relations and stability conditions.

Next we deal with all of the remaining arrows that do not have head or tail at the bottom right corner vertex $(n-2,0)$; these are the vertices on the right-hand column of $\TQ$ except for $(n-2,0)$. All of the vector bundles situated at the head or tail of these arrows have rank less than or equal to $n-2$. Similar to above, the relations and stability conditions on the matrices corresponding to these arrows have already been used elsewhere in the quiver to prove that all matrices of the same size must be of the desired form; see \textsc{Step 4} of the base case, for example. Hence we once again extend the definitions of $F_i^{(k)}$ and $G_i^{(k)}$ in \eqref{eqn:matrix_defns} to all $\lambda\in\TQ_0$ rather than only the proper subset $S_0$.

\begin{figure}[!ht]
	\begin{center}
		\begin{tikzpicture}[xscale=0.9,yscale=0.9]
		\tikzset{edge/.style = {->,> = latex'}}
		\node[thick] (A) at  (0,0) {$\Sym^{n-4}\cW$};
		\node[thick] (B) at  (4,0) {$\Sym^{n-3} \cW$};
		\node[thick] (C) at  (4,4) {\begin{tabular}{c} $\bw^2 \cW$ \\ $\otimes \Sym^{n-4}\cW$ \end{tabular}};
		\node[thick] (D) at  (9,0) {$\Sym^{n-2} \cW$};
		\node[thick] (E) at  (9,4) {\begin{tabular}{c} $\bw^2 \cW$ \\ $\otimes \Sym^{n-3}\cW$ \end{tabular}};
		\node[thick] (F) at  (9,8) {\begin{tabular}{c} $(\bw^2 \cW)^{\otimes 2}$ \\ $\otimes \Sym^{n-4}\cW$ \end{tabular}};
		\node[thick] (G) at  (0,4) {};
		\node[thick] (H) at  (4,8) {};
		\node[thick] (I) at  (-3,0) {};
		\node[thick] (J) at  (9,10) {};

		\draw[edge] (A) -- node[above]{\textcolor{black}{$F_i^{(n-3)}$}} (B);
		
		\draw[edge] (B) -- node[left]{\textcolor{black}{$G_i^{(n-2)}$}} (C);
		
		\draw[edge] (B) -- node[above]{\textcolor{black}{$F_i^{(n-2)}$}} (D);
		
		\draw[edge] (C) -- node[above]{\textcolor{black}{$F_i^{(n-3)}$}} (E);
		
		\draw[edge] (D) -- node[left]{\textcolor{black}{$G_i^{(n-1)}$}} (E);
		
		\draw[edge] (E) -- node[left]{\textcolor{black}{$G_i^{(n-2)}$}} (F);
		
		\draw[dashed, ->] (G) -- (C);
		\draw[dashed, ->] (C) -- (H);
		\draw[dashed, ->] (I) -- (A);
		\draw[dashed, ->] (F) -- (J);
		
		\end{tikzpicture}
	\end{center}
	\caption[Working for \Cref{thm:grn2_isomorphism} induction step]{The lower right corner of the tilting quiver for $\Gr(n,2)$.}
	\label{fig:grn-2_induction}
\end{figure}
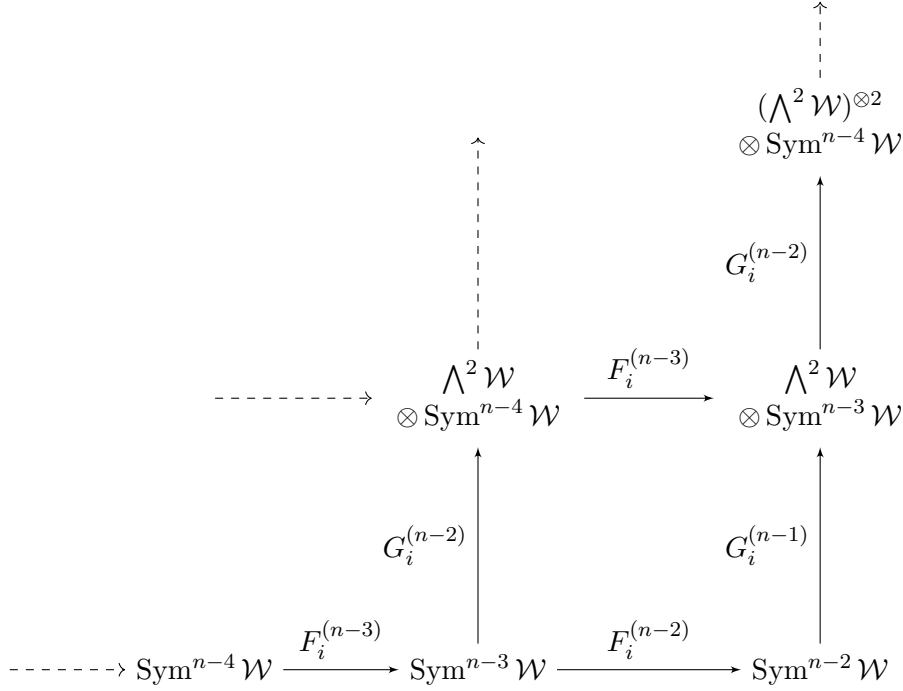

It remains to check the arrows with head or tail at the vertex $(n-2,0)$, i.e.\ the arrows between $\Sym^{n-3} \cW$, $\Sym^{n-2} \cW$ and $\bw^2 \cW \otimes \Sym^{n-3} \cW$ in the lower right corner of the tilting quiver as shown in \Cref{fig:grn-2_induction}. This forms the rest of the proof. For consistency we will denote the matrices corresponding to these arrows by $F^{(n-2)}_i:=F^\lambda_i$ for $\lambda=(n-3,0)$ and $G^{(n-1)}_i:=G^\lambda_i$ for $\lambda=(n-2,0)$. To do this, we will take inspiration from \textsc{Step 2} and \textsc{Step 3} of the base case.

\smallskip
\textsc{Step 1:} \emph{Show that $F_i^{(n-2)}$, $1\leq i \leq n$, takes the form defined in \eqref{eqn:Fmatrix}.}

To find the $F_i^{(n-2)}$ we make use of the relations
\[
F_i^{(n-2)}F_j^{(n-3)}=F_j^{(n-2)}F_i^{(n-3)}, \ 1\leq i,j \leq n.
\]
For $1\leq u \leq n-1$ and $1\leq v \leq n-2$, denote the $(u,v)$-th entry of $F_i^{(n-2)}$ for $i=1,2,3$ as follows:
\[
(F_1^{(n-2)})_{u,v}:=a_{u,v}, \ (F_2^{(n-2)})_{u,v}:=b_{u,v}, \  (F_3^{(n-2)})_{u,v}:=c_{u,v}.
\]
First we study $i=1$ and $j=2$. Since the $F_1^{(n-3)}$ and $F_2^{(n-3)}$ are just $(n-3)\times (n-3)$ identity matrices augmented by a row of zeroes at the bottom and top respectively, the relation implies that $\col_k(F_2^{(n-2)})=\col_{k+1}(F_1^{(n-2)})$ for all $1\leq k\leq n-3$, i.e.\
\begin{align}
\label{eqn:induction_cols}
b_{u,v}=a_{u,v+1}, \ 1\leq v \leq n-3,
\end{align}
thus $F_2^{(n-2)}$ is entirely determined by $F_1^{(n-2)}$ apart from its final column.

Next we set $i=3$ and $j=1$. The relation is $F_3^{(n-2)}F_1^{(n-3)}=F_1^{(n-2)}F_3^{(n-3)}$, and we analyse each side separately. The left hand side is more straightforward as multiplying by $F_1^{(n-3)}$ simply turns the last column of $F_3^{(n-2)}$ into zeroes while leaving the rest unaltered, i.e.\
\[
(F_3^{(n-2)}F_1^{(n-3)})_{u,v}=
\begin{cases}
0 & \text{if} \ v=n-2, \\
c_{u,v} & \text{otherwise}. \\
\end{cases}
\]
For the right hand side recall that $F_3^{(n-3)}$ is the $(n-2)\times (n-3)$ matrix
\[
F_3^{(n-3)}=
\begin{pmatrix}
x_1 & 0 & 0 & \cdots & \cdots & 0\\ 
x_2 & x_1 & 0 & \cdots & \cdots & 0 \\
0 & x_2 & x_1 & \cdots & \cdots & 0 \\
\vdots & \ddots & \ddots & \ddots & \ddots & \vdots \\
0 & 0 & \cdots & \cdots & x_2 & x_1 \\
0 & 0 & \cdots & \cdots & 0 & x_2 \\ 
\end{pmatrix}.
\]
Left-multiplying by $F_1^{(n-2)}$ yields 
\[
(F_1^{(n-2)}F_3^{(n-3)})_{u,v}=
\begin{cases}
0 & \text{if} \ v=n-2, \\
a_{u,v}x_1+a_{u,v+1}x_2 & \text{otherwise}, \\
\end{cases}
\]
and comparing with the above we have
\[
c_{u,v}=a_{u,v}x_1+a_{u,v+1}x_2, \ 1\leq u\leq n-1, \ 1\leq v\leq n-3.
\]
By repeating the above with $i=3$ and $j=2$ we mostly get information about $F_3^{(n-2)}$ that we already have since $F_2^{(n-3)}$ is also an identity matrix with an extra row of zeroes (this time at the top rather than the bottom), and $F_2^{(n-2)}$ is largely determined by $F_1^{(n-2)}$ by \eqref{eqn:induction_cols}. Importantly however, we do pick up the final column of $F_3^{(n-2)}$ in this calculation which is given by
\[
c_{u,n-2}=a_{u,n-2}x_1+b_{u,n-2}x_2.
\]

Combining this with the above, we can now write $F_3^{(n-2)}$ in terms of only the entries of $F_1^{(n-2)}, F_2^{(n-2)}$ and $x_1,x_2$ as follows:
\begin{align}
\label{eqn:Fn-2eqns}
(F_3^{(n-2)})_{u,v}=c_{u,v}=
\begin{cases}
a_{u,v}x_1+a_{u,v+1}x_2 & \text{if} \ 1\leq v \leq n-3, \\
a_{u,n-2}x_1+b_{u,n-2}x_2 & \text{if} \ v=n-2. \\
\end{cases}
\end{align}

We find a similar set of equations for $F_i^{(n-2)}$, $4\leq i\leq n$; simply replace $x_1,x_2$ in \eqref{eqn:Fn-2eqns} with the $i$-th column of $W_{(1,0)}$.

To finish \textsc{Step 1} we make the same observation as at the end of \textsc{Step 2} of the base case. Consider the matrix $W_{(n-2,0)}$, formed by concatenating the $F_i^{(n-2)}$, and suppose for contradiction that the $(n-1)\times(n-1)$ minor formed by taking $F_1^{(n-2)}$ and the final column of $F_2^{(n-2)}$ is not full rank. Then it is possible to use the group action (specifically, $\GL(n-1)$ acting at the vertex $\Sym^{n-2} \cW$) to produce a row of zeros in this minor; suppose for example that this is the top row (the argument for the other rows is similar). The effect this has on the rest of $W_{(n-2,0)}$ is that the entire top row becomes zero. This contradicts the stability condition that $W_{(n-2,0)}$ must be full rank, thus we conclude that the chosen minor must be full rank. As a result, we may use the group action to change the chosen minor into the identity matrix. The resulting change to $W_{(n-2,0)}$ is that each $F_i^{(n-2)}$ takes precisely the required form given in \eqref{eqn:Fmatrix}. This completes \textsc{Step 1}.

\smallskip
\textsc{Step 2:} \emph{Show that $G_i^{(n-1)}$, $1\leq i \leq n$, takes the form defined in \eqref{eqn:Gmatrix}.}

This step is slightly simpler than \textsc{Step 3} of the base case since it remains only to prove that the $G_i^{(n-1)}$ take the required forms. By \Cref{thm:complete_relns_rank2}\four, across the lower right corner square we have the relations
\begin{align}
\label{eqn:ind_corner_relns}
(n-3)G_i^{(n-1)} F_j^{(n-2)}=(n-2)F_j^{(n-3)} G_i^{(n-2)}-F_i^{(n-3)} G_j^{(n-2)}, \ 1\leq i,j\leq n.
\end{align}
We will refresh notation from \textsc{Step 1} and for $1\leq u \leq n-2$, $1\leq v \leq n-1$, denote the $(u,v)$-th entry of $G_i^{(n-1)}$ for $i=1,2,3$ as follows:
\[
(G_1^{(n-1)})_{u,v}:=a_{u,v}, \ (G_2^{(n-1)})_{u,v}:=b_{u,v}, \  (G_3^{(n-1)})_{u,v}:=c_{u,v}.
\]
We first consider the cases when $i,j\in\{1,2\}$. Note that when $i=j$, \eqref{eqn:ind_corner_relns} simplifies to $G_i^{(n-1)} F_i^{(n-2)}=F_i^{(n-3)}G_i^{(n-2)}$. Recall that $F_i^{(n-3)}, F_i^{(n-2)}$ for $i=1,2$ are identity matrices augmented by a row of zeroes at the bottom and top respectively, and that
\[
G_1^{(n-2)}=
\begin{pmatrix}
0 & 1 & 0 & \cdots & 0\\ 
0 & 0 & 2 & \cdots & 0 \\
\vdots & \vdots & \vdots & \ddots & \vdots  \\
0 & 0 & 0 & \cdots & n-3 \\
\end{pmatrix},
\phantom{aaa}
G_2^{(n-2)}=
\begin{pmatrix}
-(n-3) & \cdots & 0 & 0 & 0\\ 
\vdots & \ddots & \vdots & \vdots & \vdots \\
0 & \cdots & -2 & 0 & 0   \\
0 & \cdots & 0 & -1 & 0 \\
\end{pmatrix}.
\]
When $i=1=j$ we have $G_1^{(n-1)} F_1^{(n-2)}=F_1^{(n-3)}G_1^{(n-2)}$. The left hand side is equal to $G_1^{(n-1)}$ with the final column removed, and the right hand side is equal to $G_1^{(n-2)}$ augmented by a row of zeroes at the bottom. Comparing both sides entry-wise yields
\[
G_1^{(n-1)}=
\begin{pmatrix}
0 & 1 & 0 & \cdots & 0 & a_{1,n-1}\\ 
0 & 0 & 2 & \cdots & 0 & a_{2,n-1}\\
\vdots & \vdots & \vdots & \ddots & \vdots & \vdots  \\
0 & 0 & 0 & \cdots & n-3 & a_{n-3,n-1} \\
0 & 0 & 0 & \cdots & 0 & a_{n-2,n-1} \\
\end{pmatrix}.
\]
Repeating for $i=2=j$ yields
\[
G_2^{(n-1)}=
\begin{pmatrix}
b_{1,1} & 0 & \cdots & 0 & 0 & 0\\ 
b_{2,1} &-(n-3) & \cdots & 0 & 0 & 0\\ 
\vdots & \vdots & \ddots & \vdots & \vdots & \vdots \\
b_{n-3,1} & 0 & \cdots & -2 & 0 & 0   \\
b_{n-2,1} & 0 & \cdots & 0 & -1 & 0 \\
\end{pmatrix}.
\]
To find the remaining entries of $G_1^{(n-1)}$ and $G_2^{(n-1)}$ we use \eqref{eqn:ind_corner_relns} with $i=1$ and $j=2$. This reads
\[(n-3)G_1^{(n-1)} F_2^{(n-2)}=(n-2)F_2^{(n-3)} G_1^{(n-2)}-F_1^{(n-3)} G_2^{(n-2)},
\]
which becomes
\begin{align*}
(n-3)\begin{pmatrix}
1 & 0 & \cdots & 0 & a_{1,n-1}\\ 
0 & 2 & \cdots & 0 & a_{2,n-1}\\
\vdots & \vdots & \ddots & \vdots & \vdots  \\
0 & 0 & \cdots & n-3 & a_{n-3,n-1} \\
0 & 0 & \cdots & 0 & a_{n-2,n-1} \\
\end{pmatrix}&= \\
(n-2)\begin{pmatrix}
0 & 0 & 0 & \cdots & 0 \\
0 & 1 & 0 & \cdots & 0\\ 
0 & 0 & 2 & \cdots & 0\\
\vdots & \vdots & \vdots & \ddots & \vdots \\
0 & 0 & 0 & \cdots & n-3 \\
\end{pmatrix}&-
\begin{pmatrix}
-(n-3) & \cdots & 0 & 0 & 0\\ 
\vdots & \ddots & \vdots & \vdots & \vdots \\
0 & \cdots & -2 & 0 & 0   \\
0 & \cdots & 0 & -1 & 0 \\
0 & \cdots & 0 & 0 & 0 \\
\end{pmatrix},
\end{align*}
and so we get $a_{u,n-1}=0$ for $1\leq u \leq n-3$ and $a_{n-2,n-1}=n-2$. We now repeat the above with $i=2$ and $j=1$ to get a similar equation and ultimately discover that the first column of $G_2^{(n-1)}$ satisfies $b_{u,1}=0$ for $2\leq u \leq n-2$ and $b_{1,1}=-(n-2)$. We thus have the required forms for $G_1^{(n-1)}$ and $G_2^{(n-1)}$ as shown below:
\[
G_1^{(n-1)}=
\begin{pmatrix}
0 & 1 & 0 & \cdots & 0\\ 
0 & 0 & 2 & \cdots & 0 \\
\vdots & \vdots & \vdots & \ddots & \vdots  \\
0 & 0 & 0 & \cdots & n-2 \\
\end{pmatrix},
\phantom{aaa}
G_2^{(n-1)}=
\begin{pmatrix}
-(n-2) & \cdots & 0 & 0 & 0\\ 
\vdots & \ddots & \vdots & \vdots & \vdots \\
0 & \cdots & -2 & 0 & 0   \\
0 & \cdots & 0 & -1 & 0 \\
\end{pmatrix}.
\]

It remains to investigate $G_3^{(n-1)}$ (as usual the process of finding $G_i^{(n-1)}$ for $4\leq i \leq n$ will be identical). With $i=3$ and $j=1$, equation \eqref{eqn:ind_corner_relns} becomes
\[
(n-3)G_3^{(n-1)} F_1^{(n-2)}=(n-2)F_1^{(n-3)} G_3^{(n-2)}-F_3^{(n-3)} G_1^{(n-2)}
\]
and so we have
\begin{align*}
&(n-3)\begin{pmatrix}
c_{1,1} & c_{1,2} & \cdots & \cdots & c_{1,n-2}\\ 
c_{2,1} & c_{2,2} & \cdots & \cdots & c_{2,n-2}\\
\vdots & \ddots & \ddots & \ddots & \vdots  \\
\vdots & \ddots & \ddots & \ddots & \vdots \\
c_{n-2,1} & \cdots & \cdots & \cdots & c_{n-2,n-2} \\
\end{pmatrix}= \\
&\phantom{(n-3)}(n-2)\begin{pmatrix}
-(n-3)x_2 & x_1 & 0 & \cdots & \cdots & \cdots & 0\\ 
0 & -(n-4)x_2 & 2x_1 & \cdots & \cdots & \cdots & 0 \\
\vdots & \ddots & \ddots & \ddots & \ddots & \ddots & \vdots  \\
0 & \cdots & \cdots & \cdots & -2x_2 & (n-4)x_1 & 0 \\
0 & \cdots & \cdots & \cdots & 0 & -x_2 & (n-3)x_1 \\
0 & \cdots & \cdots & \cdots & 0 & 0 & 0 \\ 
\end{pmatrix}\\
&\phantom{(n-3)(n2)}-\begin{pmatrix}
0 & x_1 & 0 & 0 & \cdots & \cdots & 0\\ 
0 & x_2 & 2x_1 & 0 & \cdots & \cdots & 0 \\
0 & 0 & 2x_2 & 3x_1 & \cdots & \cdots & 0 \\
\vdots & \ddots & \ddots & \ddots & \ddots & \ddots & \vdots \\
0 & 0 & 0 & \cdots & \cdots & (n-4)x_2 & (n-3)x_1 \\
0 & 0 & 0 & \cdots & \cdots & 0 & (n-3)x_2 \\ 
\end{pmatrix}.
\end{align*}
This gives us most of $G_3^{(n-1)}$. For $1\leq u,v\leq n-2$, we have
\[
c_{u,v}=
\begin{cases}
-(n-2-u+1)x_2 & \text{if} \ u=v,\\
ux_1 & \text{if} \ 1\leq u\leq n-3, v=u+1, \\
0 & \text{otherwise}. \\
\end{cases}
\]
Finally, we must find the last column $c_{u,n-1}$. We repeat the above with $i=3$ and $j=2$ and, similar to previous steps, pick up largely the same set of equations but with the final column included. In particular,
\[
c_{u,n-1}=
\begin{cases}
(n-2)x_1 & \text{if} \ u=n-2, \\
0 & \text{otherwise}. \\
\end{cases}
\]

In conclusion, we have
\[
G_3^{(n-1)}=
\begin{pmatrix}
-(n-2)x_2 & x_1 & 0 & \cdots & \cdots & \cdots & 0\\ 
0 & -(n-3)x_2 & 2x_1 & \cdots & \cdots & \cdots & 0 \\
\vdots & \ddots & \ddots & \ddots & \ddots & \ddots & \vdots  \\
0 & \cdots & \cdots & \cdots & -2x_2 & (n-3)x_1 & 0 \\
0 & \cdots & \cdots & \cdots & 0 & -x_2 & (n-2)x_1 \\ 
\end{pmatrix},
\]
which is precisely the form required by the induction hypothesis. As mentioned above, the proof is identical to show that $G_i^{(n-1)}$ for $4\leq i\leq n$ takes the same form as $G_3^{(n-1)}$ except that $x_1,x_2$ are replaced by the entries of the $i$-th column of $W_{(1,0)}$. Hence, we have shown that each $G_i^{(n-1)}$ is of the form defined in \eqref{eqn:Gmatrix}. This completes \textsc{Step 2} and the proof of \Cref{thm:grn2_isomorphism}. \qed

\bibliographystyle{alpha}
\bibliography{reconstructing_grassmannian_lines}
\end{document}